\documentclass[11pt,a4paper]{article}
\usepackage{jheppub,amsmath,  amssymb,slashed,url,bm,textgreek,upgreek}
\usepackage{graphicx}
\usepackage{epstopdf}
\newtheorem{theorem}{Theorem}[section]
\newtheorem{prop}[theorem]{Proposition}

\newtheorem{definition}[theorem]{Definition}
\def\t{{ \sf t}} 
\def\KW{\bf{KW}}

\def\epsilon{\varepsilon}
\def\tt{{\frak t}}

\def\Tr{{\mathrm{Tr}}}

\def\be{\begin{equation}}
\def\ee{\end{equation}}

\def\hat{\widehat}

\def\tilde{\widetilde}

\def\frak{\mathfrak}
\def\h{\widehat}

\def\RP{{\Bbb{RP}}}

\def\S{{\mathcal S}}

\def\V{{\mathcal V}}
\def\O{{\mathcal O}}

\def\Bbb{\mathbb}

\def\A{{\mathcal A}}
\def\d{{\mathrm d}}

\def\R{{\mathbb R}}
\def\C{{\mathbb C}}

\def\[{\bigl [}

\def\]{\bigr ]}

\def\N{{\mathcal N}}
\def\T{{\mathcal T}}
\def\F{{\mathcal F}}

\def\ad{{\mathrm{ad}}}

\def\L{{\mathcal  L}}
\def\t{\widetilde }
\def\h{\widehat}

\def\V{{\mathcal V}}
\def\I{{\mathcal I}}

\def\M{{\mathcal M}}

\def\W{{\mathcal W}}
\def\P{{\mathcal P}}

\def\ss{{d}}

\def\tilde{\widetilde}

\def\Y{{\eusm Y}}

\font\teneurm=eurm10 \font\seveneurm=eurm7  \font\fiveeurm=eurm5
\newfam\eurmfam
\textfont\eurmfam=\teneurm \scriptfont\eurmfam=\seveneurm
\scriptscriptfont\eurmfam=\fiveeurm

\font\teneusm=eusm10 \font\seveneusm=eusm7 \font\fiveeusm=eusm5
\newfam\eusmfam
\textfont\eusmfam=\teneusm \scriptfont\eusmfam=\seveneusm
\scriptscriptfont\eusmfam=\fiveeusm
\def\eusm#1{{\fam\eusmfam\relax#1}}
\font\tencmmib=cmmib10 \skewchar\tencmmib='177
\font\sevencmmib=cmmib7 \skewchar\sevencmmib='177
\font\fivecmmib=cmmib5 \skewchar\fivecmmib='177
\newfam\cmmibfam
\textfont\cmmibfam=\tencmmib \scriptfont\cmmibfam=\sevencmmib
\scriptscriptfont\cmmibfam=\fivecmmib

\def\del{\partial}
\def\calO{\mathcal O}
\def\LKW{\mathcal L}

\title{The Nahm Pole Boundary Condition}

\author{Rafe Mazzeo$^a$}
\affiliation{$^a$Department of Mathematics, Stanford University, Stanford, CA 94305} 

\author{and Edward Witten$^{b}$}
\affiliation{$^{b}$School of Natural Sciences, Institute for Advanced Study,\\ 1 Einstein Drive, Princeton, NJ 08540 USA}
\abstract{The Nahm pole boundary condition for certain gauge theory equations in four and five dimensions is defined 
by requiring that a solution should have a specified singularity along the boundary.  In
the present paper, we  show that this boundary condition is elliptic and has regularity 
properties analogous to more standard elliptic boundary conditions.  We also
establish a uniqueness theorem for the solution of the relevant equations on a half-space with 
Nahm pole boundary conditions.  These
results are expected to have a generalization involving knots, with
applications to  the Jones polynomial and Khovanov homology.}

\begin{document}\maketitle

\section{Introduction}\label{intro} 

Nahm's equations are a system of ordinary differential equations for three functions $\vec\phi=(\phi_1,\phi_2,\phi_3)$ of a real variable $y$
that take values in the Lie algebra $\frak g$ of a compact Lie group $G$.  These functions  satisfy
\begin{equation}\label{tofu} \frac{\d \phi_1}{\d y}+[\phi_2,\phi_3] =0,\end{equation}
along with cyclic permutations of these equations.  More succinctly, we write
\begin{equation}\label{ofu}\frac{\d\vec\phi}{\d y}+\vec\phi\times \vec\phi = 0 \end{equation}
or
\begin{equation}\label{zofu}\frac{\d \phi_i}{\d y}+\frac{1}{2}\sum_{j,k}\epsilon_{ijk}[\phi_j,\phi_k]=0,\end{equation}
where $\epsilon_{ijk}$ is the antisymmetric tensor with $\epsilon_{123}=1$. These ways of writing the equation show that if we view $\vec\phi$
as an element of $\frak{g}\otimes \R^3$, then Nahm's equation is invariant under the action of $SO(3)$ on $\R^3$.
  
In Nahm's work on magnetic monopole solutions of gauge theory \cite{Nahm}, a key role was played by a special singular solution of Nahm's
equations on the open half-line $y>0$.  The solution reads 
\begin{equation}\label{dofo} \vec\phi(y)=\frac{\vec\tt}{y},\end{equation}
where $\vec\tt=(\tt_1,\tt_2,\tt_3)$ is a triplet of elements of $\frak g$, obeying
\begin{equation}\label{nofo} [\tt_1,\tt_2]=\tt_3, \end{equation}
and cyclic permutations thereof.  In other words, the $\tt_i$ obey the commutation relations of
the Lie algebra $\frak{su}(2)$; we can think of them as the images of a standard basis  of
$\frak{su}(2)$ under a homomorphism\footnote{We are primarily interested in the case that $\varrho$
is non-zero and hence is an embedding of Lie algebras, but our considerations also apply for $\varrho=0$.
See Appendix \ref{groups} for some background and examples concerning homomorphisms from $\frak{su}(2)$ to
a simple Lie algebra  $\frak g$.} $\varrho:\frak{su}(2)\to\frak{g}$.  
We will call this solution the Nahm pole solution.  The Nahm pole
solution has been important in many applications of Nahm's equations; 
for  example, see \cite{K}, which is also relevant as background for the present paper.

Nahm's work on monopoles was embedded in D-brane physics in \cite{Diaconescu}.  
The Nahm pole therefore plays a role in D-brane physics,
and this was explained conceptually in \cite{fuzzy}.   Results about D-branes often have 
implications for  gauge theory, and in the case at hand,
by translating the D-brane results to gauge theory language, one learns \cite{GW} that the 
Nahm pole should be used to define a natural boundary condition not just for Nahm's
1-dimensional equation but for certain gauge theory equations in higher dimensions.  
The equations in question include second order equations of supersymmetric Yang-Mills
theory, and associated first-order equations that are relevant to  the geometric Langlands correspondence
\cite{KW} and the Jones polynomial and Khovanov homology of knots \cite{WittenK,WittenKtwo}.

Our aim in this paper is to elucidate the Nahm pole boundary condition.  Though we 
will also discuss generalizations (including a five-dimensional
equation \cite{haydys,WittenK} that is important in the application to Khovanov homology), we  
will primarily study a certain system of first-order equations in four
dimensions for a pair $A,\phi$.  Here $A$ is a connection on a $G$-bundle $E\to M$, with 
$M$ an oriented Riemannian four-manifold 
with metric $g$, and $\phi$ is a 1-form on $M$ valued in
$\ad(E)$ (the adjoint bundle associated to $E$).  The equations read 
\begin{align}\label{zobo} F-\phi\wedge\phi+\star\,\d_A\phi &=0 \cr  
 \d_A\star \phi & = 0,\end{align}
 where $\star$ is the Hodge star and $\d_A=\d+[A,\cdot]$ is the gauge-covariant 
 extension of the exterior derivative.  Alternatively, in local coordinates $x^1,\dots,x^4$,
 \begin{align} \label{robo} F_{ij}-[\phi_i,\phi_j]+\epsilon_{ij}{}^{kl}D_k\phi_l &= 0 \cr   D_i\phi^i& = 0, \end{align}
 where $D_i=D/D x^i$ is the covariant derivative (defined using the connection 
 $A$ and the Riemannian connection on the tangent bundle of $M$), $\epsilon_{ijkl}$ is the
 Levi-Civita antisymmetric tensor, and indices are raised and lowered using the metric $g$.   (Summation over repeated indices is understood.)
 These equations (or their generalization to $t\not=1$; see eqn.\ (\ref{noxo}) below) have sometimes been called the KW equations and we will use this name for lack
 of another one.  For recent work on these equations, see \cite{Taubes,Taubestwo,GU}.
  
 To explain the relation to the Nahm pole, take $M$ to be the half-space $x^4\geq 0$ in a copy of $\R^4$ with Euclidean coordinates $x^1,\dots,x^4$ (oriented with
 $\epsilon_{1234}=1$). We denote this half-space as $\R^4_+$ and write $\vec x=(x^1,x^2,x^3)$ and $y=x^4$.  
 The KW equations  have a simple exact solution
 \begin{equation}\label{telmo} A=0,~~ \phi=\frac{\sum_{a=1}^3 \tt_a\,\d x^a}{y}, \end{equation}
 where the $\tt_a$ obey the $\frak{su}(2)$ commutation relations (\ref{nofo}).  
 This  gives an embedding of the basic Nahm pole solution (\ref{zobo}) in four-dimensional gauge theory, for any choice of the homomorphism 
 $\varrho:\frak{su}(2)\to \frak g$.  However, in many applications, the basic case  is that $\varrho$ defines a principal embedding of $\frak{su}(2)$ in $\frak g$,
 in the sense of Kostant.  For $G=SU(N)$, 
 this means that the $N$-dimensional representation of $G$ is an irreducible representation of $\varrho(\frak{su}(2))$; in general, the principal embedding is the closest
 analog of this for any $G$.
  
 If $\varrho $ is a principal embedding, we also say that $\varrho$ is regular or that $\phi$ has a regular Nahm pole.  The motivation for this terminology
is that if $\varrho$ is a principal embedding, then any nonzero complex linear combination of the $\tt_a$ is  a regular element  of the complex 
Lie algebra ${\frak g}_\C=\frak g\otimes _\R\C$; for instance, $\tt_1+i\tt_2$ is a regular nilpotent element.
  
For every $\varrho$, one defines \cite{GW} a natural boundary condition on the KW equations  that we call the Nahm pole boundary condition, but in this introduction,
we consider only the case of a principal embedding. (For more detail and the generalization to any $\varrho$, see section \ref{nonregular}.) 
For $M=\R^4_+$ and $\varrho$ a principal embedding, the Nahm pole boundary condition is defined by saying 
that one only allows solutions that coincide with the Nahm pole solution (\ref{telmo}) modulo terms that are less singular for $y\to 0$; the equation then implies 
that in a suitable gauge these less
singular terms actually vanish for $y\to 0$.
The Nahm pole boundary condition can be generalized, with some care, to a more general four-manifold with boundary. 
See section 3.4 of \cite{WittenK} and also section \ref{fourmanifold} of the present paper.

There are two main results of the present paper.  The first is that the Nahm pole boundary condition is elliptic.  Since the equation and solutions
contain singular terms, this is not the standard notion of ellipticity of boundary problems, formulated for example using the Lopatinski-Schapiro 
conditions, but is the analog of this in the framework of uniformly degenerate operators \cite{M-edge}. In fact, we verify the ellipticity of the 
linearization of this problem. The data prescribing the Nahm pole boundary condition are inherently discrete, so the linearization measures the 
fluctuations of the solution relative to this principal term. The boundary conditions for this linear operator simply require solutions to 
blow up less quickly than the Nahm pole; see section 2.4 for a precise statement.  The steps needed to verify that this linearization with
such boundary conditions is elliptic involve first computing the indicial roots of the problem, and then showing that the linear operator
in the model setting of the upper half-space $\R^4_+$ is invertible on a certain space of pairs $(a,\varphi)$ satisfying these boundary conditions. 
The indicial roots measure the formal rates of growth or decay of solutions as $y \to 0$. One of the key consequences of ellipticity
is that the actual solutions of this linearized problem, and eventually also the nonlinear equations, possess asymptotic expansions
with exponents determined by these initial roots. This is a strong regularity statement which allows us to manipulate solutions
to these equations rather freely.

The second main result here is a uniqueness theorem for the KW equations with Nahm pole boundary condition.  This states that a solution of these
equations on $M=\R^4_+$ which satisfies the Nahm pole condition at $y=0$ and which is also asymptotic at a suitable rate to the Nahm pole 
solution for $(\vec x, y)$ large must actually be the Nahm pole solution. This uniqueness theorem is  important  in the application
to the Jones polynomial \cite{WittenK} and corresponds to the expected  result that the Jones polynomial of the empty link is trivial.  
The proof of the uniqueness theorem involves finding a suitable Weitzenbock formula adapted to the Nahm pole solution, and showing that 
the fluctuations around the Nahm pole solution decay at a rate sufficient to justify that boundary terms in the Weitzenbock formula vanish.
Essentially the same reasoning leads to an analogous uniqueness theorem for the related five-dimensional equation that is expected to give a 
description of Khovanov homology.  In this case, the uniqueness theorem corresponds to the statement that the Khovanov homology of the 
empty link is of rank 1.  This Weitzenbock formula can be linearized, and this version of it is used to establish the second part of the proof
that the linearized boundary problem is elliptic. The uniqueness theorem and the ellipticity both hold for arbitrary $\varrho$. 

The uniqueness theorem means roughly that solutions of the KW equations with Nahm pole boundary condition do not exhibit ``bubbling'' along the boundary.
The basis for this statement is that on $\R_+^4$, the KW equations  and 
also their Nahm pole solution  are scale-invariant, that is invariant under $(\vec x,y)\to (\lambda \vec x,\lambda y)$, $\lambda>0$.  If there were a non-trivial
solution on $\R^4_+$ with the appropriate behavior at infinity, it could be ``scaled down''  by taking $\lambda$ very small and 
glued into any given solution that obeys the Nahm pole boundary conditions.  Ths would give a new approximate
solution that obeys the same boundary conditions and coincides with the given solution except in a very small region near the boundary; the behavior
for $\lambda\to 0$ would be somewhat similar to bubbling of a small Yang-Mills instanton.  

The Nahm pole boundary condition can be naturally generalized to include knots.  
 In the framework of \cite{WittenK}, this is done by modifying the boundary conditions in the equations (\ref{zobo})
along a knot or link $K\subset\partial M$.  The appropriate general procedure for this is only known if $\varrho$ is a principal embedding.
The model case is that $M=\R_+^4$ and $K$ is a straight line $\R\subset \R^3=\partial M$.  To every irreducible
representation $R^\vee$ of the Langlands or GNO dual group $G^\vee$ of $G$, one associates a model solution of eqns. (\ref{zobo}) that coincides with the Nahm pole solution
away from $K$ and has a more complicated singular behavior along $K$.  This more complicated behavior  depends on $R^\vee$.  Solutions for the model case
were found in section 3.6 of \cite{WittenK} for 
$G$ of rank 1 and  in \cite{Mikhaylov} for any  $G$.  A boundary condition on eqns. (\ref{zobo}) is then defined by saying that a 
solution should be asymptotic to this model solution along $K$, and to have a Nahm pole singularity elsewhere along $\partial M$. 

This boundary condition can again be extended naturally, with some care, to 
the case that $M$ has a product structure $W\times \R_+$ near its boundary, with an arbitrary embedded knot or link in
$W=\partial M$.  (In the case of a link with several connected components, each component can be labeled by a different representation of $G^\vee$, corresponding to
a different singular model solution.)  The Nahm pole boundary condition in the presence of a knot is again subject to a uniqueness theorem, which
says that for $M=\R^4_+$, with $K=\R\subset\partial M$, and for any representation $R^\vee$, a solution that agrees with the
model solution near $\partial M$ and has appropriate behavior at infinity must actually coincide with the model solution.  This more general type of uniqueness
theorem and the closely related ellipticity of the boundary condition in the presence of a knot will be described elsewhere.

\section{Uniqueness Theorem For The Nahm Pole Solution}\label{second}
In this section we lay out the strategy for proving uniqueness of the Nahm pole solution. The centerpiece of this is the introduction of
the Nahm pole boundary condition, and the analysis which shows that this is an elliptic boundary condition, so that solutions have well
controlled asymptotics near the boundary.  A proper statement of this boundary condition requires a somewhat elaborate calculation of the
indicial roots of the problem. These are the formal growth rates of solutions, but without further analysis, there is no guarantee
that solutions grow at these precise rates.  This further analysis rests on the verification of the ellipticity of the linearized KW operator acting on 
fields with a certain imposed growth rate at the boundary.   The second main result here is a Weitzenbock formula for these 
equations. There are a number of such formulas, in fact, and the subtlety is to choose one which is well adapted to solutions with
Nahm pole singularities.  A linearization of this formula plays an important role in understanding ellipticity of the Nahm
pole boundary condition.   

These results and ideas are somewhat intertwined, and we present them in a way that is perhaps not the most logical from a strictly
mathematical point of view, but which emphasizes the essential points as quickly as possible. Thus we first explain the Weitzenbock formula,
and then proceed to the calculation of indicial roots. We are then in a position to give a precise definition of the Nahm pole boundary conditions.
At this point, we use the Weitzenbock formula to prove the uniqueness theorem. This is only a formal calculation unless we prove
that solutions do have these asymptotic rates. This is established in the remainder of the paper.

\subsection{Solutions On A Four-Manifold Without Boundary}\label{background}

We first review how to characterize the solutions of  the KW equations when formulated on an oriented four-manifold $M$ without boundary.
(See section 3.3 of \cite{WittenK}.) The details are not needed in the rest of the paper. This material is included only to motivate the way we will search for
a uniqueness theorem in the presence of the Nahm pole.

As a preliminary, we give a brief proof of ellipticity of the KW equations.  By definition, a nonlinear partial differential equation is called elliptic if its linearization
is elliptic.  For a gauge-invariant equation, this means that the linearization is elliptic if supplemented with a suitable gauge-fixing condition.  In the case of the KW
equations  linearized around a solution $A_{(0)},\phi_{(0)}$, a suitable gauge-fixing condition is $\d_{A_{(0)}}\star (A-A_{(0)})=0$, or equivalently
\begin{equation}\label{zelob}\sum_i\frac{D}{D x^i} (A-A_{(0)})^i=0,  \end{equation}
where $D/D x^i =\partial_i +[A_{(0)\,i},\cdot]$ is the covariant derivative defined using the connection $A_{(0)}$.  Any other gauge condition
that differs from this one by lower order terms also gives an elliptic gauge-fixing condition; a convenient choice turns out to be
\begin{equation}\label{elob}\sum_i\frac{D}{ Dx^i}(A-A_{(0)}^i)+\sum_i[\phi_{(0)\,i},\phi^i-\phi_{(0)}^i]=0. \end{equation}
 It is convenient to regard the linearized
KW equations as equations for a pair  $\Phi=(A-A_{(0)},\star(\phi-\phi_{(0)}))$ consisting of a 1-form and 3-form on $M$ both
valued in $\ad(E)$.  With this interpretation, the symbol of
the linearized and gauge-fixed KW equations is the same as the symbol of the operator $\d+\d^*$ mapping odd-degree differential forms on $M$ valued in $\ad(E)$
to even-degree forms valued in $\ad(E)$.  This is a standard example of an elliptic operator, so the KW equations are elliptic. 
For future reference, we observe that  the $\d+\d^*$ operator admits two standard and very simple elliptic boundary conditions, which are much 
 more straightforward than the Nahm pole boundary
condition which is our main interest in the present paper.\footnote{As we explain in section \ref{nonregular}, these
boundary conditions
can be viewed as special cases of the Nahm pole boundary condition with $\varrho=0$.} These conditions are respectively
\begin{equation}\label{tively}i^*(\Lambda)=0 \end{equation}
and 
\begin{equation}\label{sively}i^*(\star\Lambda)=0,\end{equation}
where $i:\partial M\to M$ is the inclusion and, for a differential form $\omega$ on $M$, $i^*(\omega)$ is the pullback of $\omega$ to $\partial M$.

Now  set
\begin{equation}\label{telbo}\V_{ij}=F_{ij}-[\phi_i,\phi_j]+\epsilon_{ij}{}^{kl}D_k\phi_l,~~~~\V^0=D_i\phi^i, \end{equation}
so that the KW equations are
\begin{equation}\label{elboxoc} \V_{ij}=\V^0=0. \end{equation}
These equations arise in a  twisted supersymmetric gauge theory in which the bosonic part of the action is
\begin{equation}\label{baction}I=-\int_M\d^4x \sqrt g\Tr\left(\frac{1}{2}F_{ij}F^{ij}+D_i\phi_j D^i\phi^j+R_{ij}\phi^i\phi^j+\frac{1}{2}[\phi_i,\phi_j][\phi^i,\phi^j]\right), \end{equation}
where sums over repeated indices are understood and $R_{ij}$ is the Ricci tensor.  Also $\Tr$ is an invariant, nondegenerate, negative-definite quadratic form
on the Lie algebra $\frak g$ of $G$.  For example, for $G=SU(N)$, $\Tr$ can be the trace in the $N$-dimensional representation; the precise normalization of the
quadratic form will not be important in this paper.   

The simplest way to find a vanishing theorem for the KW equations is to form a Weitzenbock formula.  We take the sum of the squares of the equations and
integrate over $M$.  After some integration by parts, and without assuming the boundary of $M$ to vanish, we find
\begin{equation}\label{zoffbo} -\int_M\d^4x \sqrt g\Tr\left(\frac{1}{2}\V_{ij}\V^{ij}+(\V^0)^2\right)=I+\int_{\partial M}\d^3x\, \epsilon^{abc}\Tr\left(\frac{1}{3}\phi_a[\phi_b,\phi_c]-\phi_a
F_{bc}\right). \end{equation} (We write $i,j,k=1,\dots,4$ for indices tangent to $M$ and $a,b,c=1,\dots,3$ for indices tangent to $\partial M$.) In evaluating the boundary term, we 
assume that near its boundary, $M$ is a product $\partial M\times [0,1)$.  We also assume that, if $n$ is the normal vector to $\partial M$, then $n\, \llcorner \phi=0$ along 
$\partial M$, or equivalently, that the pullback of the 3-form $\star\phi$ to $\partial M$ vanishes:
\begin{equation}\label{llomigo} i^* (\star\phi) =0. \end{equation}
This condition is needed to get a useful form for the boundary contribution in the Weitzenbock formula (or alternatively because of supersymmetric
considerations explained in \cite{GW}), so it will be part of the Nahm pole boundary condition.  However, for the rest of this introductory discussion, we 
assume that $\partial M$ is empty. 

If the KW equations $\V_{ij}=\V^0=0$ hold and $\partial M$ vanishes, it follows from the formula above that $I=0$.  This immediately leads to a vanishing theorem: 
if the Ricci tensor of $M$ is 
non-negative, then each term in (\ref{baction}) must separately vanish. Thus, the curvature $F$ must vanish; $\phi$ must be covariantly constant and its components
must commute, $[\phi_i,\phi_j]=0$; and finally $\phi$ must be annihilated by the Ricci tensor, $R_{ij}\phi^j=0$.

Still on an oriented four-manifold without boundary, the KW equations  are
actually subject to a stronger vanishing theorem than we have just explained, because of a fact that is related to the underlying supersymmetry:
modulo a topological invariant,
the functional $I$ can be written as a sum of squares in multiple ways.   To explain this, we first generalize the KW equations  to depend on a real parameter $t$.
Given a two-form $\Lambda$ on $M$, we write $\Lambda=\Lambda^++\Lambda^-$, where $\Lambda^+$ and $\Lambda^-$ are 
the selfdual and anti-selfdual projections of $\Lambda$.  Then we define
\begin{align}\V^+_{ij}(t)&=(F_{ij}-[\phi_i,\phi_j]+t(D_i\phi_j-D_j\phi_i))^+ \cr \V^-_{ij}(t)&=(F_{ij}-[\phi_i,\phi_j]-t^{-1}(D_i\phi_j-D_j\phi_i))^-\cr
\V^0&= D_i\phi^i.
\end{align}
The equations \begin{equation}\label{noxo}\V^+_{ij}(t)=\V^-_{ij}(t)=\V^0=0\end{equation}
 are a one-parameter family\footnote{One can naturally
think of $t$ as taking values in $\R\cup\infty=\RP^1$.  For $t\to 0$, one should multiply $\V^-(t)$ by $t$, and for $t\to\infty$, one should multiply $\V^+(t)$ by $t^{-1}$.
The proof of ellipticity given above at $t=1$ can be extended to all $t$.  One approach to this uses the formula  (\ref{zelg}) below, supplemented by some special
arguments at $t=0,\infty$.}
 of elliptic differential equations that reduce to (\ref{robo}) for $t=1$.
All considerations of this paper can be extended to generic\footnote{The Nahm pole boundary condition is defined for generic $t$ starting with a model solution in which
the Nahm pole appears in $A$ as well as $\phi$, with $t$-dependent coefficients.}   $t$, but to keep the formulas simple and because this case has the closest relation to
Khovanov homology, we will generally focus on the case $t=1$.  

The generalization of eqn.\ (\ref{zoffbo}) to generic $t$ reads
\begin{align}\label{zelg}  -\int_M\d^4x \sqrt g&\Tr\left(\frac{t^{-1}}{t+t^{-1}}\V^+_{ij}(t)\V^{+\,ij}(t)+\frac{t}{t+t^{-1}}\V_{ij}^-(t)\V^{-\,ij}(t)+(\V^0)^2\right)\cr &
=I +\frac{t-t^{-1}}{4(t+t^{-1})}\int_M \d^4 x \,\epsilon^{ijkl}\Tr F_{ij}F_{kl}. \end{align}
In writing this formula, we have assumed that the boundary of $M$ vanishes. (For a more general formula for $\partial M$ non-empty, see eqn.\ (2.60) of \cite{WittenK}.)
Notably, the expression $I$ that appears on the right hand side on (\ref{zelg}) is the functional defined in eqn.\ (\ref{baction}), independent of $t$.  This immediately
leads to very strong  results about possible solutions.  

Suppose, for example, that we find $A,\phi$ obeying the original KW equations (\ref{robo})
at $t=1$.  Then setting $t=1$ in  (\ref{zelg}), the left hand side vanishes, and 
\begin{equation}\label{elg}\P= \frac{t-t^{-1}}{4(t+t^{-1})}\int_M \d^4 x \,\epsilon^{ijkl}\Tr F_{ij}F_{kl} \end{equation}
certainly also vanishes at $t=1$, so therefore $I=0$.   Now suppose that the integral 
 $\int_M\d^4 x \epsilon^{ijkl}\Tr\,F_{ij}F_{kl}$  -- a multiple of which is the first Pontryagin class  $p_1(E)$ -- is nonzero.  Then we can choose $t\not=1$
 to make $\P<0$, and we get a contradiction: the left hand side of (\ref{zelg}) is non-negative, and the right hand side is negative. 
 Hence any solution of the original equations at $t=1$ is on a bundle $E$ with $p_1(E)=0$.  The same is actually true for a solution of the more general eqn.\ (\ref{noxo}) at any
 value of $t$ other than 0 or $\infty$.    To show this, starting with a solution of (\ref{noxo}) at, say, $t=t_0$,  one observes that unless $p_1(E)=0$,
 one would reach the same contradiction as before by considering eqn.\ (\ref{zelg}) at a value $t=t_1$ at which $\P$ is more negative than it is at $t=t_0$.
Such a $t_1$ always exists for $t_0\not= 0,\infty$ unless $p_1(E)=\P=0$.
 
 Once we know that $\P=0$, it  follows  that the right hand side of (\ref{zelg}) is independent of $t$, and hence vanishes for all $t$ if it vanishes for any $t$. But the left
 hand side of (\ref{zelg}) vanishes if and only if the eqns. (\ref{noxo}) are satisfied.  So if $A,\phi$
 obey the eqns. (\ref{noxo}) at any $t\not=0,\infty$, they satisfy those equations for all $t$.    This leads to a simple description of all the solutions (away from $t=0,\infty$).
 Combine $A,\phi$ to a complex connection $\A=A+i\phi$.  We view $\A$ as a connection on a $G_\C$-bundle $E_\C\to M$; here $G_\C$ is a complex
 simple Lie group that is the complexification of $G$, and $E_\C\to M$ is the $G_\C$ bundle that is obtained by complexifying the $G$-bundle $E\to M$.  We also
 define the curvature of $\A$ as $\F=\d\A+\A\wedge \A$.
 The condition that eqns. (\ref{noxo}) are satisfied for all $t$ is that $\F=0$ and $\d_A\star \phi=0$.  By a well-known result \cite{corlette}, solutions of these equations are in 1-1 correspondence
 with homomorphisms $\uppsi:\pi_1(M)\to G_\C$ that satisfy a certain condition of semistability.  (This condition says roughly that if the holonomies of $\uppsi$
 are triangular, then they are actually block-diagonal.)
 
 The key to getting these simple  results was the fact that (modulo a multiple of $p_1(E)$) the same functional $I$ can be written as a sum of squares in more than
 one way.  This fact is related to the underlying supersymmetry.
 We will look for something similar to find a uniqueness theorem associated to the Nahm pole.

\subsection{A Weitzenbock Formula Adapted To The Nahm Pole}\label{beyond}

Now suppose that $M$ has a non-empty boundary, and consider a solution with a Nahm pole along $\partial M$.  The formulas above do not lead to a useful
conclusion directly because the Nahm pole causes the boundary term in (\ref{zoffbo}) to diverge. 

To make this more precise, let us specialize to the case $M=\R_+^4$.  As in the introduction, introduce coordinates  $\vec x=(x^1,x^2,x^3)$ and $y=x^4$ on $\R^4$,
with $M$ the half-space $y\geq 0$.  The familiar Nahm pole solution is given by
\begin{equation}\label{zerr}A=0,~~~\phi=\sum_{a=1}^3\frac{\tt_a\cdot \d x^a}{y}. \end{equation}
(Indices $i,j,k=1,\dots,4$ will refer to all four coordinates $x^1,\dots,x^4$, and indices $a,b,c=1,\dots,3$ will refer to $x^1,x^2,x^3$ only.)
For this solution, the commutators $[\phi_a,\phi_b]$ and covariant derivatives $D_y\phi$ are all of order $1/y^2$, hence not square-integrable 
near $y=0$ (or as $|(\vec x,y)| \to \infty$), and thus the functional $I$ in (\ref{baction}) diverges.  Accordingly, the boundary terms
in the Weitzenbock formula, which we repeat here for convenience (omitting the factor of $\sqrt{g}$ on the left because $M$ is Euclidean),
\begin{equation}\label{zoffboz} -\int_M\d^4x \Tr\left(\frac{1}{2}\V_{ij}\V^{ij}+(\V^0)^2\right)=I+\int_{\partial M}\d^3x \,\epsilon^{abc}\Tr\left(\frac{1}{3}\phi_a[\phi_b,\phi_c]-\phi_a
F_{bc}\right),\end{equation}
are also divergent.   A standard way to regularize such divergences is to replace $M$ by $M_\epsilon = \{ y > \epsilon, \ |(\vec x,y)| < 1/\epsilon\}$,
carry out the integrations by parts, and discard the terms which diverge as $\epsilon \to 0$. For the purposes of the present exposition, let us focus
only on the portion of the boundary where $y = \epsilon$; arguments are given in section \ref{zobot} to show that the contributions from the other part of the boundary 
are negligible. The bulk and boundary terms on the right hand side of (\ref{zoffboz}) are both of order $1/\epsilon^3$ near this lower boundary,
and since the left side vanishes (for the Nahm pole solution), these various diverging contributions on the right must cancel. However, when such a 
cancellation comes into play, it is very difficult to deduce any positivity of the remaining terms on the right, so this formula is not well-suited
to deduce a vanishing theorem.

It is inevitable that the boundary contribution in (\ref{zoffboz}) is at least nonzero for the Nahm pole solution, since otherwise, we could prove
that $I=0$ for this solution, contradicting the fact that $I$ is a sum of squares of quantities (such as $[\phi_a,\phi_b]$) not all of which vanish for the Nahm pole solution.
Observe that once we know that both $I$ and the boundary term are nonvanishing, scale-invariance implies that they must diverge as $\epsilon\to 0$.

To learn something in the presence of the Nahm pole, we need a different way to write the left hand side of (\ref{zoffboz}) as a sum of squares plus a boundary term,
where the boundary term will vanish for any solution that obeys the Nahm pole boundary condition. This will imply a vanishing theorem for such solutions.  
Of course, for this to be possible, the objects whose squares appear on the right hand side of the new formula must vanish in the Nahm pole solution.

So let us write down a set of quantities that vanish in the Nahm pole solution.  It is convenient to expand $\phi=\sum_{a=1}^3 \phi_a\d x ^a+\phi_y \d y$.  
The Nahm pole solution is characterized by $A=\phi_y=0$ and hence trivially
\begin{equation}\label{zon} F= D_i\phi_y =[\phi_i,\phi_y]=0. \end{equation}
Somewhat less trivially, the Nahm pole solution also satisfies
\begin{equation}\label{ozon} W_a = 0 =D_a\phi_b, \end{equation}
where we define
\begin{equation}\label{cozon} W_a=D_y\phi_a+\frac{1}{2}\epsilon_{abc}[\phi_b,\phi_c]. \end{equation}
Conversely, these equations characterize the Nahm pole solution, in the following sense.  The equations (\ref{zon}) and (\ref{ozon}) imply immediately that in a suitable
gauge $A=0$ and $\phi_a$ and $\phi_y$ are functions of $y$ only.  Moreover, $\partial_y\phi_y=0$ (in the gauge with $A=0$), so if $\phi_y$ is required to vanish at  $y=0$ 
(which will be part of the Nahm pole boundary condition) then it vanishes identically.  Finally, the condition $W_a=0$ means that the functions $\phi_a(y)$ obey 
the original 1-dimensional Nahm equation $\d\phi_a/\d y+(1/2)\epsilon_{abc}[\phi_b,\phi_c]=0$.   

This discussion motivates us to replace the functional $I$ of eqn.\ (\ref{baction}) by a new functional $I'$ which is the sum of squares of  objects which vanish for 
the Nahm pole solution: 
\begin{align}\label{longsum}I'=-\int_{\R^3\times \R_+}\d^4x\,\Tr\left(\frac{1}{2}\sum_{i,j}F_{ij}^2+\sum_{a,b}(D_a\phi_b)^2 +\sum_i(D_i\phi_y)^2+\sum_a[\phi_y,\phi_a]^2
+\sum_a W_a^2\right).\end{align}
The only difference between $I$ and $I'$ is that we have replaced $\sum_a(D_y\phi_a)^2+\frac{1}{2}\sum_{a,b}[\phi_a,\phi_b]^2$ by $\sum_aW_a^2$.
Since
\begin{equation}\label{gosum}\Tr\,\left( \sum_a(D_y\phi_a)^2+\frac{1}{2}\sum_{a,b}[\phi_a,\phi_b]^2\right)=\sum_a\Tr \,W_a^2-\frac{1}{3}\partial_y\epsilon^{abc}\Tr \,\phi_a[\phi_b,\phi_c], \end{equation}
the sole effect of this is to change the boundary term in \eqref{zoffbo}, in fact to cancel the cubic terms in $\phi$ that cause the 
divergence as $y\to 0$ in the Nahm pole solution. Eqn.\ (\ref{zoffbo}) is now replaced by the new identity:
\begin{equation}\label{zoffbox}-\int_{\R^3\times \R_+}\d^4x\, \Tr\left(\frac{1}{2}\V_{ij}\V^{ij}+(\V^0)^2\right)=I'-\left(\int_{y=0}-\int_{y=\infty}\right)\d^3x \,\epsilon^{abc}\Tr\,\phi_a F_{bc}+\Delta,
\end{equation}
where
\begin{equation}\label{mofobox}\Delta=\int_{\R^3\times\R_+}\d^4x\,\frac{\partial}{\partial x^i}\Tr\left(\phi_j D^j\phi^i-\phi^iD_j\phi^j\right).  \end{equation}
For the time being, we do not replace $\Delta$ by a boundary integral. If $M$ is a compact manifold with boundary, then $\Delta$ vanishes for a solution 
that is regular along $\partial M$ and satisfies (\ref{llomigo}), which explains why $\Delta$ does not appear in eqn.\ (\ref{zoffbo}).  However, for $M=\R^4_+$, 
the use of (\ref{llomigo}) in eliminating the boundary contribution is less simple in the presence of the Nahm pole, so the term $\Delta$ cannot be dropped trivially
and will be analyzed later.

Now there is a clear strategy for proving a uniqueness theorem for the  Nahm pole solution.  We must show that  any solution that is
asymptotic to the Nahm pole solution for $y\to 0$ and for $|(\vec x,y)| \to \infty$ approaches the Nahm pole solution quickly enough 
that the boundary terms in eqn.\ (\ref{zoffbox}) (including $\Delta$) vanish. It will then follow that $I'=0$ for any such solution. Since $I'$ is a sum of squares 
of quantities that vanish only for a solution derived from the 1-dimensional Nahm solution, the given solution will coincide with the Nahm pole solution everywhere.

\subsection{The Indicial Equation}\label{smally}
\subsubsection{Overview}\label{overview}

Our next task is to  examine in detail the possible behavior of a solution of the KW equation that is asymptotic to the Nahm pole solution (with some $\varrho$)
as $y\to 0$.  This analysis is necessary before we can properly define the Nahm pole boundary condition, and will also be essential for showing that the boundary 
terms in eqn.\ (\ref{zoffbox}) vanish.

In making this analysis, we need to supplement the KW equation with a gauge condition.  
In the Nahm pole boundary condition, we only allow gauge transformations that are trivial\footnote{At the end of section \ref{nonregular}, we explain that  in the case of a nonregular
Nahm pole, one can define a more general boundary condition in which gauge transformations are not required to be trivial at $y=0$.} 
 at $y=0$, and we are interested in a gauge condition that fixes this gauge invariance.

A gauge transformation that vanishes at $y=0$ can be chosen in a unique fashion to make $A_y=0$, and for understanding the asymptotic behavior of perturbations 
of the Nahm pole solution near $y=0$, this is a natural boundary condition.  However, for other purposes (including proving that the Nahm pole boundary condition 
is well-posed, but also studying the boundary terms at infinity in the Weitzenbock formula), it is necessary to choose an elliptic gauge condition, i.e.\ one which 
augments the KW equations to an elliptic system.  Two examples of elliptic gauge conditions were given in equations  (\ref{zelob}) and (\ref{elob}). The Nahm pole 
solution is $A_{(0)}=0$, $\phi_{(0)}=\tt\cdot \d x/y$, and we consider nearby solutions, which we write as $A=a$, $\phi=\tt\cdot \d x/y+\varphi$, so $a$ and $\varphi$ 
are the fluctuations about the Nahm pole.  The gauge conditions (\ref{zelob}) and (\ref{elob}) are $\partial_i a^i=0$ and
\begin{equation}\label{zurimo}\partial_i a^i+ \frac{1}{y}[\tt_a,\varphi_a]=0,  \end{equation}
respectively. Both of these gauge conditions are elliptic, but we use (\ref{zurimo}) as it simplifies the later analysis considerably.

Technically, we assume that $a$ and $\varphi$ admit asymptotic expansions as $y \to 0$, and consider solutions of the KW equations (with
a gauge condition) such that $a$ and $\varphi$ are less singular than $1/y$ there. 
Writing the putative expansion around the Nahm pole solution as 
\begin{equation}\label{donzo} 
A=y^\lambda a_0(\vec x)+\dots , \qquad  \phi=\frac{\sum_{a=1}^3\tt_a\,\d x^a}{y}+ y^\lambda \varphi_0(\vec x)+\dots, 
\end{equation}
where $a_0(\vec x)$, $\varphi_0(\vec x)$ depend only on $\vec x$, and the ellipses refer to terms that are less singular than $y^\lambda$ for $y\to 0$,
we ask which exponents $\lambda > -1$ are allowed if this expression satisfies the equations formally.  

In greater detail, write the KW equations along with a fixed gauge condition  as 
\begin{equation}
{\KW}(A, \phi) = 0.
\label{nlkw}
\end{equation}
Expanding this about the Nahm pole solution yields
\begin{equation}
{\KW}(a, \phi_{(0)} + \varphi) = \LKW (a,\varphi) + Q(a,\varphi),
\label{texpkw}
\end{equation}
where $\L  $ is the linearization of $\KW$ at $(0, \tt \cdot \d x/y)$ and the remainder term $Q$ vanishes quadratically in a suitable sense.  Assuming that
$a$ and $\varphi$ have expansions as above, and that these expansions may be differentiated, multiplied, etc., we see that the most singular terms in $\L(a,\varphi)$ are of order $y^{\lambda-1}$, while $Q(a,\varphi)$ is no more singular than $y^{2\lambda}$. Since $\lambda>-1$, this is less
singular than $y^{\lambda-1}$. Furthermore, only certain terms in $\L( y^\lambda a_0, y^\lambda \varphi_0)$ are as
singular as $y^{\lambda-1}$. Specifically, the terms which
include a $\del_y$ yield a singular factor $y^{\lambda - 1}$, as do the terms containing a commutator with the unperturbed Nahm pole solution.
On the other hand, terms containing $\del_{x^a}$ are $\calO(y^\lambda)$ and hence may be dropped for these considerations. What remains is 
a linear algebraic equation involving $a_0, \varphi_0$ and the exponent $\lambda$. This is known as the indicial equation for the problem. 

We have been somewhat pedantic about separating the steps of first passing to the linearization and then the indicial operator of this
linearization. The same sets of equations can be obtained by directly inserting the putative expansions for $a$ and $\varphi$ into
the nonlinear equations and retaining only the leading terms. The reason for our emphasis will become clear later. 

At this level, the dependence of these coefficients on $\vec x$ is irrelevant, and because of this, the indicial equation respects the symmetry 
$A_a\to -A_a$, $\phi_y\to -\phi_y$, with $\phi_a$ and $A_y$ left unchanged.  This means that the indicial equation uncouples into a system of equations 
for $\varphi_a, a_y$ and another for $a_a,\varphi_y$. These read
\begin{align}\label{lembo}\lambda a_a+[\tt_a,\varphi_y]-\epsilon_{abc}[\tt_b,a_c]&=0 \cr
   \lambda\varphi_y-[\tt_a,a_a]&=0, \end{align} 
and
\begin{align}\label{wembo}\lambda \varphi_a-[\tt_a,a_y]+\epsilon_{abc}[\tt_b,\varphi_c]&=0 \cr \lambda a_y+[\tt_a,\varphi_a]&=0,
\end{align}
respectively. (Had we used the gauge condition $\partial_i a^i=0$ instead of eqn.\ (\ref{zurimo}), the only difference would be that 
the $[\tt_a,\varphi_a]$ term would be missing in the second line of (\ref{wembo}).)

The determination of the indicial roots of the problem, which are defined as the values of $\lambda$ for which these equations have nontrivial solutions, requires a foray into group theory. 
    
\subsubsection{Some Useful Group Theory}\label{group}

\def\T{{\eusm T}}
\def\S{{\eusm S}}
\def\M{{\eusm M}}
\def\ss{{\frak s}}
\def\ff{{\frak f}}
\def\F{{\eusm F}}
One obvious ingredient in these equations is the $\frak{su}(2)$ subalgebra of $\frak g$ that is generated by the $\tt_a$.  This depends on the choice of
homomorphism $\varrho:\frak{su}(2)\to \frak g$; we call its image $\frak{su}(2)_\tt \subset \frak g$.
From the representation theory of $\frak{su}(2)$, we know that up to isomorphism, $\frak{su}(2)$ has one irreducible complex module of dimension $n$ for
every positive integer $n$.  It is convenient to write $n=2j+1$, where $j$ (which is a non-negative half-integer) is called the spin.
In particular under the action of $\frak{su}(2)_\tt$, the complexification $\frak g_\C=\frak g\otimes_\R\C$ of $\frak g$ decomposes as a direct sum of irreducible 
modules $\frak r_\sigma$, of dimension $n_\sigma=2j_\sigma+1$.\footnote{When the $j_\sigma$ 
are all integers, for example in the case of a principal $\frak{su}(2)$ embedding, this statement is true without having to replace $\frak g$ by its complexification
$\frak g_\C$.  The complexification is needed in case some $j_\sigma$ are half-integers.  It can happen in general that several of the $\frak r_\sigma$'s are isomorphic
and in that case the decomposition of $\frak g_\C$ as a direct sum of irreducible $\frak{su}(2)_\tt$ submodules $\frak r_\sigma$ is not unique.  This  does not affect the following analysis.}  For a principal embedding, the $j_\sigma$ are positive integers of which
precisely one is equal to 1 (the submodule of $\frak g$ of spin $j_\sigma=1$ is precisely $\frak{su}(2)_\tt\subset \frak g$).  For example,  $G=SU(N)$ has rank $N-1$
and the values of the $j_\sigma$ are\footnote{For this and additional group-theoretic background, see
Appendix \ref{groups}.}
 $1,2,3,\dots,N-1$.  At the opposite extreme, if  $\varrho=0$, so that the $\tt_a$ all vanish, then $\frak g$ is the direct sum of trivial 1-dimensional 
 $\frak{su}(2)_\tt$ modules, all of spin 0.
 
The indicial equation does not intertwine the $\frak{su}(2)_\tt$ submodules $\frak r_\sigma\subset \frak g_\C$ since none of the terms in the equation do; 
hence the equation can be restricted to any one of the $\frak r_\sigma$. For example, in \eqref{lembo}, it suffices to consider $\varphi_y$ and all the $a_a$ taking values in the same 
submodule $\frak r_\sigma$.  A general solution of the indicial equation is a sum of $\frak r_\sigma$-valued solutions over all the different $\sigma$. 
This is useful because for solutions taking values in a given $\frak r_\sigma$, the endomorphisms appearing in \eqref{lembo} and \eqref{wembo} 
reduce to diagonal operators, so that the equations then completely decouple. 
The calculation making this explicit occupies the remainder of this subsection. 
 
An important property of the algebra $\frak{su}(2)$ is the existence of a quadratic Casimir operator that commutes with the algebra.  
 In general, given any $\frak{su}(2)$ algebra with a basis $\frak b_a$, $a=1,\dots,3$, obeying the $\frak{su}(2)$ relations
 \begin{equation}\label{tildno}[\frak b_a,\frak b_b]=\epsilon_{abc}\frak b_c, \end{equation}
 we define the Casimir as
 \begin{equation}\label{dinko}\Delta=-\sum_{a=1}^3 \frak b_a^2. \end{equation}  On a module of spin $j$, one has
 \begin{equation}\label{inko}\Delta=j(j+1). \end{equation}
 
 In the case of $\frak{su}(2)_\tt\subset\frak g$, we usually write the action of the generators on $\frak g$ as $w\to [\tt_a,w]$ (rather than $w\to \tt_a(w))$.  So we can write the Casimir $\Delta_\T$ as
 \begin{equation}\label{inco}
 \Delta_\T=-\sum_{a=1}^3[\tt_a,[\tt_a,\cdot]],
\end{equation} 
or more abstractly, 
\begin{equation}\label{pilo}\Delta_\T=-\sum_{a=1}^3\tt_a^2,\end{equation}
as in (\ref{dinko}).
  
  It is often best to think of the triple $\vec a=(a_1,a_2,a_3)$ or similarly the triple $\vec\varphi=(\varphi_1,\varphi_2,\varphi_3)$
  as a single element of $\frak g\otimes N$ where $N\cong \R^3$. 
 Another useful $\frak{su}(2)$ algebra acts on the three-dimensional vector space $N$.  (This is simply inherited from invariance of the original KW equations under
 rotations of $\vec x=(x^1,x^2,x^3)$.) Explicitly, we define $3\times 3$ matrices $\ss_a$, $a=1,\dots,3$
 by
 \begin{equation}\label{donkey} (\ss_a)_{bc}=-\epsilon_{abc}.\end{equation}
 These matrices obey the $\frak{su}(2)$ commutation relations
\begin{equation}\label{monkey} [\ss_a,\ss_b]=\epsilon_{abc}\ss_c, \end{equation} and generate an $\frak{su}(2)$ algebra
that we call $\frak{su}(2)_\ss$.
We define the quadratic Casimir $\Delta_\S=-\sum_{a=1}^3\ss_a^2$ and find that $\Delta_\S=2$.
The value 2 is $j(j+1)$ with $j=1$, and reflects the fact that $N$ is an irreducible $\frak{su}(2)_\ss$ module of spin 1.

Finally, we can define a third $\frak{su}(2)$ algebra that we call $\frak{su}(2)_\ff$, generated by $\ff_a=\tt_a+\ss_a$.  
The importance of $\frak{su}(2)_\ff$ is that, since the Nahm pole solution is invariant under $\frak{su}(2)_\ff$ but not under $\frak{su}(2)_\tt$ or $\frak{su}(2)_\ss$,
it is only $\frak{su}(2)_\ff$ that is a symmetry of the indicial equation.  To be more exact, to make $\frak{su}(2)_\ff$ a symmetry of the indicial equation,
we let $\frak{su}(2)_\ff$ act on $\frak g\otimes N$ as just described, while in acting on $\frak g$ itself, we declare that $\ss_a=0$ and $\ff_a=\tt_a+\ss_a=\tt_a$.
Then interpreting $a_a$ and $\varphi_a$ as elements of $\frak g\otimes N$ and $a_y$ and $\varphi_y$ as elements of $\frak g$, with the $\frak{su}(2)_\ff$
action just described, the indicial equation is invariant under $\frak{su}(2)_\ff$.

To exploit this, it is useful to again define a quadratic Casimir $\Delta_\F=-\sum_{a=1}^3\ff_a^2$. Now we have a very useful formula for the $\frak{su}(2)_\ff$-invariant
operator $\ss\cdot \tt=\sum_a\ss_a\cdot \tt_a$:
\begin{equation}\label{mozzo} \ss\cdot \tt=-\frac{1}{2}\left(\Delta_\F-\Delta_\T-\Delta_\S\right). \end{equation}

To make this formula explicit for the module  $\frak r_\sigma\otimes N$, we need to know how to decompose this module under $\frak{su}(2)_\ff$. The answer is given
by the representation theory of $\frak{su}(2)$.
Provided that $j_\sigma\geq 1$, 
the tensor product $\frak r_\sigma\otimes N$
 decomposes under $\frak{su}(2)_\ff$ as a direct sum of modules $\frak r_{\sigma,\eta}$ of spin  $f_{\sigma,\eta}=j_\sigma+\eta$ where $\eta\in\{1,0,-1\}$.
 For $j_\sigma<1$, the decomposition is the same except that the range of values of $\eta$ is smaller; for $j_\sigma=1/2$, one has only $\eta\in\{1,0\}$, and for $j_\sigma=0$
 one has only $\eta=1$.  
 
 In any event, it follows from (\ref{mozzo}) that in acting on $\frak r_{\sigma,\eta}$, the value of $\ss\cdot\tt$ is 
\begin{equation}\label{mexico}\ss\cdot\tt=\begin{cases}-j_\sigma & {\text{if}}\;\; \eta=1\cr 1 & \text{if}\;\;\eta=0\cr j_\sigma+1 & \text{if}\;\;\eta=-1.\end{cases}\end{equation} 
This result is useful because the object $\ss\cdot\tt$  appears in the indicial equation.  For example, understanding $\vec a=(a_1,a_2,a_3)$ as an element of $\frak g\otimes N$,
so that $\ss\cdot \tt\, (\vec a)$ is also a triple of elements $(\ss\cdot \tt\,(\vec a))_a$, $a=1,2,3$ of $\frak g$,
we have from the definitions
\begin{equation}\label{exico} (\ss\cdot \tt \,(\vec a))_a=\epsilon_{abc} [\tt_b,a_c].\end{equation}
The right hand side appears in (\ref{lembo}), and now we have a convenient way to evaluate it.  Similarly, the analogous object $\epsilon_{abc}  [\tt_b,\varphi_c]$ appears
in (\ref{wembo}).  

When we decompose the $\frak r_\sigma$-valued part of the indicial equation (\ref{lembo}) under the action
of the symmetry group $\frak{su}(2)_\ff$,
modules with spin $j_\sigma\pm 1$ (in other words $\eta=\pm 1$) appear only in the $\frak{su}(2)_\ff$ decomposition of
$a_a$, while
spin $j_\sigma$ (or $\eta=0$) appears both in $a_a$ and in $\varphi_y$.  It follows that the terms in \eqref{lembo} involving 
$\varphi_y$, and likewise, the terms in \eqref{wembo} involving $a_y$, only appear when $\eta=0$. 

\subsubsection{The Indicial Roots}\label{throots}

It is now straightforward to determine the indicial roots.  First we consider the pair $a_a,\varphi_y$ and we restrict to the $\frak r_\sigma$-valued part of the equation.
For $\eta\not=0$, we can set $\varphi_y=0$, as explained at the end of section \ref{group},
and so the equation (\ref{lembo}) reduces to $\lambda a_a=\epsilon_{abc}[\tt_b,a_c]$.  The right hand side was analyzed in eqn.\ (\ref{mexico}) and (\ref{exico}),
and so $\lambda=-j_\sigma$ for $\eta=1$ and $\lambda=j_\sigma+1$ for $\eta=-1$.  For $\eta=0$, we have to work a little harder.   We solve
the second equation in (\ref{lembo}) with\footnote{This solution is not valid if $j_\sigma=0$, because of the factor of $j_\sigma$ in the denominator.
For $j_\sigma=0$, eqns. (\ref{lembo}) and (\ref{wembo}) become trivial, since all commutator terms vanish, and tell us that all modes have $\lambda=0$.
This agrees with the result we find in  eqn.\ (\ref{indroots}) below, except that some modes -- the ones with $\lambda=j_\sigma+1$ -- do not exist
for $j_\sigma=0$.} 
\begin{equation}\label{delf}a_a=-\frac{\lambda}{j_\sigma(j_\sigma+1)}[\tt_a,\varphi_y]\end{equation} 
and then after also using the Jacobi identity and the $\frak{su}(2)$ commutation relations, the first equation in (\ref{lembo}) becomes
\begin{equation}\label{urmo}-\frac{\lambda^2}{j_\sigma(j_\sigma+1)}+1+\frac{\lambda}{j_\sigma(j_\sigma+1)}=0.\end{equation}
So for $\eta=0$,  the possible values of $\lambda$ are $j_\sigma+1$ and $-j_\sigma$.  In sum for $a_a,\varphi_y$, the indicial roots are
\begin{equation}\label{indroots}\lambda=\begin{cases}-j_\sigma & {\text{if}}\;\; \eta=1\cr j_\sigma+1,\,-j _\sigma& \text{if}\;\;\eta=0\cr j_\sigma+1 & \text{if}\;\;\eta=-1.\end{cases}\end{equation} 
These results need correction for $j_\sigma<1$, since some modes are missing.  For $j_\sigma=1/2$, the $\lambda=j_\sigma+1=3/2$ mode 
with $\eta=-1$ should be dropped, and for $j_\sigma=0$, both modes with $\lambda=j_\sigma+1=1$ should be dropped.

Inspection of eqns. (\ref{lembo}) and (\ref{wembo}) shows that the indicial roots for the pair $\varphi_a,a_y$ are obtained
from those for $a_a,\varphi_y$ by just changing the sign of $\lambda$.  
So with no need for additional calculations, the indicial roots for the pair $\varphi_a,a_y$ are as follows:
\begin{equation}\label{indrootstwo}\lambda=\begin{cases}j_\sigma & {\text{if}}\;\; \eta=1\cr j_\sigma,~-j_\sigma-1& \text{if}\;\;\eta=0\cr -j_\sigma-1 & \text{if}\;\;\eta=-1.\end{cases}\end{equation} 
Again, some modes should be omitted for $j_\sigma<1$.

It is notable that all of these modes have $a_y=0$, and therefore make sense in the  gauge $A_y=0$,  
except the $\eta=0$ modes in (\ref{indrootstwo}).  Those particular modes are spurious in the sense that they are pure gauge: they are of the form
$a_y=\partial_y u$, $\varphi_a=[\tt_a/y,u]$ with $u(\vec x,y)=y^{\lambda+1}v(\vec x)$. 
After finding a solution of the KW equations, one can always make a gauge transformation that sets $A_y=0$ and eliminates these modes.
However, this can only be usefully done  after finding a global solution:
to develop a general theory of solutions of the KW equations, which can predict the existence of solutions,
one needs an elliptic gauge condition, and such a gauge condition will always allow
pure gauge modes, such as the ones we have identified.  Unlike the pure gauge modes with $a_y\not=0$, the perturbations we have found with $a_y=0$ have gauge-invariant content; they cannot
 be removed by a gauge transformation that is trivial at $y=0$, since there are no gauge transformations that
are trivial at $y=0$ and preserve the condition $a_y=0$.   So the gauge-invariant content of the possible perturbations of the Nahm pole solution near $y=0$
is precisely contained in the modes in (\ref{indroots}) and those in (\ref{indrootstwo}) with $\eta\not=0$.

One more mode in (\ref{indrootstwo})  has a simple interpretation.  Nahm's equations have the familiar Nahm pole solution $\phi_a=\tt_a/y$,
but since Nahm's equations are invariant under shifting $y$ by a constant, they equally well have a solution $\phi_a=\tt_a/(y-y_0)$ for any constant $y_0$.
Differentiating with respect to $y_0$ and setting $y_0=0$, we find that the linearization of Nahm's equations around the Nahm pole solution can be satisfied
by $\varphi_a=\tt_a/y^2$.  This accounts for the mode in (\ref{indrootstwo}) with $j_\sigma=1$, $\lambda=-2$, and $\eta=-1$.

Each value of $\lambda$ that is indicated in (\ref{indroots}) or (\ref{indrootstwo}) represents a space of fluctuations of dimension $2(j_\sigma+\eta)+1$, 
transforming as an irreducible  $\frak{su}(2)_\tt$ module.  Allowing for these multiplicities, the sum of all indicial roots is 0 for 
$a_a,\varphi_y$ and likewise for $\varphi_a,a_y$.
This is a check on the calculations: the indicial roots are eigenvalues of matrices that appear in (\ref{lembo}) and (\ref{wembo}) and are readily seen to be traceless. 

\subsection{The Nahm Pole Boundary Condition}\label{nonregular}
We are now in a position to give a precise formulation of the Nahm pole boundary condition.  This boundary condition depends on 
 the choice of homomorphism $\varrho:\frak{su}(2)\to \frak g$,   which determines the most singular behavior of the 
solution at $\del M$. As explained below, in favorable cases (for example, when $\varrho$ is a regular embedding), the Nahm pole 
boundary condition simply requires that a solution coincides with the Nahm pole solution modulo less singular terms, but in general 
(when $j_\sigma=0$ appears in the decomposition of $\frak g$) the formulation of the boundary condition involves some further details.

Without exploring any of the questions concerning existence of solutions that obey the Nahm pole boundary condition, the expectation 
is that any such solution has an asymptotic expansion at $\del M$, where the leading term is precisely the Nahm pole singularity, and 
that all fluctuations are well-behaved lower order terms in the expansion with strictly less singular rates of blowup.  As already 
suggested by the discussion in section \ref{overview}, the growth rates of these fluctuation terms are governed by the indicial roots, 
which are themselves determined by the linearization of the KW equations at the Nahm pole solution, as in \eqref{texpkw}.  In fact, 
one might try to construct solutions of ${\KW}(A,\phi) = 0$ by fixing the Nahm pole singularity, then setting the right side of 
\eqref{texpkw} to zero to solve for the fluctuation terms. This would involve inverting the linearized operator $\L  $, and it is thus
important to understand the possible invertibility properties of this operator.  In summary, the actual boundary condition we want 
to discuss is one for the linear operator $\L  $ which requires solutions to blow up at some rate strictly less than $y^{-1}$.  
The indicial root calculation above is what leads us to specify a growth (or decay) rate such that $\L  $ is as close to invertible 
as possible when acting on fields with this growth rate. 

\subsubsection{The Regular Case}\label{simple}
First assume that $\varrho$ is a principal embedding of $\frak{su}(2)$ in $\frak g$, or more generally, that in the decomposition of $\frak g_\C$ under
$\frak{su}(2)_\tt$, the minimum value of $j_\sigma$ is $1$. (The two conditions are equivalent for $G=SU(N)$ but
not in general, as explained in Appendix \ref{groups}.)  In this case there is a simplification stemming from the fact that there are no indicial roots with
$-1<\lambda<1$.  We certainly want to exclude fluctuations around the Nahm pole solution with $\lambda\leq -1$, and allow fluctuations with $\lambda>0$. 

Accordingly, for a principal embedding and more generally whenever $j_\sigma=1$ is the smallest value in the decomposition of $\frak g_\C$, we can state the Nahm pole
boundary condition in either of the following two equivalent ways:
\begin{itemize}\item[(1)] A solution satisfies the Nahm pole boundary condition if in a suitable gauge it has an asymptotic expansion as $y \to 0$
with leading term the Nahm pole solution, and all remaining terms less singular than $1/y$.
\item[(2)] A solution satisfies the Nahm pole boundary condition if in a suitable gauge it has an asymptotic expansion as $y \to 0$ with leading term the Nahm 
pole solution, and with all remaining terms vanishing as $y\to 0$. \end{itemize}
Condition (1) is {\it a priori} weaker than condition (2), but they are equivalent because  under our assumption
on the values of $j_\sigma$, there
are no indicial roots in the range $(-1,0]$.

The full explanation of ellipticity of the Nahm pole boundary condition is in section \ref{analysis}. However,
a preliminary observation that plays an important role 
 is that this boundary condition
allows half of the perturbations near $y=0$:  those with $\eta=1$ in $\varphi_a,a_y$, those with $\eta=-1$ in $a_a,\varphi_y$, and 
half of the $\eta=0$ perturbations in
both sets of fields.  (One of the two pure gauge modes that appear at $\eta=0$ in eqn.\ (\ref{indrootstwo}) vanishes
at $y=0$ and one diverges, so the pure gauge modes
did not affect this counting.)

\subsubsection{The General Case}\label{gencase}
As explained in Appendix \ref{groups}, if any value $j_\sigma<1$  occurs in the decomposition of $\frak g_\C$
under $\frak{su}(2)$, then in fact $j_\sigma=0$ occurs in this decomposition.  (If $j_\sigma=0$ occurs
in the decomposition, then $j_\sigma=1/2$ may or may not occur.)
So 
 let us consider the case that $j_\sigma=0$ does occur in the decomposition.  This simply means that there is a nonzero subspace $\frak c$ of $\frak g$ that commutes
with the $\tt_a$; $\frak c$ is automatically a Lie subalgebra of $\frak g$, the Lie algebra of a subgroup $C\subset G$.   We write $\varphi^{\frak c}$, $a^{\frak c}$ for the $\frak c$-valued
parts of $\varphi$ and $a$.  The formulas (\ref{indroots}), (\ref{indrootstwo}) for the indicial roots show that all
modes of $j_\sigma=0$ have $\lambda=0$. Actually, this is clear without a detailed calculation; for $j_\sigma=0$, the commutator terms can be dropped
in (\ref{lembo}) and (\ref{wembo}), which just say that $\lambda=0$. 
This also shows that for $j_\sigma=0$, 
the equations for $a_a$ and $\varphi_y$ decouple and the modes with $\eta=1$ or $0$ describe fluctuations of only $a_a$ or
$\varphi_y$, respectively; a similar remark applies, of course, for $\varphi_a$ and $a_y$.  Exactly what one means by the Nahm pole boundary conditions depends on
how one treats these $\lambda=0$ modes; this will be discussed momentarily.

When both $j_\sigma=1/2$ and $j_\sigma=0$ occur in the decomposition of $\frak g_\C$, there is a further subtlety.
In this case, to preserve the counting mentioned at the end of section \ref{simple},
we have to define the Nahm pole boundary condition to not allow the perturbation with an indicial root $\lambda=-1/2$ that appears in (\ref{indroots})
for $j_\sigma=1/2$ and $\eta=0$.  In other words, we have to require that a solution departs from the Nahm pole solution by
a correction that is less singular than $1/y^{1/2}$.

As for the modes with $j_\sigma=0$, 
there are different physically motivated choices  of how to treat them  \cite{GW}, and these correspond to different boundary conditions.
The general possibility 
is explained at the end of this subsection.     However, for every $\varrho$, there is  a natural boundary condition that we call the strict Nahm pole
boundary condition; for $G=SU(N)$, this is the half-BPS boundary condition that can be naturally realized via D-branes.  For this boundary condition,
we want to leave $\phi_a^{\frak c}$ unconstrained at $y=0$ but to constrain $\phi_y^{\frak c}$, $A_a^{\frak c}$ to vanish at $y=0$.

Thus, we can state the strict Nahm pole
boundary condition for general $\varrho$ in either of the following two equivalent ways:
\begin{itemize}\item[(1)] A solution satisfies the strict Nahm pole boundary condition if in a suitable gauge it has an asymptotic expansion as $y \to 0$ 
with leading term the Nahm pole solution and with remaining terms less singular than $1/y^{1/2}$ (or $1/y$  if $j_\sigma=1/2$ does not occur in the decomposition of $\frak g_\C$), 
with the further restriction that $\phi_y^{\frak c}$ and $A_a^{\frak c}$ vanish at $y=0$.
\item[(2)] A solution satisfies the strict Nahm pole boundary condition if in a suitable gauge it has an asymptotic expansion with leading term the Nahm pole solution 
and with remaining terms vanishing as $y\to 0$, except that $\phi_a^{\frak c}$ is regular at $y=0$ but does not necessarily vanish there.  \end{itemize}

To illustrate the strict Nahm pole boundary condition when $\frak c\not=0$, let us consider the extreme case $\varrho=0$, so that there is no Nahm pole and $\frak c
=\frak g$.  In this case, the strict Nahm pole boundary condition just means that $\phi_a$ is regular at $y=0$ while $\phi_y$ and $A_a$ are constrained to vanish.  This is
an elementary elliptic boundary condition on the KW equation, already formulated in eqn.\ (\ref{tively}) above.  For general $\varrho$, the strict Nahm pole boundary condition
is a sort of hybrid of this case with the opposite case of a principal embedding.  Of course, the well-posedness of the strict Nahm pole boundary
condition is elementary  at $\varrho=0$, while understanding it for $\varrho\not=0$ is the main goal of the present
paper.

Finally, we describe a generalized Nahm pole boundary condition  associated to 
a more general treatment of the $j_\sigma=\lambda=0$ modes.  For this, we pick an arbitrary subalgebra $\frak h\subset \frak c$, corresponding to a subgroup $H\subset C\subset G$, and we denote
as $\frak h^\perp$ the orthocomplement of $\frak h$ in $\frak c$ ($\frak h^\perp$ is a linear subspace of $\frak c$ but generically not a subalgebra).
Then, in addition to allowing only perturbations that are less singular than $1/y^{1/2}$,
 we declare that the $\frak h^\perp$-valued parts of $A_a$ and $\phi_y$ vanish at $y=0$, while the $\frak h$-valued part of these fields is unconstrained; and reciprocally,
we place no constraint on the $\frak h^\perp$-valued part of $\phi_a$, but require the $\frak h$-valued part of $\phi_a$ to vanish at $y=0$. 
  The strict Nahm pole
boundary condition is the case that $\frak h=0$.  To get the generalized Nahm pole boundary condition for 
$\frak h\not=0$, we relax the requirement that a gauge transformation should be trivial at $y=0$; instead, we allow gauge transformations that are $H$-valued at
$y=0$.  
  (For $G=SU(N)$, any $\varrho$, and some specific choices of $\frak h$,
this boundary condition can be realized via a combination of D-branes with an NS5-brane \cite{GW}.)
   For a simple illustration of this more general boundary condition, take $\varrho=0$ and
$\frak h=\frak c=\frak g$. Then the boundary condition is simply that $i^* \phi=0$.
 This is actually a second elementary elliptic boundary condition on the KW equations, which was formulated 
in eqn.\ (\ref{sively}).  We will show in section \ref{genebc} that the boundary condition described in this paragraph is
well-posed for all $\varrho$ and $\frak h$.

\subsection{The Boundary Terms And The Vanishing Theorem}\label{largey}
We can now easily show the vanishing of  the boundary terms at $y=0$ in the Weitzenbock-like formula (\ref{zoffbox}). 
(In doing this, we can ignore the spurious modes that can be removed by a gauge transformation -- the $\eta=0$ modes 
in (\ref{indrootstwo}). Since the boundary terms are gauge-invariant, they do not receive  a contribution from the spurious modes.)   

For example, let us first look at the boundary contribution $\int\d^3 x \,\epsilon_{abc} \Tr\, F_{ab}\phi_c$.  The dominant part of $\phi_c$ 
is the Nahm pole term $\tt_c/y$. To avoid a contribution involving this term, we need the  part of $F_{ab}$ that is valued in 
$\frak{su}(2)_\tt$ (and thus not orthogonal to the coefficient of this Nahm pole)  to vanish faster than $y$ for $y\to 0$.  
The relevant part of $F_{ab}$ is of order $y^2$ for $y\to 0$,  since for example the part of $F$ that is $\frak{su}(2)_\tt$-valued and  
linear in the connection $A$ is a $j_\sigma=1$ mode with an indicial root $j_\sigma+1=2$.  The part of $F$ that is quadratic in $A$ 
vanishes equally fast. (Note that since all of this appears in a boundary integral, we discard the component $F_{ay}$, which only vanishes like $y$.)
 So  the contribution to $\epsilon_{abc}\Tr\,F_{ab}\phi_c$ involving the Nahm pole in $\phi_c$ is of order $y^2\cdot y^{-1}=y$. 
We can also consider contributions to $\epsilon_{abc}\Tr \,F_{ab}\phi_c$ that come from fluctuations in $\phi_c$ around the Nahm 
pole solution (as well as fluctuations in the connection $A_a$).  As long as we do not consider modes with $j_\sigma=0$, all fluctuations 
in either $\phi_c$ or $A_a$ are controlled by strictly positive indicial roots, so the fluctuations vanish at $y=0$ and their contribution 
to $\epsilon_{abc}\Tr \,F_{ab}\phi_c$ vanishes.  The last case to consider is the case of $j_\sigma=0$ fluctuations in both $\phi_c$ 
and $A_a$.  The general boundary condition formulated in the last paragraph of section \ref{nonregular} ensures that the contributions 
of these fluctuations to $\epsilon_{abc}\Tr\,F_{ab}\phi_c$ vanishes, because $F_{ab}$ when restricted to $y=0$ is valued in a subalgebra 
$\frak h\subset\frak c$, while $\phi_c$ is valued in an orthogonal subspace $\frak h^\perp\subset \frak c$.  Thus the general construction 
with arbitrary $\varrho$, $\frak h$ ensures the vanishing of $\epsilon_{abc}\Tr\,F_{ab}\phi_c$ at $y=0$.

The other possible boundary term in the Weitzenbock formula at $y=0$  comes from the expression
 $\Delta$, defined in eqn.\ (\ref{mofobox}).  Here the integral we have to consider at $y=0$ is $\int \d^3x \Tr\, (\phi_a D_a\phi_y-\phi_y D_a\phi_a).$
Again, we first consider the Nahm pole contribution with $\phi_a=\tt_a/y$.  A contribution from this term is avoided for reasons similar to what we found in the
last paragraph.  Indeed, the only contribution to $\Tr\, (\phi_a D_a\phi_y-\phi_y D_a\phi_a)$ that is linear in fluctuations around the Nahm pole solution comes from the $j_\sigma=1$
part of $\phi_y$.  The corresponding indicial root is 2, so the relevant piece of $\phi_y$ vanishes as $y^2$, too quickly to contribute for $y\to 0$ even when multiplied
by the Nahm pole $\tt_a/y$.   Alternatively, we can consider contributions to $\Tr\,(\phi_a D_a\phi_y-\phi_y D_a\phi_a)$ that are bilinear in fluctuations around the
Nahm pole.  Here, for group-theoretic reasons, only modes with $j_\sigma>0$ are relevant.  Modes with $j_\sigma>0$ have indicial roots of at least 
3/2 for $a_a,\varphi_y$
or $1/2$ for $\varphi_a,a_y$, so an expression bilinear in
such  modes and linear in the Nahm pole part of
$\phi_a$ vanishes at least as fast as $y^{3/2}y^{1/2}\cdot y^{-1}\sim y$.  Finally,
we can consider contributions to $\Tr\, (\phi_a D_a\phi_y-\phi_y D_a\phi_a)$ that do not involve the Nahm pole part of $\phi_a$ at all.  Since all relevant indicial
roots are nonnegative, the only possible contribution come from the $j_\sigma=0$ modes for which the indicial root vanishes.  These contributions vanish for
much the same reason as in the last paragraph:  the restriction of $A_a$ and $\phi_y$ to $y=0$ is valued in a subalgebra $\frak h$, while the $j_\sigma=0$ part of
$\phi_a$ is valued in an orthogonal subspace $\frak h^\perp$.

To complete the proof of the uniqueness theorem for the Nahm pole solution,  we need to know that the surface terms in (\ref{zoffbox}) 
also vanish for $\vec x$ and/or $y$ going to $\infty$. Let $r=\sqrt{|\vec x|^2+y^2}$.  To ensure vanishing of the surface terms for 
$\vec x ,y\to\infty$, we need $\epsilon_{abc}\Tr\,F_{ab}\phi_c$ and $\Tr\, (\phi_i D_i\phi_j-\phi_j D_i\phi_i)$ to vanish at infinity 
faster than $1/r^3$.  For example, this is so if the deviation of $A$ and $\phi$ from the Nahm pole
solution $A=0$, $\phi=\tt\cdot \d \vec x/y$ vanishes at infinity faster than $1/r$ and the curvature $F$ and the covariant derivatives 
of $\phi$ vanish faster than $1/r^2$. For a solution of the KW equations with this property, the boundary terms for 
$\vec x,y\to\infty$ vanish, so if such a solution obeys the Nahm pole boundary condition, it actually is the Nahm pole solution.

\subsection{Behavior At Infinity}\label{zobot}

To decide if the uniqueness result stated in section \ref{largey} 
is strong enough to be useful, we need to know if the rate of approach to the Nahm pole solution 
that we had to assume
for $\vec x,y\to\infty$
 is natural.  The goal of the following analysis is to show that it is.  

We start with the same reasoning with which we began the analysis for $y\to 0$.  If a solution of the KW equations does approach
the Nahm pole solution for $r\to\infty$, then its leading deviation from that solution satisfies, at large $r$,
 the linear equation obtained by linearizing  around the Nahm pole solution.  We will show that 
 any solution of that linear equation that vanishes for $r\to\infty$ vanishes at least as fast as $1/r^2$ (and its derivative vanishes at least as
 fast as $1/r^3$).  This holds irrespective of what singularities the solution might have if continued in to small $r$ (where we do not assume
 the linearized KW equations to be valid).  Vanishing of the perturbations as $1/r^2$ for $r\to\infty$
  is more than was needed in section \ref{largey}.  

As before, we will write $a_i$ and $\varphi_i$ for the perturbations of $A_i$ and $\phi_i$ around the Nahm pole solution.
To explain the idea of the analysis, we first describe the simplest case, which is the behavior of $\varphi_y$.  The linearized KW equations
imply a linear equation for $\varphi_y$ independent of all other modes.  This perhaps surprising fact can be proved by taking a certain linear
combination of derivatives of the KW equations.  However, a quicker route is to go back to eqn.\ (\ref{zoffbo}).  At a solution of the KW
equations $\V_{ij}=\V^0=0$, the left hand side of (\ref{zoffbo}) is certainly stationary under variations of $A$ and $\phi$, so the right hand
side is also.  If we consider variations of $A$ and $\phi$ whose support is in the interior of $M$, the boundary term can be ignored
and therefore the functional $I$ is stationary at a solution of the KW equations.  In other words, the KW equations imply the Euler-Lagrange equations
$\delta I/\delta A=\delta I/\delta\phi=0$.  (This statement is part of the relation of the KW equations to
a four-dimensional supersymmetric gauge theory.)

It is straightforward to work out the Euler-Lagrange equation for
$\varphi_y$ from the explicit formula for the action $I$ in eqn.\ (\ref{baction}). A convenient way to do this is to expand the action in powers
of the fluctuations $a$ and $\varphi$, around the Nahm pole solution on $\R^4_+$.  There are no linear terms -- since the Nahm pole solution is a solution -- but there are quadratic terms.
The part of $I$ that is second order in the fluctuations and
has a nontrivial dependence on $\varphi_y$ is
\begin{equation}\label{udz}-\int \d^4x \Tr\left(\sum_i (\partial_i\varphi_y)^2+\sum_a [\phi_a,\varphi_y]^2\right).\end{equation}
Importantly, there are no terms that are bilinear in $\varphi_y$ and the other fluctuations $a$ and $\varphi_a$; that is why at the linearized
level one finds an equation that involves only $\varphi_y$ and not the other fields.  
This equation is just the Euler-Lagrange equation that arises in varying the functional (\ref{udz})
with respect to $\varphi_y$.  We write
\begin{equation}\label{tenso}\Delta=-\sum_{i=1}^4\partial^2_{x^i} \end{equation}
for the Laplacian on $\R_+^4$ (with the gauge connection $A$ taken to vanish, as in
the Nahm pole solution), and of course we set $\phi_a=\tt_a/y$.  The equation for $\varphi_y$ is
\begin{equation}\label{benso}\left(\Delta -\frac{1}{y^2}\sum_a[\tt_a,[\tt_a,\cdot]]\right)\varphi_y=0.   \end{equation}
In four dimensions, 
\begin{equation}\label{lenso}\Delta=-\frac{\partial^2}{\partial r^2}-\frac{3}{r}\frac{\partial}{\partial r}+\frac{\Delta_{S^3}}{r^2},\end{equation}
where $\Delta_{S^3}$ is the Laplacian on the three-sphere $r=1$.  Since we are working on the half-space $\R_+^4$ rather
than all of $\R^4$, we consider $\Delta_{S^3}$ as an operator defined on a hemisphere in $S^3$.  It is convenient
to introduce the polar angle $\psi$ where $y/r=\cos \psi$, so that $\psi=0$ along the positive $y$-axis where $\vec x=0$, and the hemisphere
is defined by $\psi\leq \pi/2$.  One has
\begin{equation}\label{dsthree}\Delta_{S^3}=-\frac{1}{\sin^2\psi}\partial_\psi \sin^2\psi \partial_\psi +\frac{1}{\sin^2\psi}\Delta_{S^2},\end{equation}
where $\Delta_{S^2}$ is the Laplacian on a unit two-sphere.
The boundary condition at $\psi=\pi/2$ is determined by the four-dimensional boundary condition
at $y=0$ and hence
can be read off from section \ref{nonregular}.  In our context, this generally means that $\varphi_y$ vanishes at $\psi=\pi/2$,
except possibly if $j_\sigma=0$, in which case $\varphi_y^{\frak h}$ is not required to vanish at $y=0$ or $\psi=\pi/2$.  (Rather,
the KW equation $D_y\phi_y+D_a\phi_a=0$, where $\phi_a^{\frak h}=0$ at $y=0$, and $A$ is $\frak h$-valued at $y=0$,
implies that $D_y\phi_y^{\frak h}=0$
at $y=0$, so for $\varphi_y^{\frak h}$, eqn.\ (\ref{benso}) should be supplemented with Neumann boundary conditions
at $y=0$.)  Actually, as we will see in a moment, for $j_\sigma>0$, there is a potential that enforces the vanishing of $\varphi_y$ at $\psi=\pi/2$.

To make (\ref{benso}) more explicit, we   replace $-\sum_a[\tt_a,[\tt_a,\cdot]]$ with $j_\sigma(j_\sigma+1)$,
where $j_\sigma$ is defined as in section \ref{group}.    The equation
for $\varphi_y$ becomes
\begin{equation}\label{enso}\left(-\frac{\partial^2}{\partial r^2}-\frac{3}{r}\frac{\partial}{\partial r}+\frac{W}{r^2}\right)\varphi_y=0,
\end{equation}
where 
\begin{equation}\label{tensox} W=\Delta_{S^3} + \frac{j_\sigma(j_\sigma+1)}{\cos^2\psi} = -\frac{1}{\sin^2\psi}\partial_\psi \sin^2\psi \partial_\psi +\frac{1}{\sin^2\psi}\Delta_{S^2}+  \frac{j_\sigma(j_\sigma+1)}{\cos^2\psi} . \end{equation}
$W$ is a self-adjoint operator on the hemisphere with a discrete and non-negative spectrum 
 (strictly positive except for $j_\sigma=0$
and  $\varphi_y\in\frak h$).  
Any solution of (\ref{enso}) is a linear combination of solutions of the form $\varphi_y=r^s f$, where $f$ is an eigenfunction of $W$,
obeying  $Wf=\gamma f$ for some $\gamma\geq0$, and   $s(s+2)-\gamma=0$ or
\begin{equation}\label{zom} s=-1\pm \sqrt{1+\gamma}.\end{equation}
Actually, the spectrum of the operator $W$ can be found in closed form.  In particular, for $j_\sigma>0$, the ground state (which is the unique
everywhere positive
eigenfunction) is $f=\cos^{j_\sigma+1} \psi$,
with eigenvalue $\gamma=(j_\sigma+1)(j_\sigma+3)$.  
So from (\ref{zom}), if  $s$ is negative -- as it must be if $\varphi_y$ is to vanish for $r\to\infty$ --  then $s=-3-j_\sigma$.  The perturbations thus decay for
large $r$ as $r^{-3-j_\sigma}$.   Thus for example
if $\varrho$ is principal so that $j_\sigma\geq 1$ for all modes, then in a solution that is asymptotic to the Nahm pole solution,
$\varphi_y$ vanishes for $r\to\infty$ as $1/r^4$.  

This analysis of the fluctuations is valid at large $r$ even for $y\to 0$, so we can
compare to our study of the indicial roots.  The wavefunction $\varphi_y=\cos^{j_\sigma+1}\psi/r^{j_\sigma+3}$  vanishes for $y\to 0$
as $y^{j_\sigma+1}$, in agreement with (\ref{indroots}), where the positive  indicial roots are $\lambda=j_\sigma+1$.  

As usual, for $j_\sigma=0$, there is more to say as there are actually two types of mode.  For $\varphi_y\in \frak h^\perp$,  $\varphi_y$ obeys Dirichlet boundary
conditions at $y=0$, and all the previous formulas are valid, including the asymptotic behavior $\varphi\sim 1/r^{3+j_\sigma}=1/r^3$.
But for $\varphi_y\in \frak h$, we want Neumann boundary conditions at $\psi=\pi/2$.  The lowest eigenvalue of $W$ is $\gamma=0$,
with eigenfunction 1, leading  to $s=-2$ and $\varphi_y\sim 1/r^2$. 

We have analyzed here a second order equation, not all of whose solutions are necessarily associated to solutions of the linearized KW equation, which is first order.
In practice, we need not explore this issue here in detail since all modes we have found decay more rapidly at infinity than was needed for the vanishing
argument of action \ref{largey}.  The same remarks apply in what follows.

Fluctuations in the other fields can be analyzed along the same lines.  For this, it is convenient to write a general formula for the expansion of the action
$I$ around the Nahm pole solution.  We write $I_2$ for the part of $I$ that is quadratic in the fluctuations $a,\varphi$, and compute that
\begin{equation}\label{monx} I_2=I_{2,0}+I_{2,1}+I_{2,2} \end{equation} with
\begin{align}\label{tonx}  I_{2,0}&=-\int \d^4x\,\Tr\left(\sum_{i,j}\left((\partial_i a_j)^2+(\partial_i\varphi_j )\right)^2+\frac{1}{y^2}\sum_{a,i}\left([\tt_a,a_i]^2+[\tt_a,\varphi_i]^2\right)\right)\cr
 I_{2,1} &=-\int \d^4 x \,\Tr \left(\frac{2}{y^2}\epsilon_{abc} [\tt_a,\varphi_b]\varphi_c]+\frac{4}{y^2}a_y[\tt_a,\varphi_a]
    \right) 
 \cr   I_{2,2}&= ~~\int \d^4x\, \Tr\left(\sum_i\partial_i a^i+\frac{1}{y}[\tt_a,\varphi_a]\right)^2=\int \d^4x\, \Tr\,S^2,  \end{align} 
 where the gauge condition (\ref{zurimo}) was $S=0$.   This illuminates one of the advantages of that gauge condition: because $I_{2,2}$ is
 homogeneous and quadratic in $S$, it is automatically stationary when  $S=0$  and hence
 does not contribute to the Euler-Lagrange equations.  This significantly simplifies the analysis.

Since the spatial part $a_a$ of the connection does not appear in $I_{2,1}$, it obeys an Euler-Lagrange equation that comes entirely from $I_{2,0}$. This
equation coincides with the equation for fluctuations of $\varphi_y$, which we have already analyzed.  This is in keeping with the fact that
$\varphi_y$ and $a_a$ have the same indicial roots, so their fluctuations must have the same behavior for $\psi\to \pi/2$.

The equation for fluctuations of $\varphi_a,a_y$ does receive a contribution from $I_{2,1}$.  This contribution, which only arises for $j_\sigma>0$ (since
$I_{2,1}=0$ for $j_\sigma=0$), slightly modifies
 the behavior of the perturbations for $r\to\infty$.  It can 
 be analyzed using methods similar to those that we used in computing the indicial roots.   
 
 For the same reason as in that analysis, the term in $I_{2,1}$ that involves $a_y$  contributes only for $\eta=0$.  
 For $\eta=\pm 1$, we only have to consider the term involving
$\epsilon_{abc} \Tr\,[\tt_a,\varphi_b]\varphi_c$ which we evaluate  using (\ref{exico}) and (\ref{mexico}), to find that
\begin{equation}\label{pilon}I_{2,1}=-\int \d^4 x \frac{1}{r^2\cos^2\psi}\Tr \,\phi_a \phi_a \cdot \begin{cases} -2j_\sigma & \eta=1\cr
              2(j_\sigma+1) & \eta=-1. \end{cases}\end{equation}
              The equation (\ref{enso}) is modified only by a shift in $W$,
              \begin{equation}\label{dosox}W\to W+\frac{1}{\cos^2\psi}\begin{cases} -2j_\sigma & \eta=1\cr 2(j_\sigma+1) &\eta=-1. \end{cases}\end{equation}
              This is equivalent to replacing $j_\sigma$ by $j_\sigma-\eta$ in the definition (\ref{tensox}) of $W$, so that the perturbations in $\phi_a$ with $\eta=\pm 1$
               vanish
              at infinity as $1/r^{3+j_\sigma-\eta}$.  The modes that decay most slowly              are those with $\eta=1$; they decay for $r\to\infty$
              as  $1/r^{2+j_\sigma}$ and for $\psi\to \pi/2$ as $\cos^{j_\sigma}\psi$.  The last statement is in accord with the value found in (\ref{indrootstwo})
              for the indicial root at $\eta=1$.  

 For $\eta=0$, we can express $\varphi_a$ in terms of a new field $u$ by
 $\varphi_a=[\tt_a,u]/\sqrt{j_\sigma(j_\sigma+1)}$  (recall that
 we can assume that $j_\sigma>0$, since otherwise the perturbation vanishes).  In terms of these variables, the Euler-Lagrange equation turns out to be
 \begin{equation}\label{zongo} \left( -\frac{\partial^2}{\partial r^2}-\frac{3}{r}\frac{\partial}{\partial r}+\frac{W}{r^2}+\frac{2}{r^2\cos^2\psi}M\right)\begin{pmatrix}
 u \cr a_y\end{pmatrix}=0,\end{equation}
 where the contribution of $I_{2,1}$ is the term proportional to 
 \begin{equation}\label{ofty}M =\begin{pmatrix} 1 & \sqrt{j_\sigma(j_\sigma+1) }\cr \sqrt{j_\sigma(j_\sigma+1)}&0\end{pmatrix}.\end{equation}
 The eigenvalues of $M$ are $j_\sigma+1$ and $-j_\sigma$.  Upon substituting one of these eigenvalues for $M$ in (\ref{zongo}), one gets precisely
 the same shifts of $W$ as described in eqn.\ (\ref{dosox}).  So the two modes with $\eta=0$ obey precisely the same equations as the two
 modes with $\eta=\pm 1$.  In particular, the mode that decays most slowly for $r\to\infty$ again decays as $1/r^{2+j_\sigma}$, while vanishing
 as $\cos^{j_\sigma}\psi$ for $\psi\to \pi/2$.  The last statement corresponds to the indicial root
 $\lambda=j_\sigma$ found at $\eta=0$ in (\ref{indrootstwo}).

 Though we motivated this analysis by asking if the conditions needed to get a uniqueness theorem for the Nahm pole solution are reasonable,
the results are applicable more widely.  For example, we may be interested in a solution of the KW equations in which the Nahm pole
boundary condition is modified by inclusion of knots at $y=0$.  As long as the knots are compact, one can look for solutions of the KW equation
that coincide with the Nahm pole solution for $r\to \infty$.  Their rate of approach to that solution will be as we have just described.

We observed in section \ref{background} that the symbol $\sigma$ of the linearized KW equations is the symbol of the operator $D=\d+\d^*$ mapping odd
degree forms to even degree forms.  Let $D^\dagger$ be the adjoint of $D$ and let $\sigma^{\dagger}$ be the adjoint of $\sigma$.
Since $D^\dagger D$ is the Laplacian on odd degree differential forms, it follows that $\sigma^\dagger\sigma$ is the symbol of
the Laplacian.  This is reflected  in the above formulas: the fluctuations are annihilated by a second order differential operator that is
equal to the Laplacian plus corrections of lower order.

\subsection{Extension To Five Dimensions}\label{exfive} 

The four-dimensional KW equations are closely related to a certain system of elliptic differential equations in five dimensions \cite{haydys}, \cite{WittenK,WittenKtwo}
and this relationship is crucial in the application to Khovanov homology.

To explain the relationship, we first specialize the KW equations to a four-manifold of the form $M=W\times \I$, with $W$ an oriented 
three-manifold,  and
$\I$ an oriented one-manifold, possibly with boundary.  We endow  $M$ with a product metric $g_{ab}(x)\d x^a \d x^b+\d y^2$, where the $x^a$, $a=1,\dots,3$, parametrize $W$ and $y$ parametrizes $\I$, and as usual we expand
$\phi=\sum_a \phi_a \d x^a+\phi_y\d y$.  The KW equations have the property
that $\phi_y$ enters only in commutators -- either covariant derivatives $D_a\phi_y=[D_a,\phi_y]$, or commutators $[\phi_a,\phi_y]$ with other components
of $\phi$.  This enables us to do the following.  We replace the four-manifold $M$ by the five-manifold $Y=\R\times M=\R\times W\times \I$, where 
$\R$ is parametrized by a new
coordinate $x^0$.  Then wherever a commutator with $\phi_y$ appears in the KW equations, we simply replace it by a commutator with $D_0=D/Dx^0$.
So we replace $D_a\phi_y$ with $[D_a,D_0]=F_{a0}$, and $[\phi_a,\phi_y]$ with $[\phi_a,D_0]=-D_0\phi_a$.

In this way, we get some partial differential equations in five dimensions.  Most
of what we have said in this paper about the KW equations carries over to them.  For example,
the five-dimensional equations have  a Weitzenbock formula quite analogous to (\ref{zelg}).  
One simply has to replace $[\phi_y,\cdot]$ with $[D_0,\cdot]$ in the formula
for the action functional $I$ of equation (\ref{baction}).  This Weitzenbock formula can be used to prove that the five-dimensional equations
are elliptic for $t\not=0,\infty$.  (See eqn.\ (5.44) of \cite{WittenK} for the Weitzenbock formula at $t=1$.)  The proof of ellipticity 
in the interior amounts to showing that -- similarly to what 
was explained for the KW equations at the end of section \ref{zobot} --  if $\sigma$ is
the symbol of the five-dimensional equations, then $\sigma^\dagger\sigma$ is the symbol of the Laplacian (times the identity matrix) and in particular is invertible; hence $\sigma$
is invertible and the equations are elliptic.  

Though the equations are elliptic for generic $t$, something nice happens precisely for $t=1$ (or $t=-1$, which is equivalent to $t=1$ modulo $\phi\to -\phi$).   Just in this case,
the equations acquire four-dimensional symmetry.  From the way we described these equations, they are formulated on a five-manifold of the particular
form $Y=\R\times W\times \I$.  However, at $t=1$, there is more symmetry, a fact that is essential in the application to Khovanov homology.
  One can replace $\R\times W$ by a general oriented  Riemannian four-manifold $X$
with no additional structure, and formulate
the equations on\footnote{Still more generally, one can formulate these equations -- and they remain elliptic -- 
 on an arbitrary five-manifold $Y$ with an everywhere
non-zero vector field \cite{haydys}.}
 $Y=X\times \I$. We take on $Y$ a product metric
$\sum_{\mu,\nu=0}^3 g_{\mu\nu}\d x^\mu\d x^\nu+\d y^2$, where $x^\mu$, $\mu=0,\dots,3$ are local coordinates on $X$ and $\I$ is parametrized
by $y$.  The equations on $Y$ can be described
as follows.  Let $\Omega^{2,+}\to X$ be the bundle of self-dual two-forms, and using the natural projection $X\times \I\to X$, pull this bundle
back to a bundle over $Y=X\times \I$ that we also denote as $\Omega^{2,+}$.  
The fields appearing in the five-dimensional equations are a connection
$A$ on a $G$-bundle $E\to Y$, and a section $B$ of $\Omega^{2,+}(\ad(E))=\Omega^{2,+}\otimes \ad(E)$.   (For $X=\R\times W$, the relation between $B$ and the object
$\vec\phi=\sum_a\phi_a\d x^a$ that appears in the KW equations is $B_{0a}=\phi_a$, $B_{ab}=\epsilon_{abc}\phi_c$.)  The five-dimensional equations can be written
\begin{equation}\label{moniko}
\begin{split}F^+-  \frac{1}{4}B\times B-\frac{1}{2}D_y B &=0 \\
   F_{y\mu}+ D^\nu B_{\nu\mu}&=0.  \end{split}\end{equation}
Here $F^+$ is the orthogonal projection of the curvature $F$ onto the part valued in $\Omega^{2,+}(\ad(E))$, and $B\times B$ is defined
as follows.  Since $\Omega^{2,+}$ is a rank 3 real bundle with structure group $SO(3)$, there is a natural isomorphism\footnote{For a vector space $V$,
we denote the symmetric and antisymmetric parts of $V\otimes V$ as $\mathrm{Sym}^2V$ and $\wedge^2V$, respectively.}  $\wedge^2\Omega^{2,+}\cong
\Omega^{2,+}$.  By composing this with the Lie bracket $\wedge^2\frak g\to \frak g$, we get a natural map $\mathrm{Sym}^2\Omega^{2,+}(\ad(E))
\to \Omega^{2,+}(\ad(E))$.  The image of $B\otimes B$ under this map is what we call $B\times B$.  An explicit formula, viewing $B$ as a self-dual
two-form valued in $\ad(E)$, is
\begin{equation}\label{omen}(B\times B)_{\mu\nu}=\sum_{\sigma,\tau=0}^3 g^{\sigma\tau}[B_{\mu\sigma},B_{\nu\tau}]. \end{equation}

In view of all this, any solution of the KW equations on $W\times \R_+$ can be viewed as a ``time''-independent solution of the five-dimensional
equations on $\R\times W\times \R_+$ (where we identify $x^0$ as time, as is natural in the application to Khovanov homology), with 
$\phi_y$ reinterpreted as $A_0$.  In particular, the Nahm pole solution on $\R^3\times \R_+$ can be regarded
as a solution of the five-dimensional equations on $\R^4\times \R_+$.    The solution is simply \begin{equation}\label{ford} A=0,~~~B_{0a}=\tt_a/y,~~~
~B_{ab}=\epsilon_{abc}\tt_c/y.\end{equation}
Given the classical Nahm pole solution, we then ask if we can define a boundary condition on the five-dimensional equations by allowing only solutions
that are asymptotic for $y\to 0$ to the Nahm pole.  The first step is to compute the indicial roots.  These are precisely the same for the five-dimensional
equation as for the four-dimensional KW equations, since the indicial roots are defined in terms of solutions that depend only on $y$ and so in particular
are time-independent.  The only real difference between the computation of indicial roots in five dimensions and the four-dimensional computation described in
section \ref{throots} is that the five-dimensional interpretation, with $\phi_y$ reinterpreted as $A_0$, gives a better explanation of the symmetry of the equations
between $a_a$ and $\varphi_y$.  (This symmetry is visible in the formula (\ref{tonx}) for $I_2$, and accounts for the fact that in eqn.\ (\ref{indroots}), the indicial roots
for different values of $\eta$ are pairwise equal.)   Given the indicial roots, one can imitate the discussion in section \ref{nonregular} to define precisely
the Nahm pole boundary condition, and as in section \ref{largey} it follows that the boundary terms at $y=0$ in the Weitzenbock formula vanish.
The Nahm pole solution on $\R^4\times \R_+$ is therefore unique if one requires sufficiently fast convergence as $r=\sqrt{x^2+y^2}$ becomes large.
The analysis in section \ref{zobot} can be repeated to show that the expected rate of convergence at infinity of a solution that does converge to the Nahm
pole solution is fast enough to make the uniqueness result concerning the Nahm pole solution relevant.   Just a few modifications are needed.  In
equation (\ref{benso}), $\Delta$ should now be the five-dimensional Laplacian on a half-space, 
\begin{equation}\label{umongo}\Delta =-\frac{\partial^2}{\partial r^2}-\frac{4}{r}\frac{\partial}{\partial r}+\frac{\Delta_{S^4}}{r^2}. \end{equation}
We expand the fluctuation in the connection around the Nahm pole as $a=\sum_{s=0}^3 a_s\d x^s +a_y \d y$.  The equation obeyed by $a_s$ is now
\begin{equation}\label{tumongo}\left(-\frac{\partial^2}{\partial r^2}-\frac{4}{r}\frac{\partial}{\partial r}+\frac{W}{r^2}\right)a_s=0,~~s=0,\dots,3\end{equation}
where 
\begin{equation}\label{rumongo}W=\Delta_{S^4}+\frac{j_\sigma(j_\sigma+1)}{\cos^2\psi}=-\frac{1}{\sin^3\psi}\partial_\psi \sin^3\psi\partial_\psi+\frac{1}{\sin^2\psi}\Delta_{S^3}
+\frac{j_\sigma(j_\sigma+1)}{\cos^2\psi}. \end{equation}
The lowest eigenvalue of $W$ is now  $(j_\sigma+1)(j_\sigma+4)$, again with eigenfunction $\cos^{j_\sigma+1}\psi$, and now the fluctuations  decay as 
$r^{-4-j_\sigma}$ for $r\to\infty$, with precisely one extra power of $1/r$ compared to the four-dimensional case.
The corresponding formulas for $B$ and $a_y$ are similar to the analysis of $\varphi_a$ and $a_y$ in the four-dimensional case,
 and again, the fluctuations in five dimensions 
decay  with one extra power of $r$ compared to what we found in four dimensions.

\section{The Linearized Operator On A Half-Space}\label{third}
\subsection{Overview}\label{overview2}
A nonlinear partial differential equation is said to be elliptic if its linearization is elliptic.  An important property of a linear elliptic differential operator on
a closed manifold is that its kernel and cokernel are finite-dimensional.  

If a linear elliptic differential equation is considered on a manifold $M$ with a nonempty boundary, then it is necessary to impose some sort of boundary condition 
to make the problem similarly well posed. A boundary condition such that the accompanying problem has a finite dimensional kernel and cokernel, and so that in addition
solutions enjoy optimal regularity properties, is called an elliptic boundary condition. As before, for a nonlinear elliptic differential equation on a manifold $M$ 
with boundary, a choice of boundary condition is called elliptic if it (or its linearization if the boundary condition is also nonlinear) is an elliptic boundary
condition for the linearized operator. For standard nondegenerate elliptic operators, the theory of elliptic boundary conditions is now classical, and the criterion 
for ellipticity of a boundary condition (which involves both the interior and boundary operators) is called the Lopatinski-Schapiro condition. The linearized operator 
in our setting cannot be treated with this classical theory since its coefficients of order $0$ are singular at $y=0$. It is, however, a uniformly degenerate 
elliptic operator, as introduced in \cite{M-edge}. There is a suitable notion of ellipticity for boundary conditions in this setting as well, to which we shall be appealing here. 

As a way to motivate the definition of ellipticity of a given boundary condition, consider the model case where $M$ is the half-space $\R^n_+ = \{x^n\geq 0\}$.
If a linear differential operator on this half-space and a boundary condition on the boundary $\R^{n-1}$ are both invariant under rotations and translations 
in the boundary variables, then the kernel and cokernel of this boundary problem (amongst tempered solutions on $\R^n_+$) are finite-dimensional if and 
only if they are actually trivial. Conversely, a boundary condition with this property, and so that the  accompanying linear operator
as a map between appropriate Sobolev spaces has closed range, is elliptic.  We have stated  this formulation  for operators and boundary conditions 
with substantial symmetry; more generally, if $\L  $ is a general  linear elliptic operator with variable coefficients  and $B$ a possibly non-constant operator giving the boundary conditions, then we may apply 
this condition to the constant coefficient problem on $\R^n_+$ obtained by freezing the coefficients of $\L  $ and $B$ at $q \in \R^{n-1}$. 

These remarks are relevant to the Nahm pole boundary condition on the KW equations.  Our task now is to show that the linearization $\L  $ of the KW operator 
around a solution with Nahm pole boundary data on a four-manifold $M$ with boundary satisfies this ellipticity condition. By the remarks above, 
this is actually equivalent to proving the corresponding property for the linearization of the KW operator around the actual model Nahm pole solution 
on the half-space $\R^4_+$.  We shall explain in section \ref{analysis} how this fits into the analytic theory which justifies the main consequences 
of this paper, namely the regularity at $y = 0$ of more general solutions satisfying Nahm pole boundary conditions and the uniqueness theorem. 

Thus the aim of this section is to study the linearized KW operator on the half-space, and to show that it is an isomorphism on the space of $L^2$
fields which satisfy the Nahm pole boundary conditions.  As suggested above, this consists of two rather separate parts: one involves showing
that the kernel and cokernel of $\L  $ vanish, while the other requires showing that $\L  $ has closed range as a map between the appropriate
function spaces. We undertake  the first of these in the present section. Section \ref{vanishing} contains a proof that the kernel of $\L  $ vanishes. This turns
out to be a rather direct consequence of the formulas that were used in section \ref{second} to establish the uniqueness theorem. As for vanishing of 
the cokernel, a standard strategy, once the kernel is known to vanish, is to show, after reducing to an ODE via a Fourier transform, that the index of $\L  $ (defined as the difference 
in dimension between the kernel and cokernel) vanishes.   We do this in two essentially separate ways. The first involves some fairly elementary
algebraic considerations; see section \ref{index}.
The second is more direct (and much more useful in generalizations).  It  turns out, see section \ref{pseudo}, that the adjoint of $\L  $ is 
conjugate to $-\L  $, a property that we call pseudo skew-adjointness. This immediately gives an isomorphism between the kernel and 
cokernel of $\L  $, so the cokernel vanishes if the kernel does.  Finally, in section \ref{extension}, we show that these arguments carry over more or 
less immediately to the five-dimensional extension of the KW equations that is relevant to Khovanov homology.

The remaining task, to show that the range of $\L  $ is closed, turns out to follow using general machinery that will be explained in section \ref{analysis}.

\subsection{Vanishing Theorem For The Kernel}\label{vanishing}

The uniqueness theorem for the Nahm pole solution
on a half-space was deduced from an identity (\ref{zoffbox}) which reads schematically
\begin{equation}\label{terrfo}-\int \d^4x \,\Tr\,\sum_\lambda\V_\lambda^2 =-\int \d^4x \,\Tr\,\sum_\sigma \W_\sigma^2 ,\end{equation}
where we omit surface terms  since we have shown them to vanish in a solution of the KW equations that is asymptotic to the Nahm pole solution for $y\to 0$ and at infinity. The $\V_\lambda$ are the $\frak g$-valued quantities $\V_{ij}$ and $\V^0$ that appear on the left hand side
of (\ref{zoffbox}).  The $\W_\sigma$ are the objects  $F_{ij}$, $D_a\phi_b$, $D_i\phi_y$, $[\phi_y,\phi_a]$, and $W_a$ whose squares appear on the right hand
side of the definition (\ref{longsum}).    
If the KW equations $\V_\lambda=0$ are obeyed, then the identity (\ref{terrfo}) shows
that  the $\W_\sigma$ vanish, which implies that the solution is
constructed from a solution of Nahm's equations.

Now let us see what this formula tells us about the linearization of the KW equations about the Nahm pole solution.  Schematically, let us combine
the fields $A,\phi$ to an object $\Phi$ (one can think of $\Phi=A+\star\phi$ as an odd-degree differential form on $\R^4_+$ valued in $\mathrm{ad}(E)$).
We write $\Phi_0$ for the Nahm pole solution, and we consider a family of fields $\Phi_s=\Phi_0+s \Phi_1$, where $s$ is a parameter and $\Phi_1$ is a perturbation.  
Expanding the identity (\ref{terrfo}) in powers of $s$, the linear term vanishes because $\V_\lambda = \W_\sigma=0$ at $s=0$.  Taking the second derivative
with respect to $s$, terms such as $\V_\lambda \partial_s^2 \V_\lambda$ vanish at $s=0$ since $\V_\lambda=0$ at $s=0$.  So we get
\begin{equation}\label{terrf}-\int \d^4x \,\Tr\,\sum_\lambda\left(\frac{\partial\V_\lambda}{\partial s}\right)^2 =-\int \d^4x \,\Tr\,\sum_\sigma\left(\frac{\partial \W_\sigma}{\partial s}
\right)^2.\end{equation}
Hence the equations $\partial \V_\lambda/\partial s=0$ are satisfied if and only if the equations $\partial \W_\sigma/\partial s=0$ are satisfied.

The equations $\partial \V_\lambda/\partial s=0$ are the linearization of the KW equations around the Nahm pole solution.  In other words, these equations
are $\L  \Phi_1=0$, where $\L  $ is the linearization of the KW equations and $\Phi_1$ is the perturbation around the Nahm pole solution.

The equations $\partial \W_\sigma/\partial s=0$ imply that $\Phi_1$ actually vanishes if it vanishes at infinity.  For example, the equation $\partial_s F_{ij}=0$
implies that the fluctuation in the connection $A$ can be gauged away, and upon doing so, the equations $\partial_s(D_a\phi_b)=\partial_s(D_a\phi_y)=0$ imply
that the perturbation in $\phi$ is independent of $\vec x$ and so vanishes if it vanishes at infinity.  The other conditions $\partial_s\W_\sigma=0$ imply
that $\Phi_1$ actually comes from a solution of the linearization of Nahm's equation.  Of course, such a perturbation (or any perturbation that is independent of $\vec x$)
is not square-integrable in four dimensions.

Hence if the linearization $\L  $ of Nahm's equations is understood as an operator acting on a Hilbert space of square-integrable wavefunctions, its kernel vanishes.

\subsection{Index}\label{index}

Here we will sketch a standard strategy to prove that the cokernel of $\L  $ vanishes once one knows that the kernel vanishes.  We only provide a sketch
since in the particular case of the KW equations, there is a more powerful and direct  method that we explain in section \ref{pseudo}.  

First of all, using the translation symmetries of $\R^3=\partial(\R^4_+)$, we can look for a momentum eigenstate, that is a 
perturbation of the form $\Phi_1(\vec x,y)=e^{i\vec k \cdot \vec x}F(y)$
where $F$ depends on $y$ only and $\vec k$ is a real ``momentum'' vector.  $F(y)$ is a function on $\R_+$ with values in a finite-dimensional complex vector space $Y$.  Let $d=\dim\,Y$;
for the KW equations, 
$d=8\dim \frak g$.  To show that the cokernel of $\L  $ is trivial for square-integrable wavefunctions, it suffices to show that it vanishes for  momentum eigenstates with
$\vec k\not=0$.

On momentum eigenstates, the linearized KW equation $\LKW \Phi_1=0$ reduces to an equation $\L  _1(\vec k)F(y)=0$,  with
\begin{equation}\label{redy}\L  _1(\vec k)=\frac{d}{d y} + B(y,\vec k), \end{equation}
where $B(y,\vec k)$ is a self-adjoint matrix-valued function of $y$ and $\vec k$. In fact,
\begin{equation}\label{edy} B(y,\vec k)=\frac{B_0}{y}+B_1(\vec k), \end{equation}
where $B_0$ is a constant matrix (independent of $y$ and $\vec k$) and $B_1$ is independent of $y$ and homogeneous and linear in $\vec k$.  Actually, $B_0$
is the matrix whose eigenvalues are the indicial roots, which we computed in section \ref{throots}.  $B_1(\vec k)$ is the symbol of the operator $\d+\d^*$ on 
$\mathrm{ad}(E)$-valued
dfferential forms on $\R^3$ of all possible degrees; in other words, $B_1(\vec k)$ is the momentum space version of the $\d+\d^*$ operator. 

Let $\Y$ be the space of all solutions of the linear equation $\L  _1F(y)=0$ on $\R_+$, with no condition on the behavior near $y=0$ or $\infty$.  
We can identify $\Y$ with $Y$ by, for example, mapping a solution $F(y)\in \Y$ to its
value $F(1)\in Y$, so $\Y$ has dimension $d$.  Let $\Y_0$ be the subspace of $\Y$ consisting of solutions  that obey the Nahm pole boundary condition at $y=0$ (and
in particular are square-integrable near $y=0$),
and let $\Y_\infty$ be the subspace consisting of solutions that are square-integrable at infinity.  Also, let $d_0=\dim\,\Y_0$, $d_1=\dim\,\Y_1$.  
Finally, let $\Y^*$ be the space of solutions of (\ref{redy}) that obey the Nahm pole boundary condition at $y=0$, and in addition are square-integrable,
and set $d^*=\dim \,\Y^*$.  

$\Y^*$ is simply the intersection $\Y_0\cap\Y_\infty$ of the spaces of solutions that are well-behaved at 0 and at $\infty$.  So 
if $d_0+d_\infty-d\geq 0$ and the subspaces $\Y_0,\Y_\infty\subset \Y$ are generic, then the dimension of $\Y^*$ is $d^*=d_0+d_\infty-d$.  Even if these conditions
do not hold, the index of the operator $\L  _1$ is $d_0+d_\infty-d$.

In the case of the KW equations with Nahm pole boundary conditions, $d_0=d_\infty=d/2$ and therefore the index of $\L  _1$ is 0.  Hence if the kernel of $\L  _1$
vanishes -- as shown in section \ref{vanishing} -- then the cokernel also vanishes.   The fact that $d_0=d/2$ was explained in section \ref{nonregular}, basically as
a consequence of the fact that the indicial roots of $\varphi_a,a_y$ are the negatives of those of $a_a,\varphi_y$ (with some care when some roots vanish).
  Since $B(y,\vec k)$ can be approximated by $B_1(\vec k)$ for $y\to\infty$, solutions of $\L  _1F(y)=0$ that are square-integrable for $y\to \infty$ correspond
  to positive eigenvalues of $B_1(\vec k)$.  So to show that $d_\infty=d/2$, we must
show that $B_1(\vec k)$ has equal numbers of positive and negative eigenvalues.  This follows from the fact that $B_1(\vec k)$ is the symbol of the $\d+\d^*$ operator; its
square is $|\vec k|^2$, so its eigenvalues are $\pm |\vec k|$, and as it is traceless, precisely half of the eigenvalues are positive.
Alternatively, by rotation symmetry, the number of positive eigenvalues of $B_1(\vec k)$ is invariant under $\vec k\to -\vec k$; but since $B_1(-\vec k)=-B_1(\vec k)$,
it has equally many positive and negative eigenvalues.

We have omitted various details here, 
since we turn next to a more direct (and much more widely applicable) proof of the vanishing of the cokernel of the linearized KW operator $\LKW$.  
There are two specific reasons to have included the above material.
First, this explains why the counting of positive 
indicial roots in section \ref{nonregular} is important.  Second, even without a direct calculation of the index of the operator $\L  _1$,
the fact that the problem can be formulated in terms of the vanishing of
the cokernel of this 1-dimensional operator will make it easy to go from 4 to 5 dimensions in section \ref{extension}.

\subsection{Pseudo Skew-Adjointness}\label{pseudo}

Inspection of the equations (\ref{lembo}) and (\ref{wembo}) that determine the indicial roots shows that these roots (in the gauge (\ref{zurimo}), which we assume in what
follows) are odd under exchange of $a_a,\varphi_y$ with $\varphi_a,a_y$.  We make this exchange via the linear transformation
\begin{align}\label{roggo}N\begin{pmatrix}a_a \cr \varphi_y\end{pmatrix}&=\begin{pmatrix} \varphi_a \cr a_y \end{pmatrix} \cr
   N\begin{pmatrix}  \varphi_a \cr a_y\end{pmatrix}&=-\begin{pmatrix} a_a \cr \varphi_y \end{pmatrix}.\end{align}
   The minus sign in the second line does not affect what we have said so far, but will be important shortly.  Taking this minus sign into account, we have
   \begin{equation}\label{donzott}N^2=-1,~~~~N^\dagger=-N ,\end{equation}
   where $N^\dagger$ is the transpose of $N$ in the usual basis given by $a_i$ and $\varphi_i$, or more invariantly the 
   adjoint of $N$ with respect to the quadratic form
   \begin{equation}\label{elbo} -\Tr\,\sum_{i=1}^4\left( a_i^2+\varphi_i^2\right).\end{equation}
      
   Perturbations of the Nahm pole solution that depend only on $y$ are governed by an equation that we schematically write
   $\L  _1(0)\Phi=0$, where $\Phi$ combines all the fields and
   \begin{equation}\label{zorbo}\L  _1(0)=\frac{\d}{\d y}+\frac{B_0}{y} \end{equation}
   is obtained from $\L  _1(\vec k)$ (eqn.\ (\ref{redy})) by setting $\vec k=0$.  The matrix $B_0$ can be read off from (\ref{lembo}) and (\ref{wembo}) (which are derived from the equation $\L  _1(0)\Phi=0$
   by replacing $\d/\d y$ with $\lambda/y$), and by inspection we see that $B_0$ is a  real symmetric matrix in the usual basis, or in other words is real and self-adjoint
   for the quadratic form (\ref{elbo}).  On the other hand, $\d/\d y$ is skew-adjoint.  The adjoint of $\L  _1(0)$ is thus
   \begin{equation}\label{morbo} \L  _1(0)^\dagger=-\frac{\d}{\d y} +\frac{B_0}{y}.\end{equation}  Since the indicial roots are the eigenvalues of $-B_0$,
   the statement that the matrix $N$ reverses the sign of the indicial roots is equivalent to $N B_0=-B_0N$.  We can combine this with
   (\ref{morbo}) as the statement that
   \begin{equation}\label{orbo} \L  _1(0)^\dagger=- N \L  _1(0)N^{-1}. \end{equation}
    
    So far we have just reformulated the symmetry that changes the sign of the indicial roots.  It turns out, however, that (\ref{orbo}) holds without change
    for $\vec k\not=0$:
    \begin{equation}\label{porbo}\L  _1(\vec k)^\dagger = -N\L  _1(\vec k)N^{-1}. \end{equation}
    Once one knows that the kernel of $\L  _1(\vec k)$ is trivial, it immediately follows from (\ref{porbo}) that the cokernel of this operator is also trivial.
    Indeed, the cokernel of $\L  _1(\vec k)$ is the kernel of $\L  _1(\vec k)^\dagger$, but (\ref{porbo}) implies that the kernel of $\L  _1(\vec k)^\dagger$ is obtained
    by acting with $N$ on the kernel of $\L  _1(\vec k)$.
    
    Since eqn.\ (\ref{porbo}) holds for any $\vec k$, this statement can be formulated without introducing momentum eigenstates.  If $\L  $ is the linearization of the KW equations
    around the Nahm pole solution, then
    \begin{equation}\label{homebo}\L  ^\dagger =-N \L   N^{-1}, \end{equation}
    a property that we describe by saying that $\L  $ is pseudo skew-adjoint. 
    We can write $\L  =\partial_y +B$, where $B$ is a self-adjoint\footnote{Self-adjointness of $B$ is not hard to verify by inspection, and is clear in the relation \cite{WittenKtwo} 
    of the KW equations on $W\times \I$, for a one-manifold $\I$, to gradient
    flow equations for the complex Chern-Simons functional on $W$.  Linearization of the gradient flow equation for any Morse function $h$ on a Riemannian manifold $X$
    always produces
    a differential operator $\d/\d y+B$, where $B$ (which is derived from the matrix of second derivatives of the function $h$ and the metric of $X$) is self-adjoint.  In the present example, $X$ is
    essentially the space of complex-valued connections on the bundle $E_\C\to W$, where $E_\C$ is the complexification of $E$, and $h$ is the imaginary 
    part of the Chern-Simons
    functional of such a connection.} first-order differential operator that contains derivatives only along $W$. Then $\L  ^\dagger=-\partial_y+B$
    and (\ref{homebo}) amounts to the statement that
    \begin{equation}\label{zomebo} N BN^{-1}=-B. \end{equation}
   Another equivalent statement is that $\t \L  =N\L  $ is actually self-adjoint.
   
    These assertions hold in much greater generality than perturbing around the Nahm pole solution on $\R^3\times \R_+$.  They hold, as we will see, in perturbing
    around any solution of the KW equations on $W\times \I$ for any oriented three-manifold $W$ and one-manifold $\I$ (endowed with a product metric
    $g_{ab}\d x^a \d x^b+\d y^2$), provided only that $\phi_y$ (which as usual is the component of $\phi$ in the $\I$ direction) vanishes.     
      
    This claim can be verified by inspection of the linearized KW equations.  In doing this, to minimize clutter,
    we write simply $A, \phi$ (rather than $A_{(0)}, \phi_{(0)}$ as in section
    \ref{background}) for a solution of the KW equations about which we wish to perturb.  As usual, we  denote the perturbations about this solution
     as $a,\varphi$, so we consider the condition that $A+s a, \phi+s\varphi$ obeys the KW equations to first order in the small parameter $s$.  The symbol $D_i$ will
    denote a covariant derivative defined using the unperturbed connection $A$ (and the Levi-Civita connection of $W$ if $W$ is not flat).  As already explained,
    we assume that the solution about which we expand obeys $\phi_y=0$, and we describe the background in the gauge $A_y=0$.  However, the gauge condition
    that we impose on the fluctuations is that of eqn.\ (\ref{zurimo}):
    \begin{equation}\label{consoo} D_i a^i+[\phi_a,\varphi^a]=0. \end{equation}
    Given that $A_y=0$, this is equivalent to
    \begin{equation}\label{bonso} \frac{\partial a_y}{\partial y} + D_a a^a+[\phi_a,\varphi^a]=0. \end{equation}
    Now let us compare this gauge condition to the linearization of one of the KW equations, namely the condition $D_i\phi^i=0$.  With $A_y=\phi_y=0$, the linearization
    of this equation gives
    \begin{equation}\label{conso}\frac{\partial \varphi_y}{\partial y} + D_a\varphi^a -[\phi_a,a^a]=0. \end{equation}
    When we transform $a$ and $\varphi$ via (\ref{roggo}), these two equations are exchanged except that the signs are reversed
    for all terms not proportional to $\partial_y$, as predicted in eqn.\ (\ref{zomebo}).

 The other KW equations $F_{ij}-[\phi_i,\phi_j]+\epsilon_{ijkl}D^k\phi^l=0$ behave similarly.  We  write the linearization of 
 these equations in detail\footnote{Our orientation convention
 is such that the antisymmetric
 tensors $\epsilon_{ijkl}$ and $\epsilon_{abc}$ obey $\epsilon_{abcy}=\epsilon_{abc}=-\epsilon_{yabc}$. }  in a way adapted to the split $M=W\times \I$:
 \begin{align}\label{tellme}\partial_y a_a-D_a a_y-[\varphi_y,\phi_a]-\epsilon_{abc}D^b\varphi^c-\epsilon_{abc}[a^b,\phi^c] & = 0 \cr
 \partial_y\varphi_a+[a_y,\phi_a]-D_a\varphi_y-\epsilon_{abc}D^b a^c+\epsilon_{abc}[\phi^b,\varphi^c] &=0. \end{align}
   Again  when we transform $a$ and $\varphi$ via (\ref{roggo}), these two equations are exchanged except that the signs are reversed
    for all terms not proportional to $\partial_y$.                                      

\subsection{Extension To Five Dimensions}\label{extension}

As explained in section \ref{exfive}, the Nahm pole solution can also be used to define a boundary condition on certain elliptic differential equations in five
dimensions that are expected to be relevant to Khovanov homology.  We simply replace $\R^3\times \R_+$ by $\R\times \R^3\times\R_+$, where the first
factor is parametrized by a new ``time'' coordinate $x^0$.  We
reinterpret $\phi_y$ as the component $A_0$ of the connection in the $x^0$ direction, and replace $[\phi_y,\,\cdot\,]$ with $D/D x^0$.  
The equations acquire a rotation symmetry
in $\R^4=\R\times \R^3$.

For the Nahm pole boundary condition in five dimensions to be elliptic,  the linearization $\h \L  $ of the five-dimensional equations
around the Nahm pole solution must have trivial kernel and cokernel.  Using the translation symmetries of $\R^4$, we can consider momentum eigenstates,
proportional to $\exp\left(i\sum_{j=0}^3 k_jx^j\right)$ with a real four-vector $k=(k_0,\dots,k_3)$.  Acting on wavefunctions of this kind, $\h \L  $ becomes a 1-dimensional operator
$\h \L  _1(k)$, acting on functions that depend only on $y$, and it suffices to show that for $k\not=0$, $\h \L  _1(k)$ has trivial kernel and cokernel.  

Using the rotation symmetries of $\R^4$, we can assume that the ``time'' component $k_0$ of $k$ vanishes.  In that case, we are dealing with a time-independent perturbation.
By definition, the five-dimensional equations reduce to the KW equations in the time-independent case (with $A_0$ interpreted as $\phi_y$), so for $k_0=0$,
$\h \L  _1(k)$ coincides precisely with the
corresponding operator $\L  _1(\vec k)$ of the four-dimensional poblem.  
But we already know that the kernel and cokernel of $\L  _1(\vec k)$ vanish;
the kernel vanishes by the vanishing result of section \ref{vanishing}, and the cokernel vanishes because of pseudo skew-adjointness.  So in expanding around the
Nahm pole solution in five dimensions, the kernel and cokernel of $\hat \L  $ vanish, as we aimed to show.

\section{The Nahm Pole Boundary Condition On A Four-Manifold}\label{fourmanifold}

So far, we have described the Nahm pole boundary condition for certain four-  or five-dimensional equations 
on a half-space $\R^4_+$ or $\R^5_+$.  The purpose of the present
section is to explore the Nahm pole boundary condition on a general manifold with boundary.  
In section \ref{induced}, we formulate the Nahm pole boundary condition for the KW equation on an oriented
 four-manifold $M$ with boundary $W$.
 The five-dimensional case is similar but will not be described
here.

Once the KW equation with
Nahm pole boundary condition is defined on a four-manifold with boundary, one can inquire about the index of the linearization $\L$
of this equation. Assuming certain foundational results that we postpone to section \ref{analysis} (such as the fact that $\L$ is Fredholm),
 a simple formal computation that we present in section \ref{indexcalc} determines this index.
The analogous index problem on a five-manifold with boundary will not be treated in the present paper.
The index of an elliptic operator on a five-manifold without boundary is always 0, but this is not necessarily the case on a five-manifold with boundary.  

\subsection{Boundary Conditions on the Connection}\label{induced}

For a homomorphism $\varrho:\frak{su}(2)\to \frak g$, let us say that $\varrho$ is quasiregular if $j_\sigma=0$ does
 not occur in the decomposition of $\frak g$  or equivalently if the commutant $C$ is a finite group.  (This is so if $\varrho$ is principal; for additional examples see Appendix \ref{groups}.)
 On a half-space, for quasiregular $\varrho$,
 the Nahm pole boundary condition on the KW equations implies  that the connection $A$ vanishes along the boundary, since
 the relevant indicial roots are strictly positive.

More generally, on any four-manifold $M$ with boundary $W=\del M$, 
the leading order behavior of the connection $A$ along $W$ is
coupled by the KW equations to the leading order behavior of $\phi$. In particular, if ${\KW}(A,\phi) = 0$ and
$\phi$ and $A$ have expansions 
$\phi = y^{-1}(\sum \frak t_a \d x^a) + \dots$, $A = A_{(0)} + \dots$ near the boundary, then for quasiregular
$\varrho$,  the restriction $A_{(0)}$ of the connection  to 
$W$ is uniquely determined, as we will explain. In expanding the solution near $y=0$ and analyzing $A_{(0)}$, we 
 implicitly use the regularity theorem of section \ref{nonlinreg}, that solutions $(A,\phi)$ have asymptotic
expansions as $y \to 0$. In the calculations below we only use the first few terms of this expansion.
The KW equations determine an entire sequence of relationships between the higher coefficients in the expansion of $\phi$ and $A$, with a certain 
number of these left undetermined because of the gauge freedom. However, we focus
here on the leading order relationships, which signify the rigidity of the Nahm pole boundary condition. 

The main subtlety is to understand the generalization of the formula $\phi=\sum_a\frak t_a\d x^a/y+\dots$ to the case that
the boundary manifold $W$ is not flat.  For this, 
we view $\phi$ as a section of the bundle $\mathrm{Hom}(TW, \ad(E))$.  For a point $\vec x\in W$, 
let $\{e_a\}$ be any orthonormal basis of $T_{\vec x}W$.
Suppose that the leading term in the expansion of $\phi$ is  of order $y^{-1}$.  This leading term must be $\sum_a \frak t_a e_a^*/y$
for some $\frak t_a\in {\mathrm{ad}}(E)_{\vec x}$, and the $y^{-2}$ term in the first of eqns. (\ref{kwe}) below implies that the $\frak t_a$
satisfy the commutation relations of $\frak {su}(2)$.  The classification of $\frak{su}(2)$ subalgebras of $\frak g$ up to conjugacy is discrete
and hence the conjugacy class of $\frak t_a$ is independent of $\vec x$; this is the conjugacy class of some homomorphism $\varrho:\frak{su}(2)\to
\frak g$.  It also follows from the theory of $\frak{su}(2)$ that $\phi_{\varrho}=\sum_a \frak t_a e_a^*$, viewed as a homomorphism from $T_{\vec x}W$ to
the subspace of $\mathrm{ad}(E)_{\vec x}$ spanned by the $\frak t_a$, is an isometry.  
 For $G ={SO}(3)$ or ${SU}(2)$, assuming that $\varrho\not=0$,
 the $\frak t_a$ span all of $\mathrm{ad}(E)$ and we have learned that the polar
 part of $\phi$ determines an isomorphism between $\mathrm{ad}(E)$ and $TW$.  (We assume that $\varrho\not=0$ to avoid many
 exceptions in the following remarks, but the discussion below of the non-quasiregular case applies in particular to $\varrho=0$.)
 
  For $G={SO}(3)$, a knowledge of
 $\mathrm{ad}(E)$ is equivalent to a knowledge of the principal bundle $E$.  For $G=SU(2)$, this is not quite true; the possible
 choices of $E$ -- once the identification of $\mathrm{ad}(E)$ with $TW$ is known -- correspond to spin structures on $W$.
 For $G$ of higher rank, the full story is more complicated, and includes the possibility of twisting by a $C$-bundle,
 where $C$ is the commutant of $\varrho(\frak{su}(2))$ in $G$.  The possibility of this twisting will be reflected in eqn.\ (\ref{zub}) below.
 
We henceforth fix $\phi_{\varrho}$ and consider any solution pair $(A,\phi)$ with Nahm pole given by $\phi_\varrho=\sum_a \frak t_a e_a^*$. 
We assume that both $\phi$ and $A$ are  polyhomogeneous (this will be proved in section \ref{nonlinreg}),
where $\phi$ has leading term $\phi_{(0)} = y^{-1} \phi_{\varrho}$ and $A$ has leading term $A_{(0)}$. 
Insert the expansions for $A$ and $\phi$ into the two equations in \eqref{zobo} and collect the terms with like powers. 
We discuss the coefficients of the powers $y^{-2}$ and $y^{-1}$ in turn.
We first discuss the quasiregular case, which means that 
\begin{equation}
\phi \sim y^{-1}\phi_{\varrho} + y\varphi_1 + \ldots, \qquad A \sim A_{(0)} + y a_1 + \ldots.
\end{equation}
We compute
\begin{equation}
\begin{split}
F_A &= \mathcal O(1), \\ 
\phi \wedge \phi & = y^{-2} \phi_{\varrho} \wedge \phi_{\varrho} + \mathcal O(1), \\ 
\star \d_A \phi & = - y^{-2} \star (\d y \wedge \phi_{\varrho}) + y^{-1} \star \d_{A_{(0)}} \phi_{\varrho} + \mathcal O(1),  \\
\d_A \star \phi & = - y^{-2} \d y \wedge \star \phi_{\varrho} + y^{-1} \d_{A_{(0)}}(\star \phi_\varrho)+ \mathcal O(1), 
\end{split}
\end{equation}
so the KW equations become 
\begin{equation}\label{kwe}
\begin{split}
- y^{-2}( \star \d y \wedge \phi_\varrho + \phi_\varrho \wedge \phi_\varrho) + y^{-1} (\star \d_{A_{(0)}} \phi_\varrho) 
+ \ldots & = 0 \\ - y^{-2} \d y \wedge \star \phi_{\varrho} + y^{-1} \d_{A_{(0)}}(\star \phi_\varrho) + \ldots & = 0.
\end{split} 
\end{equation}
The gauge condition \eqref{elob} does not include any terms with $y^{-2}$ or $y^{-1}$. The coefficient of $y^{-2}$ in the 
first of eqns. (\ref{kwe}) is just Nahm's equation, forcing the $\frak t_a$ to generate an $\frak{su}(2)$ subalgebra of $\mathrm{ad}(E)$,
 while in the second, $\star \phi_{\varrho}$ already contains a $\d y$,
so both terms vanish. 

The coefficients of $y^{-1}$ reduce to
\begin{equation}
i)\ \d_{A_{(0)}} \phi_\varrho = 0, \qquad \mbox{and}\quad ii)\ \d_{A_{(0)}} \star \phi_{\varrho} = 0 
\label{harmonicphi}
\end{equation}
Since  $\star \phi_\varrho$ along $W$ is proportional to $\d y$,  eqn.\ $ii)$ involves only tangential derivatives and does not
involve the $ y$ component $A_{(0)}$. Equation $i)$, on the other hand, may have a $\d y$ term, but let us first 
examine its pullback to $W$.  This restricted equation implies that $\phi_\varrho$ 
intertwines the Levi-Civita connection on $TW$ and the connection $A_{(0)}$ on $\ad(E)$. Indeed, this equation shows
that the part of this connection induced on the image of $\phi_\varrho$ is torsion-free, and since it is a $\frak g$ connection, 
it is also compatible with the Killing metric on $\ad(E)$. Hence its pullback to $TW$ must be the Levi-Civita connection:  
i.e.\ in terms of the orthonormal frame $e_a$, 
\begin{equation*}
\phi_{\varrho}\left( \nabla_{e_a} e_b \right) = \nabla_{e_a} \left(\phi_\varrho(e_b)\right)
\end{equation*}
for all $a, b$. Equivalently, with respect to the product connection on $T^*W \times \ad(E)$, $\nabla \phi_{\varrho} = 0$. 
The second equation in \eqref{harmonicphi} is then automatically satisfied.   eqn.\ $i)$ further implies that $A_{(0)}$ is 
valued in $\varrho(TW)\subset \mathrm{ad}(E)$; its projection onto the orthocomplement of this space vanishes.  Thus,
if $\varrho$ is quasiregular or equivalently if the commutant $C$ is a finite group, the restriction of
 $\mathrm{ad}(E)$ to $W$ and its connection $A_{(0)}$    are uniquely determined  locally in terms of $TW$ with its Levi-Civita connection.
(Globally, depending on $C$,  there may be some discrete choices that generalize
the choice of a spin structure for $G=SU(2)$.)

The discussion up to this point is independent of how $\phi_\varrho$ or $A_{(0)}$ is extended into $M$. The $\d y$ component 
of eqn.\ $i)$ in \eqref{harmonicphi} states that $\nabla^{A_{(0)}}_{\del_y} \phi_\varrho = 0$ at $y=0$.  This is not only a gauge-invariant
condition, but it is also independent of the extension of $A_{(0)}$ into the interior.  In particular, if we choose the gauge so that $(A_{(0)})_y = 0$,
then $\del_y \phi_\varrho = 0$ (which vindicates the fact that we have omitted the $y^0$ term in the expansion for $\phi$). 
In any case, the normal derivative of $\phi_\varrho$ at $W$ is pure gauge. The remaining coefficients in the expansions of $A$ 
and $\phi$ are then determined by these leading coefficients 
and the successive equalities determined by the vanishing of the coefficients of each $y^\lambda$. However, it is increasingly hard 
to extract information from the higher terms; even the coefficients of $y^0$ are not so easy to interpret. 

Finally, consider the non-quasiregular case.  Let $C$, with Lie algebra $\frak c$,  be the
commutant in $G$ of $\varrho(\frak{su}(2))$. The KW equations
would allow the above description of $\phi$ and $A$ to be modified by $\frak c$-valued terms that would be $\O(1)$ for $y\to 0$.
However, as in section \ref{gencase}, we pick an arbitrary subgroup $H\subset C$, with Lie algebra $\frak h$, and we write $\frak h^\perp$
for the orthocomplement of $\frak h$ in $\frak c$.  Then 
we restrict the expansion to take the form
\begin{equation}\label{zub}
\phi \sim y^{-1}\phi_{\varrho} +  \varphi_0 + \ldots, \qquad A \sim A_{(0)} + a_0 + \ldots,
\end{equation}
where $(\varphi_0)_a \in \frak h^\perp$ and $(\varphi_0)_y, a_a \in \frak h$, and $\phi_{\varrho}$, $A_{(0)}$ are as above.
(In general, in the non-quasiregular case, the next term in the expansion is of order $y^{1/2}$.)  In this expansion, we have set to
0  the $\frak h^\perp$-valued part of $a_a, \,(\varphi_0)_y$ and the $\frak h$-valued part of $(\varphi_0)_a$,  and
then, in an appropriate global setting, the KW equations will determine globally the $\frak h$-valued part of $a_a,\,(\varphi_0)_a$
and the $\frak h^\perp$-valued part of $(\varphi_0)_a$.  The justification for this assertion is provided in section \ref{analysis}, where we show that the boundary problem just stated
is well-posed.

Calculating the first few terms in the expansions of $F_A$, $\phi \wedge \phi$, $\star \d_A \phi$ and $\d_A \star \phi$ as before,
we see that the coefficient of $1/y^2$ in ${\KW}(A,\phi)$ vanishes just as before. The $y^{-1}$ coefficient in the term $\phi \wedge \phi$ now
equals $\phi_\varrho \wedge \varphi_0$. This involves terms of the form 
$[\frak t_a, (\varphi_0)_b]$ or $[\frak t_a, (\varphi_0)_y]$, and these vanish since 
$\varphi_0$ is valued in the commutant $ \frak c$ of $\varrho(\frak s \frak u(2))$. Similarly, $A_{(0)}$ must be replaced by $A_{(0)} + a_0$
in each of the two equations in \eqref{harmonicphi}, but 
again this does not modify the above considerations, since $a_0$ is valued in $\frak c$.

\subsection{The Index}\label{indexcalc}
Here, assuming foundational results proved in section \ref{analysis}, we calculate the index of the linearization $\LKW$ of the KW equations.
We do this first on a closed four-manifold, and then on a compact
four-manifold with boundary, with arbitrary Nahm pole boundary conditions. The first case is included 
simply to isolate the contribution from the interior topology of $M$. The computation in the second
case relies on a short argument to show that the index is actually independent of the choice of Nahm pole
boundary condition, or in other words, of the representation $\varrho$. This reduces the computation
to one for the special case $\varrho = 0$, where the computation reduces to a well-known one.

\begin{prop}
Let $(M^4,g)$ be closed, and suppose that ${\KW}(A_{(0)} ,\phi_{(0)}) = 0$. Writing ${\LKW}$ for the linearization of $\KW$ at 
this solution, then 
\begin{equation*}
\mathrm{ind} ({\LKW}) = - (\dim \frak g ) \, \chi(M).
\end{equation*}
\label{index1}
\end{prop}
We have already remarked in section \ref{background} that the symbol of $\LKW$ is the same as that of the twisted 
Gauss-Bonnet operator $(\d + \d^*) \otimes \mathrm{Id}_{\ad(E)} : \wedge^{\mathrm{odd}}M \otimes \ad(E) \to 
\wedge^{\mathrm{even}}M \otimes \ad(E)$. The formula here follows directly from the fact that the 
index of $\d + \d^*$, acting from even forms to odd forms, equals $\chi(M)$, the Euler characteristic of $M$. (Twisting by $E$ does
not affact the index of the twisted Gauss-Bonnet operator even
if $E$ is topologically non-trivial.)

\begin{prop}
Let $(M^4,g)$ be a manifold with boundary, with $g$ cylindrical near $\del M$. 
Then fixing the Nahm pole boundary conditions at $\del M$
with $\varrho = 0$ and $\frak h = \{0\}$, the index of the linearization about any solution is given by
\begin{equation*}
\mathrm{ind}({\LKW}) = - (\dim \frak g) \, \chi(M). 
\end{equation*}
\label{index2}
\end{prop}
This choice of Nahm pole condition is the same as the absolute boundary condition for $\d + \d^*$.
This is again a classical formula.  We could equally well have chosen $\frak h = \frak g$, corresponding
to relative boundary conditions for the Gauss-Bonnet operator, in which case the index equals $- (\dim \frak g) \, \chi(M, \del M)$,
but by Poincar\'e duality, $\chi(M, \del M) = \chi(M)$. This is a special case of the next result, 
which is the main one of this section.

\begin{prop}
Let $(M^4,g)$ be an arbitrary compact manifold with boundary, with $g$ cylindrical near $\del M$, and fix any choice of Nahm pole boundary
condition $\varrho$ at $W = \del M$. Then once again
\begin{equation*}
\mathrm{ind}({\LKW}) = - (\dim \frak g) \, \chi(M). 
\end{equation*}
\label{index3}
\end{prop}
We prove this in two steps. First consider the special case where $M = W \times I$ with a product metric, and with Nahm pole 
boundary condition given by any $\varrho$ at $y = 0$ and with trivial Nahm pole boundary condition ($\varrho = 0$, 
$\frak h = 0$) at $y = 1$. We see immediately, using the pseudo skew-hermitian property of $\LKW$ in this product setting,
that the index vanishes. 

Now let $M$ and $\varrho$ be arbitrary and denote by ${\LKW}_\varrho$ the linearized KW operator
about any solution satisfying the Nahm pole boundary condition associated to $\varrho$, and similarly let
${\LKW}_{\mathrm{rel}}$ denote the linearized KW operator relative to $\varrho = 0$ and $\frak h = \frak g$,
i.e.\ with relative boundary conditions. We apply a standard excision theorem for the index, for example \cite[Prop 10.4]{BBW}, 
which shows that 
\begin{equation*}
\mathrm{ind}({\LKW}_\varrho) - \mathrm{ind}({\LKW}_{\mathrm{rel}}) = \mathrm{ind}({\LKW}_{\varrho, \mathrm{rel}}),
\end{equation*} 
where the operator on the right is the linearized KW operator on the cylinder $\del M \times I$ with Nahm 
pole boundary condition $\varrho$ at one end and with relative boundary condtions on the right. 
We have already shown that the index on the right vanishes, whence the claim.

\section{Analytic theory}\label{analysis}
We now turn to a more careful description of the analytic theory which underpins many of the preceding
considerations. More specifically, we briefly describe some aspects of the theory of linear elliptic uniformly degenerate
equations, all taken from \cite{M-edge} and \cite{MVer}, explain how the results and calculations obtained above 
fit into this theory, and then prove a regularity theorem for solutions of the nonlinear KW equations which 
justifies the calculations in the uniqueness theorem. 

Let us make two comments before proceeding.  The first is that in the special case where $\varrho = 0$, the
Nahm pole boundary condition reduces to a classical elliptic boundary problem, and it is well-known that 
solutions are then smooth up to $W$. The theory described below is the natural extension of those ideas 
which allows us to handle the case where $\varrho \neq 0$ and ${\LKW}$ is no longer a uniformly elliptic operator. 
The second is that to keep the exposition simpler, we focus on the problem in four dimensions.
All of the theory below generalizes immediately to the linearized problem in five dimensions, as does the
application of these results to the regularity of solutions of the nonlinear equations as in section \ref{nonlinreg}. 
The changes required are strictly notational. 

\subsection{Uniformly Degenerate Operators}\label{unifdegg}
Let $M$ be a manifold with boundary, and choose coordinates $(\vec x, y)$ near a boundary point, where 
$\vec x \in U \subset \R^{n-1}$ and $y \geq 0$. A differential operator ${\LKW}_0$ is called uniformly degenerate if in any
such coordinate chart near the boundary it takes the form 
\begin{equation}\label{unifdeg} 
{\LKW}_0 = \sum_{j + |\alpha| \leq m} A_{j \alpha}(\vec x, y) (y\del_y)^j (y\del_x)^{\alpha}; 
\end{equation}
here $(y\del_x)^\alpha = (y\del_{x^1})^{\alpha_1} \ldots (y\del_{x^{n-1}})^{\alpha_{n-1}}$. Such operators are also 
called $0$-differential operators. The key point in this definition is that every derivative is accompanied by a factor 
of $y$. In our setting,  the order $m$ equals $1$ and the coefficients $A_{j \alpha}$ are matrices. Observe that a 
uniformly degenerate can never be uniformly elliptic in the standard sense at $\partial M$
because all the coefficients of the highest order terms vanish when $y=0$.  However, there is an extended notion of ellipticity 
for such operators: ${\LKW}_0$ is said to be an elliptic uniformly degenerate operator if it is elliptic in the standard sense at 
interior points, where $y > 0$, and if in addition, near points of the  boundary, the matrix-valued polynomial obtained by replacing 
each $y\del_{x^a}$ and $y\del_y$ with multiplication by the linear variables $-i k_a$ and $-i k_n$ is invertible when 
$(k_1, \ldots, k_n) \neq 0$. (This formal replacement actually has an invariant meaning; see \cite{M-edge}.) 

The linearized KW operator $\LKW$ is not quite of the form \eqref{unifdeg}; instead, $y\LKW  = \LKW_0$ is an elliptic 
uniformly degenerate operator. This is close enough so that the methods described here can be applied to its analysis equally well. 
To put this into perspective, note that if $\Delta$ is the standard Laplacian on a half-space, then $y^2 \Delta$ is elliptic 
uniformly degenerate, which indicates that the uniformly degenerate theory for the latter operator must therefore reflect the 
well-known properties of the former. In other words, the theory described below subsumes and generalizes the classical 
theory of boundary problems for nondegenerate elliptic operators.  Unlike $\Delta = y^{-2}( y^2 \Delta)$, however,
the operator $\LKW = y^{-1}\LKW_0$ is still degenerate because of the presence of terms involving $1/y$; hence (contrary
to the study of $\Delta$), it is necessary to draw on this uniformly degenerate theory. 

The mapping and regularity properties of solutions of an elliptic uniformly degenerate operator ${\LKW}_0$ hinge on the study of 
two simpler model operators. The first of these: 
\begin{equation}\label{normalop}
N({\LKW}_0) = \sum_{j + |\alpha| \leq m} A_{j \alpha}( \vec x, 0) (s \del_s)^j (s\del_{\vec w})^\alpha, 
\end{equation}
is called the normal operator. This is invariantly defined (up to a linear change of variables) as an operator on the half-space $\R^n_+$, 
naturally identified with the inward-pointing tangent space at the boundary point $(\vec x, 0) \in \partial M$.
To emphasize that it acts on functions defined on this entire half-space, rather than just on a coordinate chart, we write it 
using the linear variables $s \geq 0$, $\vec w \in \R^{n-1}$, which are globally defined on this half-space. The global 
behavior of $N({\LKW}_0)$ on $\R^n_+$ plays a central role in the analysis below.  As a matter of notation, we 
define the normal operator of the linearized KW operator as
\begin{equation}
\label{normalopL}
N({\LKW}) = s^{-1} N({\LKW}_0). 
\end{equation}

Notice that $N({\LKW}_0)$ only depends on $\vec x \in \del M$ as a parameter; there is a different normal operator 
at each point of the boundary, each of which is again uniformly degenerate and elliptic in this sense.  This whole collection of 
operators is called the normal family of $\LKW_0$. In some 
cases, certain crucial features of each $N_{\vec x}(\LKW_0)$ vary with the parameter $\vec x$. Fortunately this does not
happen in our setting; the normal operators at different boundary points all `look the same', and so we shall usually
omit the dependence on $\vec x$.  The normal operator $N({\LKW}_0)$ enjoys considerably more symmetries than 
${\LKW}_0$ itself; namely, it is translation invariant in $w \in \R^{n-1}$ and invariant under dilations $(s,w) \mapsto (\lambda s, 
\lambda w)$, $\lambda > 0$.  Because of these symmetries, it is relatively elementary to study directly, and the goal of 
this entire theory is to show that key properties of these normal operators carry over to ${\LKW}_0$ itself. 

The second model operator is a further reduction, called the indicial operator
\begin{equation}\label{indicialop}
I({\LKW}_0) = \sum_{j \leq m} A_{j 0} (\vec x, 0) (s\del_s)^j,
\end{equation} 
obtained from $N({\LKW}_0)$ by dropping all of the terms $(s\del_{w_j})^{\alpha_j}$ with $\alpha_j > 0$.  Following
the convention above, we also write $I({\LKW}) = s^{-1}I({\LKW}_0)$, when $\LKW$ is the KW operator. 

There is a further, purely algebraic, reduction of the indicial operator obtained by letting $I({\LKW}_0)$ act on the 
elementary functions $s^\lambda$. This yields the indicial family: 
\begin{equation}\label{indfam}
I({\LKW}_0, \lambda) = \sum_{j \leq m} A_{j0}(\vec x, 0) \lambda^j = s^{-\lambda} I({\LKW}_0) s^\lambda,
\end{equation}
where each $s\del_s$ has been replaced by a factor of $\lambda$. 

The reader will notice that the normal operator $N({\LKW})$ was effectively already introduced in section \ref{third}
when we considered the linearized KW operator at the special Nahm pole solution on $\R^4_+$. Moreover, 
we also encountered the indicial family of the linearized KW operator in the matrices on the left hand side of \eqref{lembo} 
and \eqref{wembo}. The indicial roots of ${\LKW}_0$ (or equivalently, of ${\LKW} = s^{-1}{\LKW}_0$) are the finite set of 
values of $\lambda$ for which $I({\LKW}_0, \lambda)$ is not invertible. As we have seen, their computation
in our setting involves nontrivial algebraic subtleties. 

For the rest of this discussion, let us consider only the case where ${\LKW}$ is the linearized KW operator.
Everything we say here has analogues for operators of higher order. With this assumption, one obvious 
simplification is that the equation characterizing the indicial roots
is a simple (generalized) eigenvalue problem, namely that the matrix
\begin{equation}\label{indrts}
A_{10} \lambda + A_{00}
\end{equation}
has nontrivial kernel. Writing the fields on which ${\LKW}_0$ acts as $\Psi$, then corresponding to each
indicial root $\lambda$ there is an eigenvector $\Psi_\lambda$. Equivalently, there is a solution of
the indicial operator of the form $\Psi_\lambda s^\lambda$.   In general the indicial roots
may depend on the basepoint $\vec x \in \partial M$, and this introduces substantial analytic
complications. Fortunately, in our case, this does not occur and we assume henceforth that
the indicial roots are constant in $\vec x$. 

The importance of these indicial roots can be explained at various levels.  At the simplest level,
they provide the expected growth or decay rates of solutions to the equation ${\LKW}_0 \Psi = 0$. 
There is no a priori guarantee that actual solutions to this linear PDE actually do grow or
decay at these precise rates, and the fact that they do in some cases is a regularity theorem.
The discussion in the earlier part of this paper assumes that these growth rates
are legitimate, and we are now in the process of showing that this is so for the linearized KW equation
with the Nahm pole boundary condition. 

To proceed further, we pass from the normal operator $N({\LKW})$ to the same operator conjugated
by the Fourier transform in $\vec w$, just as in section \ref{index}. This leads to the matrix-valued
ordinary differential operator 
\begin{equation}\label{ftnormal}
\widehat{N}({\LKW}) = A_{10} \del_s - i A_{0a} k^a + \frac{1}{s}A_{00}
\end{equation}
The factor of $i$ appears since $e^{i \vec k \cdot \vec w} (s\partial_{\vec w}) e^{-i \vec k\cdot \vec w} = -ik^a$. There are two key 
facts about solutions of this operator, each following from elementary considerations:
\begin{itemize}
\item[i)] Any solution $\hat \Psi(s,\vec k)$ to $\widehat{N}({\LKW}) \hat \Psi = 0$ either 
decays exponentially or else grows exponentially as $s \nearrow \infty$;
\item[ii)] Any solution of this equation near $s = 0$ has a complete (and in fact convergent) asymptotic  
expansion 
\begin{equation}\label{expnorsoln}
\hat \Psi(s,\vec k) = \sum_\lambda \sum_{j=0}^\infty \hat \Psi_{\lambda j} s^{\lambda + j},
\end{equation}
where the first sum is over indicial roots of ${\LKW}$. (In exceptional cases, where the difference between different indicial
roots is an integer, this sum may include extra logarithmic factors. This actually happens in the case of the KW equations.) 
\end{itemize}

The second of these assertions is an immediate consequence of the classical theory of Frobenius series of
solutions of equations with analytic coefficients near regular singular points. The first assertion is slightly more subtle
in that it depends on the ellipticity of the normal operator $N({\LKW})$.  The dominant terms in \eqref{ftnormal} as $s \to \infty$ 
are the first two on the right. Dropping the third term $A_{00}$, we obtain the constant coefficient operator 
$A_{10} \del_s - i A_{0a}k^a$, which has solutions of the form $\tilde \Psi_\lambda e^{\lambda s}$ where $\lambda$ and $\tilde \Psi_\lambda$ 
satisfy the algebraic eigenvalue equation $(A_{10} \lambda - i A_{0a}k^a )\tilde \Psi_\lambda = 0$. The fact that there are no solutions
of this equation with 
purely imaginary $\lambda$ (or with $\lambda=0$, $\vec k\not=0$)  is a restatement of the ellipticity, in the
ordinary sense, of the operator $A_{10} \del_s + A_{0a} \del_{w^a}$. 

Before proceeding with the formulation of boundary conditions, we recall the general notions of conormality and polyhomogeneity
of a field $\Psi$ near $\partial M$; these are simple and useful extensions of the notion of smoothness up to the
boundary of $\Psi$.   We say that $\Psi$ is conormal of order $\lambda_0$, and write $\Psi \in \mathcal A^{\lambda_0}$,  
if $y^{-\lambda_0}|\Psi| \leq C$ with a similar estimate for all its derivatives, i.e.\ $y^{-\lambda_0}| (y\del_y)^j (\del_{\vec x})^\alpha 
\Psi| \leq C_{j \alpha}$ for all $j, \alpha$.  Any such field is smooth in the interior of $M$, 
but these estimates give only a very limited sort of smoothness near the boundary: for example, both $y^{\sqrt{-1}}$ and 
$1/\log y$ lie in $\mathcal A^0$. A more tractable subclass consists of the space of polyhomogeneous fields $\Psi$. 
Here $\Psi$ is said to be polyhomogeneous at $y=0$ if it is conormal and in addition has an asymptotic expansion 
\begin{equation}
\Psi \sim \sum y^{\gamma_j} (\log y)^p \Psi_{jp}(\vec x).
\label{phg}
\end{equation}
The exponents $\gamma_j$ on the right lie in some discrete set $E \subset \mathbb C$, called the index set of $\Psi$, 
which has the following properties: $\mathrm{Re}\, \gamma_j \to \infty$ as $j \to \infty$, the powers $p$ of
$\log y$ are all nonnegative integers, and there are only finitely many such log terms accompanying any $y^{\gamma_j}$.  
Notice that the conormality of $\Psi$ implies that each 
$\Psi_{jp}(\vec x) \in \mathcal C^\infty(\del M)$. The meaning of $\sim$ is the classical one for an asymptotic
expansion: namely, 
\begin{equation*}
| \Psi - \sum_{j \leq N} y^{\gamma_{j}} (\log y)^p \Psi_{jp} | \leq C y^{\mathrm{Re} \gamma_{N+1}}(\log y)^q,
\end{equation*}
where the term on the right is the next most singular term in the expansion.  The corresponding statement must
hold for the series obtained by differentiating any finite number of times.  If the $\gamma_j$ are all nonnegative integers 
and the log terms are absent, this is just the standard notion of smoothness up to $y = 0$. 
Solutions of uniformly degenerate equations ${\LKW}_0 \Psi = 0$ are typically polyhomogeneous (at least in favorable 
circumstances), but essentially never smooth in the classical sense. In our specific problem, the exponents
$\gamma_j$ are of the form $\gamma_j = j/2$, $j = 0, 1, 2, \ldots$; 
log  terms, if they appear at all, do not occur in the leading terms. 

\subsection{Elliptic Weights}\label{ellweights}
We now turn to the various sorts of boundary conditions that can be imposed on the operator ${\LKW}$ and a description of what 
makes a boundary condition elliptic (relative to ${\LKW}$).  General types of boundary conditions can be local and of `mixed' Robin 
type, or nonlocal, such as an Atiyah-Patodi-Singer type boundary condition. We shall focus, however, on the particular
local algebraic boundary conditions which arise in the Nahm pole setting. In this section we describe the simplest of these
boundary conditions, where ${\LKW}$ acts on fields with a prescribed rate of vanishing or blowup at $y = 0$. This is 
analogous to a homogeneous Dirichlet condition (which is tantamount to considering solutions 
which vanish like $y^\epsilon$ at the boundary for any $0 < \epsilon < 1$).  This type of boundary condition is relevant in our 
setting only when none of the $j_\sigma = 0$; in particular, this is the correct type of boundary condition when $\varrho$ is a 
regular representation.  This case is simpler to state, and considerably simpler to analyze, than the more general one when 
some of the $j_\sigma=0$, which we come to only in  section \ref{genebc}.  As before, we continue to focus exclusively on the 
linearized KW operator ${\LKW}$, or its uniformly degenerate associate ${\LKW}_0 = y{\LKW}$, although all the results below 
have analogs for more general elliptic uniformly degenerate operators. 

We shall study the action of $\LKW$ on weighted $0$-Sobolev spaces $y^{\lambda_0 + 1/2}H^k_0(M)$, and so we 
start by defining these. First consider $y^{\lambda_0 + 1/2}L^2(M)$, which consists of all fields $\Psi = y^{\lambda_0 + 1/2} \Psi_1$ 
where $\Psi_1 \in L^2(M)$. The measure is always assumed to equal $\d\vec x \,\d y$ up to a smooth nonvanishing multiple. 
Next, for $k \in \mathbb N$, let
\begin{equation*}
H^k_0(M) = \{ \Psi \in L^2(M): (y \del_{\vec x})^\alpha (y\del_y)^j \Psi \in L^2(M),\ \forall\ j + |\alpha| \leq k\},
\end{equation*}
and finally, define $y^{\lambda_0 + 1/2} H^k_0 = \{\Psi = y^{\lambda_0 + 1/2}\Psi_1: \Psi_1 \in H^k_0\}$. The spaces $s^{\lambda_0 + 1/2} H^k_0(\R^n_+)$ 
are defined similarly. The subscript $0$ on these Sobolev spaces indicates that they are defined relative to the $0$-vector fields
$y\del_y$ and $y\del_{x^a}$; it does not connote that the fields have compact support. A key feature of these spaces is that their norms have 
a scale invariance coming from the  invariance of $y\partial_y$ and $y\partial_{x^a}$ under dilations $(y,\vec x) \mapsto
(cy, c\vec x)$, $c > 0$.  In fact, $N({\LKW}_0)$ does not act naturally on the more standard Sobolev spaces defined
using the vector fields $\partial_y$, $\partial_{x^a}$. The shift by $1/2$ in these weight factors is for notational convenience only and
corresponds to the fact that the function $y^{\lambda_0}$ lies in $y^{\lambda_0 + 1/2 + \epsilon}L^2$ locally near $y = 0$ 
when $\epsilon < 0$ but not when $\epsilon \geq 0$. In other words, $y^{\lambda_0}$ just marginally fails to lie 
in $y^{\lambda_0+1/2}L^2$ (near $y=0$).  

It is evident that 
\begin{equation}
{\LKW}: y^{\lambda_0 + 1/2}H^1_0(M) \longrightarrow y^{\lambda_0 - 1/2} L^2(M)
\label{map11}
\end{equation}
is a bounded mapping for any $\lambda_0 \in \R$. Notice that the weight on the right has dropped by $1$, reflecting that the operator
${\LKW}$ involves the  terms $\del_y$ and $1/y$; if we were formulating this using ${\LKW}_0$, then it would be appropriate to
use the same weight on the left and the right.  Our main concern is whether this mapping is Fredholm, i.e.\ has closed
range and a finite dimensional kernel and cokernel, and to describe the regularity of solutions of ${\LKW} \Psi = 0$ (or ${\LKW} \Psi = f$ 
for fields $f$ which have better regularity and decay as compared to $y^{\lambda_0 - 1/2}L^2$).  It is not hard to show that \eqref{map11} does 
not have closed range when $\lambda_0$ is an indicial root of $\LKW$. Indeed, in this case, an appropriate sequence of compactly
supported cutoffs of the function $y^{\lambda_0}$ can be used to create a Weyl sequence $\Psi_j$, i.e.\ an orthonormal sequence of fields 
such that 
\begin{equation*}
||\Psi_j||_{y^{\lambda_0+1/2}H^1_0} =1, \qquad || {\LKW} \Psi_j||_{y^{\lambda_0 - 1/2}L^2} \to 0.
\end{equation*}
Hence \eqref{map11} is certainly not Fredholm then. When $\lambda_0 \ll 0$, then \eqref{map11} has
an infinite dimensional kernel, while if $\lambda_0 \gg 0$, its cokernel is infinite dimensional. Thus the only chance for \eqref{map11}
to be Fredholm is when $\lambda_0$ is nonindicial and lies in some intermediate range. In some cases (as described in 
section \ref{genebc}), it is not Fredholm for any weight $\lambda_0$. Closely related is the fact that fields in the kernel of \eqref{map11} 
when $\lambda_0 \ll 0$ are, in general, not regular, i.e.\ polyhomogeneous; indeed, for such $\lambda_0$, most solutions are quite 
rough at $\del M$. All of this motivates the following definition.

\begin{definition}\label{ellbc}
The weight $\lambda_0$ is called elliptic for the linearized KW operator ${\LKW}$ if its normal operator defines an invertible mapping:
\begin{equation}\label{noropmap}
N({\LKW}): s^{\lambda_0 + 1/2} H^1_0( \R^n_+; \d\vec w \,\d s) \longrightarrow s^{\lambda_0 - 1/2} L^2(\R^n_+; \d\vec w\, \d s). 
\end{equation}
\end{definition}

The two main consequences of the ellipticity of a weight are stated in the following propositions.
\begin{prop}\label{fred1}
Let ${\LKW}$ be the linearized KW operator on a compact manifold with boundary $M$, and suppose that $\lambda_0$ 
is an elliptic weight. Then the mapping  \eqref{map11} is Fredholm. 
\end{prop}

Recalling that we are in the case where no $j_\sigma = 0$, let $\underline{\lambda}$ and $\overline{\lambda}$ be the 
largest negative and smallest positive indicial roots of $\LKW$, respectively. Thus (see Appendix \ref{groups}) $\underline{\lambda} = -1$ 
and $\overline{\lambda} =1$. We assert that any $\lambda_0 \in (\underline{\lambda}, \overline{\lambda})$ 
is an elliptic weight. We shall prove this later; in fact, all of the necessary facts for the proof come from the 
considerations in sections \ref{second} and \ref{third}.  It will follow from this proof that if some $\lambda_0$ is an elliptic
weight, then so is any other $\lambda_0'$ which lies in a maximal interval around $\lambda_0$ containing
no indicial roots. 

\begin{prop}\label{regularity1}
With all notation as above, let $\lambda_0$ be an elliptic weight for ${\LKW}$, and suppose that ${\LKW}\Psi = f \in y^{\lambda_0 - 1/2}L^2$ 
where $\Psi \in y^{\lambda_0 + 1/2} L^2$. If $f$ is smooth in a neighborhood of some point $q \in \del M$, or slightly more generally, 
if it has a polyhomogeneous asymptotic expansion as $y \to 0$, then in that neighborhood, $\Psi$ admits a polyhomogeneous expansion 
\begin{equation}
\Psi \sim \sum y^{\gamma_j} \Psi_j,
\label{expansion}
\end{equation}
where the exponents $\gamma_j$ are of the form $\lambda + \ell$, $ \ell \in \mathbb N$, where either $\lambda$ 
is an indicial root of ${\LKW}$ or else is an exponent occurring in the expansion of $f$. 
\end{prop}
The expansion \eqref{expansion} may contain terms of the form $y^{\gamma_j} (\log y)^p$, $p > 0$. These can only appear when 
the differences between indicial roots are integers (as happens in our setting), or when there is a coincidence between the 
indicial roots of ${\LKW}$ and the terms in the expansion of the inhomogeneous term $f$.  However, the key fact is simply 
that $\Psi$ has an expansion at all; once we know this, then the precise terms in its expansion can be determined by matching 
like terms on both sides of the equation ${\LKW} \Psi = f$.

The existence of such an asymptotic expansion for solutions should be regarded as a satisfactory replacement for smoothness
up to the boundary. For an ordinary (nondegenerate) elliptic operator ${\LKW}_0$, if the standard Dirichlet condition (requiring 
solutions to vanish at $y=0$) is an elliptic boundary condition in the classical sense, then solutions of ${\LKW}_0 \Psi = 0$ with
$\Psi(0,\vec x) = 0$ are necessarily smooth and vanish to order $1$ at the boundary. Proposition \ref{regularity1} is the exact 
analogue of this. For the linearized KW operator, the nonnegative indicial roots lie in the set $\{j/2: j = 0, 1, 2, \ldots\}$,
so if ${\LKW} \Psi = 0$ in some neighborhood of a boundary point, then 
\begin{equation*}
\Psi \sim  \sum y^{j/2} \Psi_j,
\end{equation*}
(we are neglecting log terms which might appear); in general, there are half-integral exponents, so $\Psi$ is genuinely not smooth at $y=0$. 

We now verify the invertibility of \eqref{noropmap} for every $\lambda_0 \in (\underline{\lambda}, \overline{\lambda})$ 
in our particular example, which proves the assertion that every such $\lambda_0$ is an elliptic weight. First 
conjugate $N({\LKW})$ with the Fourier transform in $\vec w$, thus passing to the simpler ordinary differential operator 
$\widehat{N}({\LKW})$ as in \eqref{ftnormal}.  We must show that
\begin{equation}
\widehat{N}({\LKW}): s^{\lambda_0+1/2} H^1_0(\R^+; \d s) \longrightarrow s^{\lambda_0 - 1/2} L^2( \R^+; \d s),
\label{noropmapft}
\end{equation}
is invertible for each $\vec k$, and that the norm of its inverse is bounded independently of $\vec k$.

The first step is to use the scaling properties of $\widehat{N}({\LKW})$ to reduce to the case $|\vec k| = 1$. Indeed, 
set $t = s|\vec k|$ and write $B({\LKW}) = A_{10} \partial_t - it A_{0a} k^a/|\vec k| + A_{00}$. (The ``$B$'' refers to the
fact that this operator which has many features in common with the Bessel equation, and so we call $B({\LKW})$ 
the model Bessel operator of $\LKW$.)  Applying this change of variables replaces $B({\LKW})$ by $|\vec k|^{-1} \widehat{N}({\LKW})$.

Suppose that we have already shown that the version of \eqref{noropmapft} with $B({\LKW})$ replacing $\widehat{N}({\LKW})$ 
is invertible for every $\vec k$ with $|\vec k| = 1$, and let $B(G)(t, t', \vec k)$ denote the Schwartz kernel of this inverse. 
We then recover the Schwartz kernel $G(s,s', \vec k)$ of the inverse of \eqref{noropmapft} for any $\vec k \neq 0$ as  
\begin{equation*}
\widehat{N}(G)(s, s', \vec k) = B(G)(s|\vec k|, s' |\vec k|, \vec k/|\vec k|). 
\end{equation*}
To see that this is the case, we first compute that 
\begin{equation*}
\begin{split}
\widehat{N}({\LKW}) \int & B(G)(s|\vec k|, s'|\vec k|, \vec k/|\vec k|) f(s', \vec k)\, \d s' = 
\int (B({\LKW}) B(G)) (s|\vec k|, s'|\vec k|, \vec k/|\vec k|) |\vec k| f(s', \vec k) \, \d s' \\
& = \int \delta(s|\vec k| - s' |\vec k|) |\vec k| f(s', \vec k)\, \d s' = f(s, \vec k),
\end{split}
\end{equation*}
since the $\delta$ function in one dimension is homogeneous of degree $-1$. This result may seem counterintuitive 
since one expects that $||\widehat{N}({\LKW})^{-1}|| \sim 1/|\vec k|$, but that expectation is false because we are
letting $\widehat{N}({\LKW})$ act between spaces with different weight factors. In fact, the norm of 
$\widehat{N}(G)$ is bounded uniformly in $\vec k$, but does not decay as $|\vec k| \to \infty$. To this end,
observe that we must estimate the norm of 
\begin{equation*}
H(s, s', \vec k) := s^{-\lambda_0 - 1/2} \widehat{N}(G)(s,s', \vec k) (s')^{\lambda_0 - 1/2}: L^2 \to L^2.
\end{equation*}
This is done by calculating
\begin{multline*}
\int \left| \int H(s, s', \vec k) f(s', k)\, \d s' \right| ^2 \, \d s  \\
= \int \left| \int B(G)(s|\vec k|, s'|\vec k|, \vec k/|\vec k|) (s|\vec k|)^{-\lambda_0 - 1/2} (s' |\vec k|)^{\lambda_0 - 1/2} 
|\vec k| f(s', \vec k)\, \d s'\right|^2 \, \d s \\ = \int \left| \int B(G)(t, t', \vec k/|\vec k|) 
t^{-\lambda_0 - 1/2}(t')^{\lambda_0 - 1/2} |\vec k|^{-1/2} f(t'/|\vec k|, \vec k)\, \d t' \right|^2 \\
\leq C || |\vec k|^{-1/2} f( t/|\vec k|, \vec k) || = C ||f (\cdot, \vec k)||.
\end{multline*}
The inequality in the fourth line reflects the boundedness of $B(G): t^{\lambda_0 - 1/2} L^2 \to t^{\lambda_0 + 1/2}L^2$. 

Beyond all this, compactness of the unit sphere in $\vec k$ shows that the norm of $B(G)$ can be bounded
independently of $\vec k$. 

As for showing that \eqref{noropmapft} is invertible for each $\vec k$, we first show that it is Fredholm. This can be
done by a standard ODE analysis of the operator.  First construct approximate local inverses near $s=0$ and $s=\infty$;
the existence of these shows that \eqref{noropmapft} is Fredholm precisely when $\lambda_0$ is {\it not} 
an indicial root of ${\LKW}_0$. (As noted earlier, when $\lambda_0$ is an indicial root, this mapping does not have closed range.)  
Now fix $\lambda_0$ to be any nonindicial value and recall the fact i) that any element of the kernel either grows or decays exponentially 
as $s \to \infty$. Then injectivity of this map means precisely that the solutions which decay exponentially as $s \to \infty$ 
must blow up faster than $s^{\lambda_0}$ as $s \to 0$.  One can perform a similar analysis for the adjoint operator,
or by showing  by other methods that the index vanishes, to show that \eqref{noropmapft} is surjective too.

For the linearized KW operator ${\LKW}$, the discussion in section \ref{nonregular} implies directly that the kernel of 
$\widehat{N}({\LKW})$ has only trivial kernel on $s^{\lambda_0 + 1/2}L^2$  when $\lambda_0 > 0$.  Indeed, 
perform the integrations by parts (which are now only in the $s$ variables), using the decay of solutions both 
as $s \to 0$ and as $s \to \infty$ to rule out contributions from the boundary terms. We can extend this to allow
any $\lambda_0 > \underline{\lambda}$ simply by observing using the fact ii) about solutions that if $\widehat{N}({\LKW}) 
\hat\Psi = 0$ and $\hat \psi \in s^{\lambda_0 + 1/2}L^2$, then $\hat\Psi$ vanishes like $s^{\overline{\lambda}}$, so we may
integrate by parts as before. The fact that the index vanishes when $\lambda \in (\underline{\lambda}, \overline{\lambda})$
follows by using the pseudo skew-adjointness established in section \ref{pseudo}. Note that those arguments are for ${\LKW}$ on the model 
space $\R^4_+$, which is canonically identified with the normal operator $N({\LKW})$ of the linearized KW equations on 
any manifold with boundary, and the pseudo skew-adjointness passes directly to $\widehat{N}({\LKW})$ as well.

We have now proved that in the quasiregular case, when no $j_\sigma  = 0$, any $\lambda_0 \in (\underline{\lambda}, \overline{\lambda})$ 
is an elliptic weight for $\LKW$. 

The results just stated are not well suited for our nonlinear problem simply because these weighted $L^2$ spaces do not behave 
well under nonlinear operations. One hope might be to use $L^2$- (or $L^p$-) based scale-invariant Sobolev spaces with 
sufficiently high regularity. These do have good multiplicative properties locally in the interior, but not near the boundary. 
This leads us to introduce several related H\"older-type spaces, and then describe the mapping properties of ${\LKW}$ acting on them. 

We start with the spaces $\mathcal C^k_0$, which consist of all fields $\Psi$ such that
$(y\partial_y)^j (y\partial_{\vec x})^\alpha \Psi$ is bounded on $M$ and continuously differentiable in
the interior of $M$ for all $j + |\alpha| \leq k$. The H\"older seminorm is defined by
\begin{equation*}
[ \Psi]_{0; 0,\gamma} = \sup_{(y,\vec x) \neq (y', \vec x\, ')} \frac{ |\Psi(y, \vec x) - \Psi(y', \vec x\, ')|(y+y')^\gamma}{
|y-y'|^\gamma + |\vec x - \vec x\, '|^\gamma} 
\end{equation*}
Then $\mathcal C^{k,\gamma}_0$ consists of all $\Psi \in \mathcal C^k_0$ such that 
$ [ (y\partial_y)^j (y\partial_{\vec x})^\alpha \Psi]_{0; 0,\gamma} < \infty$. Finally, $y^{\lambda_0} \mathcal C^{k,\gamma}_0$ 
consists of all $\Psi = y^{\lambda_0}\Psi_1$ with $\Psi_1 \in \mathcal C^{k,\gamma}_0$. 

These spaces capture no information about regularity in the $\vec x$ directions at $y=0$, so we also introduce hybrid spaces
\begin{equation*}
\mathcal C^{k,\ell, \gamma}_0 = \{ \Psi \in \mathcal C^{k,\gamma}_0: (\partial_{\vec x})^\alpha \Psi \in 
\mathcal C^{k-|\alpha|,\gamma}_0,\ \mbox{for all}\ |\alpha| \leq \ell\}.
\end{equation*}
Note that all of these spaces contain elements like $y^\lambda$ or $y^\lambda (\log y)^p$, provided
$\lambda > \lambda_0$ (or $\lambda = \lambda_0$ if $p = 0$). 

The mapping property of $\LKW$ on these spaces is much the same as in Proposition \ref{fred1}.
\begin{prop}
Let ${\LKW}$ be the linearized KW operator. Suppose that no $j_\sigma = 0$ and let $\lambda_0 \in (\underline{\lambda},
\overline{\lambda})$ be an elliptic weight. Then the mapping 
\begin{equation}
{\LKW}: y^{\lambda_0} \mathcal C^{k,\ell, \gamma}_0 \longrightarrow y^{\lambda_0 - 1} \mathcal C^{k-1, \ell, \gamma}_0
\label{Fredholder}
\end{equation}
is Fredholm for $0 \leq \ell \leq k-1$ and $k \geq 1$. 
\label{maphold}
\end{prop}
We explain in the next section how this result is essentially a corollary of Proposition~\ref{fred1}. More specifically, 
both results are proved by parametrix methods; this parametrix is constructed using $L^2$ methods and it is initially
proved to be bounded between weighted Sobolev spaces, but it is also bounded between certain of these weighted 
H\"older spaces, which leads directly to the proof of Proposition \ref{maphold}.

\subsection{Structure of the Generalized Inverse}\label{strgeninv}
We now briefly describe the technique behind the proofs of these results. The main step in each is the 
construction and use of the generalized inverse $G$ for \eqref{map11}. The ellipticity of the weight enters
directly into this construction.  By definition, a generalized inverse for \eqref{map11} is a bounded operator 
$G: y^{\lambda_0 - 1/2}L^2 \to y^{\lambda_0 + 1/2}H^1_0$ which satisfies
\begin{equation}
G {\LKW} = \mbox{Id} - R_1, \qquad {\LKW} G = \mbox{Id} - R_2,
\label{geninv}
\end{equation}
where $R_1$ and $R_2$ are finite rank projections onto the kernel and cokernel of $\LKW$, respectively. 
The nonuniqueness here is mild and results only from the different possible choices of projectors. Since
we are working on a specific weighted $L^2$ space, it is natural to demand that $R_1$ and $R_2$ 
be the orthogonal projectors onto the kernel and cokernel with respect to that inner product, and
we make this choice henceforth. 

If we already know that \eqref{map11} is Fredholm, then general functional analysis tells us that a generalized
inverse exists. Conversely, the existence of an operator $G$ with these properties (the boundedness of 
$G$ is particularly important) implies that \eqref{map11} is Fredholm. In fact, it is only necessary to find a 
bounded operator $\tilde{G}$ such that the `error terms' $R_1$ and $R_2$ defined as in \eqref{geninv} are 
compact operators, for then standard abstract arguments imply that \eqref{map11} is Fredholm and show that $\tilde{G}$ 
can be corrected to an operator such that \eqref{geninv} holds, with $R_1$ and $R_2$ the actual projectors.  
This observation is important because it is certainly easier to construct an intelligently designed approximation 
to the generalized inverse than to construct the precise generalized inverse directly.  The criterion by which
one judges the approximation to be good enough is simply that the remainder terms $R_1$ and $R_2$ are compact.
An approximation of this type is called a parametrix for $\LKW$. 

A parametrix can be constructed within the framework of geometric microlocal analysis, as carried out in full
detail in \cite{M-edge}. The key point is to work within a class of pseudodifferential operators on $M$ 
adapted to the particular type of singularity exhibited by $\LKW$. This is the class of $0$- (or uniformly degenerate) 
pseudodifferential operators, $\Psi_0^*(M)$. We wish that $\Psi^*_0(M)$ is sufficiently large to contain parametrices 
of elliptic uniformly degenerate differential operators, but not so large that the individual operators in $\Psi^*_0(M)$ 
are too unwieldy to analyze.  We describe these operators in sufficient 
detail for the present purposes, but refer to \cite{M-edge} for further details.

The elements of $\Psi^*_0(M)$ are characterized by the singularity structure of their Schwartz kernels. Thus, an operator 
$A \in \Psi^*_0(M)$ has a Schwartz kernel $\kappa_A(y, \vec x, y', \vec x\,')$, which is a distribution on $M^2$. 
We expect it to have a a standard pseudodifferential singularity (generalizing that of the Newtonian potential,
for example) along the diagonal $\{y = y', \vec x = \vec x\,'\}$, but we also require a very precise regularity along 
the boundaries of $M^2$, $\{y = 0,\ \mbox{or}\ y' = 0\}$, and at the intersection of the diagonal with the boundary. 
To formalize this, introduce the space $M^2_0$ obtained by taking the real blowup of the product $M^2$ 
at the boundary of the diagonal. In local coordinates this means that we replace each point $(0, \vec x, 0, \vec x)$
in the boundary of the diagonal with its inward-pointing normal sphere-bundle. Alternatively, in polar coordinates
\begin{equation*}
R = (y^2 + (y')^2 + |\vec x - \vec x\,'|^2)^{1/2}, \ \omega = (\omega_0, \omega_0', \hat\omega) = 
(y, y', \vec x - \vec x\,')/R \in S^4_+,
\end{equation*}
where $S^4_+$ consists of the unit vectors in $\R^5$ with $\omega_0, \omega_0' \geq 0$, we replace each point 
$(0,\vec x, 0, \vec x)$ by the quarter-sphere at $R = 0$. We can then use $(R, \omega, \vec x\,')$ as a full set of coordinates. 
This new space is a manifold with corners up to codimension three, and has a new hypersurface boundary at $R=0$, 
which is called the front face. Its two other codimension one boundaries
$\omega_0 = 0$ and $\omega_0' = 0$ are called its left and right faces.  There is an obvious blowdown map $M^2_0 \to M^2$.
We now say that $A$ is a $0$-pseudodifferential operator if $\kappa_A$ is the pushforward under blowdown of a distribution
on $M^2_0$ (which we denote by the same symbol) which decomposes in the following fashion
as a sum $\kappa_A(R,\omega,\vec x\,') = \kappa_A' + \kappa_A''$.
Here $\kappa_A'$ is supported away from the left and right faces and has a pseudodifferential singularity of some order $m$ 
along the lifted diagonal $\{\omega_0 = \omega_0', \hat \omega = 0\}$, and if we factor $\kappa_A' = R^{-4} \hat\kappa_A'$,
then $\hat\kappa_A'$ (along with its conormal diagonal singularity) extends smoothly across the front face of $M^2_0$.
This exponent $-4$ is dimensional; in general it should be replaced by the dimension of $M$. On the other hand $\kappa_A''$ is
smooth in the interior of $M^2_0$ and has polyhomogeneous expansions at the left, right and front faces of this space,
with product type expansions at the corners. Altogether, if the expansions at these faces commence with the terms
$\omega_0^a$ (at the left face), $(\omega_0')^b$ (at the right face) and $R^{-4 + s}$ (at the front face), then we write
\begin{equation*}
A \in \Psi_0^{m, s, a, b}(M).
\end{equation*}
Slightly more generally we could replace the superscripts $a$, $b$, denoting the leading exponents of the polyhomogeneous 
expansions at the left and right faces, by index sets, but we do not need this more refined notation here.  

This elaborate notation simply specifies the precise vanishing or blowup properties of $\kappa_A$ in each of these regimes.
We have introduced it out of some necessity since at least some features of this precise structure will be used in an
important way below.  Before proceeding, note one very special case: the identity operator $\mbox{Id}$
is an element in this class, and lies in $\Psi_0^{0, 0, \emptyset, \emptyset}$.  The fact that it has order $0$ along
the diagonal is expected, and since its Schwartz kernel $\delta(y-y') \delta(\vec x - \vec x\,')$ is supported
on the diagonal, its expansion is trivial at the left and right faces, which explains the third and fourth
superscripts. Finally, the second superscript is explained by noting that in polar coordinates
\begin{equation*}
\delta(y-y') \delta(\vec x - \vec x\,') = R^{-4} \delta(\omega_0 - \omega_0')\delta(\hat\omega)
\end{equation*}

Having introduced this general class of pseudodifferential operators, we now explain the parametrix construction.
We aim to find an operator $\tilde{G} \in \Psi_0^*(M)$ such that $\LKW \tilde{G}$ is equal to the identity
up to some compact remainder terms. Rewriting this as the distributional equation 
\begin{equation*}
{\LKW} \kappa_{\tilde G} = R^{-4} \delta(\omega_0 - \omega_0') \delta(\hat\omega)
\end{equation*}
we see that the singularity of $\kappa_{\tilde{G}}$ along the diagonal can be obtained by classical methods
(the symbol calculus), and this construction is uniform as $R \to 0$ once we have removed the appropriate
factors of $R$. In fact, writing $\LKW$ in these same polar coordinates and noting that it lowers homogeneity
in $R$ by $1$, we expect $\kappa_{\tilde{G}}$ to only blow up like $R^{-3}$ at the front face.  In addition,
$\LKW$ must kill the terms in the expansion of $\kappa_{\tilde{G}}$ at the left face, which means that
the terms in the expansion in this face should involve the indicial roots; in particular, the leading exponent at 
this face must equal $\overline{\lambda}$. 

It is not apparent here where the ellipticity of the weight $\lambda_0$ enters. The answer is as follows. After
first solving away the diagonal singularity using the symbol calculus, we must then improve the initial
guess for the parametrix to another one for which ${\LKW} \kappa_{\tilde{G}} - \delta_{\mathrm{Id}}$ vanishes at the front 
face as well. Because the lift of $\LKW$ to $M^2_0$ acts tangentially to the boundary faces of that space,
this equation restricts to an elliptic equation on the front face. Using a natural identification of each quarter-sphere
fiber of the front face with a half-Euclidean space, and a few other steps which we omit, we are led to having to find 
the {\it exact} solution to $N({\LKW}) \Psi = f$ where $f$ is some smooth compactly supported function on $\R^4_+$. 
If we are able to do this, we can then correct the parametrix to all orders so that the remainder terms are clearly compact. 
The natural identification used here is that the quarter-sphere $S^4_+$ fiber in the front face over each point $(0,\vec x, 0, \vec x)$
of the boundary of the diagonal can be identified with the half-space $\R^4_+$, where this identification is unique
up to a projective map, and the restriction to the front face of the lift of $\LKW$ to $M^2_0$ is transformed to 
$N({\LKW})$ in this identification. Section 2 of \cite{M-edge} (especially around eqn.\ (2.10)) explains more about
these identifications. This explains why the exact invertibility of the normal operator plays a crucial role.  

The vindication that this all works is that by carrying out this parametrix construction, one proves that the generalized 
inverse $G$ for $\LKW$ is an element of $\Psi^{-1,1, \overline{\lambda}, b}_0(M)$, where the final index $b$ is some positive 
number related to the indicial roots of the adjoint of $\LKW$, and that the remainder terms in \eqref{geninv} satisfy 
$R_1 \in \Psi^{-\infty, \overline{\lambda},  b}(M)$, $R_2 \in \Psi^{-\infty, b, \overline{\lambda}}(M)$, where this notation (note the 
absence of the subscript $0$) means that their Schwartz kernels are smooth in the interior and polyhomogeneous at 
the two boundary hypersurfaces of $M^2$, rather than being polyhomogeneous on the blown up space $M^2_0$. 

We now explain how to use all of this for our purposes.  Granting this structure of the Schwartz kernel of the
generalized inverse $G$, the boundedness of the map 
\begin{equation}
G: y^{\lambda_0 - 1/2} L^2 \longrightarrow y^{\lambda_0 + 1/2} H^1_0
\label{GL2}
\end{equation}
can be deduced from the standard local boundedness of pseudodifferential operators of order $-1$ between
$L^2$ and $H^1$, and the following inequalities for the orders of vanishing of $\kappa_G$ at the 
various boundary faces. The fact that $\kappa_G$ blows up like $R^{-3}$ at the front face, one order 
better than the Schwartz kernel of the identity, partly explains why $G$ raises the exponent in the weight factor by $1$.
The other aspect which affects the weight on the right in \eqref{GL2} is the leading exponent $\overline{\lambda}$
at the left face, since it is clear that the decay profile of 
\begin{equation}
\int_M \kappa_G( y, \vec x, y', \vec x\,') f(y', \vec x\,') \, \d y' \d\vec x\,'
\label{intexp}
\end{equation}
as $y \to 0$ must incorporate that exponent. Recalling the earlier discussion that $\kappa_G$ decomposes
into the near-diagonal and off-diagonal parts, $\kappa_G'$ and $\kappa_G''$, a close analysis of the
integrals from each of these terms proves that
\begin{equation}
G: y^{\lambda_0 - 1/2}L^2 \longrightarrow y^{\lambda_0 + 1/2} H^1_0 + \bigcap_{\epsilon > 0 \atop k \geq 0 } y^{\overline{\lambda}+1/2 - \epsilon} H^k_0.
\label{refregsob}
\end{equation}
In other words, the near-diagonal part raises regularity and the weight parameter by exactly $1$, while
the off-diagonal part improves the $0$-regularity to an arbitrarily large amount, and has growth/decay
rate of the outcome dictated by the leading term of $\kappa_G$ on the left face. We must intersect
over all $\epsilon>0$ here simply because $y^{\overline{\lambda}} \in y^{\overline{\lambda} + 1/2-\epsilon}L^2$ when $\epsilon > 0$,
but not otherwise. 

From all of this, and observing that when $\overline{\lambda} > \lambda_0$, the range of \eqref{refregsob} is contained 
in $y^{\lambda_0 + 1/2}H^1_0$, we see that \eqref{map11} is Fredholm, 
which is Proposition \ref{fred1}. The proof of Proposition \ref{regularity1} is obtained by a more refined examination of the 
mapping properties of $G$, in particular the fact that if $f$ is polyhomogeneous, then so is the outcome of the integral \eqref{intexp}.
Finally, the proof of Proposition \ref{maphold} can be explained as follows.  The preceding discussion has been 
based on $L^2$ considerations, which is natural since, for example, Fourier analysis has been used
at several points. However, the precise pointwise behavior of the Schwartz kernel of $G$ makes it possible to 
read off its mapping properties on other function spaces.  In particular, we obtain the analog of \eqref{refregsob}
on weighted H\"older spaces: 
\begin{equation}
G: y^{\lambda_0-1} \mathcal C^{k,\gamma}_0  \longrightarrow y^{\lambda_0}\mathcal C^{k+1,\gamma}_0 + 
\bigcap_m y^{\overline{\lambda}} \mathcal C^{m,\gamma}_0.
\label{refreghold}
\end{equation}
As before, this range lies in $y^{\lambda_0} \mathcal C^{k+1,\gamma}_0$ when $\lambda_0 < \overline{\lambda}$.
There is a slight refinement of this which we shall need later, namely that
\begin{equation}
G: y^{\mu-1} \mathcal C^{k,\gamma}_0  \longrightarrow y^{\mu}\mathcal C^{k+1,\gamma}_0 + 
\bigcap_m y^{\overline{\lambda}} \mathcal C^{m,\gamma}_0
\label{refreghold2}
\end{equation}
for any $\mu > \lambda_0$ (and for simplicity of the statement, $\mu$ not an indicial root). 
Since the equations \eqref{geninv} are satisfied as distributions, and every operator in them is bounded between 
the appropriate H\"older spaces, we see that \eqref{Fredholder} is in fact Fredholm, at least for $\ell = 0$.   
To prove that \eqref{Fredholder} is Fredholm for $\ell \geq 0$, we need one extra fact, which is that each of 
the commutators $[ G, \del_{x^a}]$ lies in $\Psi^{-1, 1, \overline{\lambda}, b}_0$, i.e.\ is a $0$-pseudodifferential operator 
with the {\it same} indices, hence has the same mapping properties as $G$ itself. 

This simple transition from Sobolev  to H\"older spaces is a good exemplar of the parametrix method; if we were working 
solely with a priori estimates, then it is no simple matter to deduce Fredholmness on one type of function space from 
the corresponding property on another type of function space.  


\subsection{Algebraic Boundary Conditions and Ellipticity}\label{genebc}
There are many natural operators, however, for which there are no ellliptic weights, i.e.\ so that for any nonindicial $\lambda_0$, 
the map \eqref{noropmapft} has either nontrivial kernel or cokernel, or both. This is the case for the linearized KW operator ${\LKW}$
when some of the $j_\sigma = 0$. We now describe the somewhat more complicated formulation of the ellipticity criterion
for boundary conditions in these cases. As before, we consider only the parts of this story relevant to the Nahm pole  boundary 
conditions for ${\LKW}_0$.  

The argument in the last section (slightly after \eqref{noropmapft}) shows that even when the lowest nonnegative indicial 
root $\lambda$ is $0$, then $\hat{N}({\LKW}) : s^{\lambda_0 + 1/2}L^2 \to s^{\lambda_0-1/2}L^2$ has no kernel if $\lambda_0 > 0$,
although the cokernel of this mapping is nontrivial. On the other hand, when $\lambda_0 < 0$, the nullspace has positive dimension,
though the map is surjective. To be definite, suppose that $-1/2 < \lambda_0 < 0$, which rules out solutions which blow up like 
$s^\lambda$ where $\lambda$ is any one of the strictly negative indicial 
roots of ${\LKW}_0$. Using the fact ii) that solutions of $\widehat{N}({\LKW}) \hat{\Psi} = 0$ have (convergent) 
expansions at $s = 0$, the leading coefficient $\hat\Psi_0 = (\hat{a}_a, \hat{\varphi}_y, \hat{a}_y, \hat{\varphi}_a)$, 
i.e.\ the coefficient of $s^0$, is well-defined. This coefficient lies in the eigenspace corresponding to the indicial root 
$\lambda = 0$, and hence, following the language at the end of section \ref{nonregular}, is an element of 
$\frak c^{8}$. We call this the Cauchy data of $\hat{\Psi}$ and write it as $\mathcal C(\hat{\Psi})$. 

At this ODE level, the Nahm pole boundary condition intermediates between the spaces $y^{\lambda_0 + 1/2}L^2$
when $-1/2 < \lambda_0 < 0$ and $0 < \lambda_0 < 1/2$. Recall that for this boundary condition, we
fix a subalgebra $\frak h \subset \frak c$ and its orthogonal complement $\frak h^\perp$ in $\frak c$, and
then consider the linear map 
\begin{equation} 
\mathcal B_{\frak h} (\hat{a}_a, \hat{\varphi}_y, \hat{a}_y, \hat{\varphi}_a) = (\hat{a}_a^{\, \frak h^\perp}, \hat{\varphi}_y^{\, \frak h^\perp} , 
\hat{a}_y^{\, \frak h} , \hat{\varphi}_a^{\,\frak h}) \in (\frak h^\perp)^3 \oplus \frak h^\perp \oplus \frak h \oplus \frak h^3. 
\label{bcs0}
\end{equation}
This determines a boundary condition for $\hat{N}(\LKW)$ and we shall study the problem
\begin{equation}
\hat{N}(\LKW) \hat \Psi = f, \qquad \hat\Psi \in s^{\lambda_0 + 1/2}H^1_0,\ \ \mathcal B_{\frak h}(\mathcal C(\hat \Psi)) = 0. 
\label{bcs}
\end{equation}

As a first observation, following the same integrations by parts as above, there are no nontrivial fields $\hat{\Psi}$ which decay 
at infinity and satisfy \eqref{bcs} with $f = 0$; indeed, the conditions $a_a^{\frak h^\perp} = 0$, $\varphi_y^{\frak h^\perp} = 0$, 
$a_y^{\frak h} = 0$, $\varphi_a^{\frak h} = 0$ make all boundary terms at $s=0$ in this integration by parts vanish. 
On the other hand, if we only assume that $f \in s^{\lambda_0 - 1/2}L^2$ for some $-1/2 < \lambda_0 < 0$, then this problem is 
not well posed. Indeed, although there is a solution to the first equation in \eqref{bcs}, there is no reason for the 
leading coefficient $\mathcal C(\hat \Psi)$ to have any meaning, so the boundary condition may not have any sense.
Thus it is necessary to suppose that $f$ lies in a slightly better space, as we now describe. 

With $\lambda_0 \in (-1/2,0)$ as before, define
\begin{equation}
\hat{\mathcal H}_{\lambda_0} = \{\hat \Psi \in s^{\lambda_0 + 1/2}L^2(\R^+): \hat{N}(\LKW) \hat \Psi \in s^{\lambda_0 + 1/2}L^2(\R^+)\}.
\label{Hnor}
\end{equation}
The right hand side is one order less singular than might be expected, and using standard ODE techniques, one sees
that $\hat \Psi = s^0 \hat \Psi_0 + \mathcal O(s^{\lambda_0 + 3/2})$, and hence the leading coefficient $\hat \Psi_0$ is 
well defined. We may now legitimately consider the mapping
\begin{equation}
\hat{N}(\LKW): \hat{\mathcal H}_{\lambda_0} \cap \{\hat \Psi: \mathcal B_{\frak h} (\mathcal C(\hat \Psi)) = 0\} 
\longrightarrow s^{\lambda_0 + 1/2} L^2. 
\end{equation}
Note that the weighted $L^2$ restriction on $\hat \Psi$ when $s$ is large precludes the exponentially growing solutions of 
$N(\LKW) \hat\Psi = 0$.

We have so far suppressed the dependence of $\hat{N}({\LKW})$, and hence also of the map $\mathcal B_{\frak h}$, on
the parameters $\vec x \in \partial M$ and $\vec k \neq 0$. As observed earlier, the only essential part of the dependence on
$\vec k$ is on the direction $\vec k/|\vec k|$, and hence we assume for the remainder of this discussion
that $|\vec k| = 1$, i.e.\ $\vec k \in S^2$. For this particular operator, the dependence on $\vec x$ is not 
very serious in that the operator `looks the same' in appropriate local coordinates at any point of $\del M$. 
However, formally we should be considering $(\vec x, \vec k)$ as a point
in $S^* \partial M$, the cosphere bundle of $\del M$. If $\pi: S^* \partial M \to \partial M$
is the natural projection, then $\mathcal B_{\frak h}$ is a bundle map between 
$\pi^* \frak c^{8}$ and $\pi^* \left((\frak h^\perp)^3 \oplus \frak h^\perp \oplus \frak h \oplus (\frak h)^3\right)$.
The fact that $\mathcal B_{\frak h}$ is independent of $\vec k$ allows us to call this an algebraic boundary condition. 

The additional ingredient we need in this discussion is the space 
\begin{equation}
V = V_{\vec x, \vec k} = \{ \hat\Psi \in \hat{\mathcal H}_{\lambda_0}:\hat N(\LKW)(\hat\Psi) = 0\}
\label{defV}
\end{equation}
of homogeneous solutions in $s^{\lambda_0 + 1/2}L^2$ which do not necessarily 
satisfy the boundary condition. The dependence
of $V$ on $\vec x$ is again negligible, but its dependence on $\vec k$ is genuine since $\vec k$ appears in the 
coefficients of $N(\LKW)$ and these kernel elements do vary nontrivially with $\vec k$.  However, the dimension
of $V_{\vec x, \vec k}$ does not depend on $\vec k$, and in fact $V$ varies smoothly with $\vec k$ (and $\vec x$).
This means that we can regard $V$ as a vector bundle over $S^* \partial M$. The injectivity of $\hat N({\LKW})$
on $y^{\lambda_0 + 1/2}L^2$ when $\lambda_0 > 0$ shows that the restriction of the Cauchy data map $\mathcal C$
to $V$ is injective. This means that $\mathcal C(V)$ is a subbundle of $\frak c^{8}$; this is sometimes
called the Calderon subbundle. 

We can now finally state the property which makes $\mathcal B_{\frak h}$ an elliptic boundary condition. 
\begin{definition}
Let $\mathcal B$ be any bundle map from $\pi^* \frak c^{8}$ to another vector bundle $W$ over $S^* \partial M$. 
Then $\mathcal B$ is said to be an elliptic boundary condition for $\hat{N}(\LKW)$ (and hence ultimately for $\LKW$) if
the restriction of $\mathcal B$ to the subbundle $V$ is bijective onto $W$, i.e.\ so that 
\begin{equation*}
\left. \mathcal B \right|_V: V  \longrightarrow W
\end{equation*} 
is a bundle isomorphism. 
\label{bcs1}
\end{definition}
This condition can be phrased in various obviously equivalent ways; the one we use below is to require 
that $\mathcal B$ is injective and that the ranks of $V$ and $W$ are the same. However, in the end, this
condition is precisely what is 
needed to construct a good parametrix for the actual boundary problem on $M$.  We can see this rather
easily at the level of this ODE.  If $f \in s^{\lambda_0 + 1/2}L^2 \subset s^{\lambda_0 - 1/2}L^2$, then there is
a solution $\hat \Psi \in s^{\lambda_0 + 1/2}H^1_0$ to $\hat{N}({\LKW})\hat\Psi = f$. This is not unique,
since there is a nullspace. In any case, this solution has a leading coefficient $\hat \Psi_0$, which
is however unlikely to satisfy $\mathcal B_{\frak h}(\hat \Psi_0 ) = 0$. Now modify $\hat\Psi$ by subtracting 
an element $\hat\Phi \in V$. By definition, $\hat{N}({\LKW})(\hat\Psi - \hat \Phi) = f$, and provided we choose 
$\hat \Phi$ so that $\mathcal B_{\frak h}(\hat\Phi_0) = \mathcal B_{\frak h}(\hat\Psi_0)$, 
then $\hat \Psi - \hat \Phi$ satisfies the boundary condition too. The fact that there is a unique 
such choice of $\hat \Phi$ is precisely the content of Definition~\ref{bcs1}.  

Let us now check that this condition holds for the map $\mathcal B_{\frak h}$ which appears in the general Nahm
pole boundary condition.  We have proved above that $\mathcal B_{\frak h}$ is injective on $V$.  Furthermore,
it follows from the results of section \ref{third} that the rank of $V$ is half the rank of $\frak c^{8}$, i.e.\ 
$\mathrm{rk}(V) = 4\dim \frak c$. Since this is the same as $\dim ( (\frak h^\perp)^3 \oplus \frak h^\perp \oplus 
\frak h \oplus \frak h^3)$, we see that $\left. \mathcal B_{\frak h}\right|_V$ is also surjective, and hence an isomorphism.

As noted earlier, $\mathcal B$ is called an algebraic boundary condition if $W = \pi^* W'$ where $W'$ is a bundle over 
$\del M$. If this is the case, then the analytic theory of the boundary problem for the actual operator $\LKW$ is simpler 
because the boundary conditions are local (of mixed Dirichlet-Neumann type), rather than nonlocal (pseudodifferential). 
It is clear that $\mathcal B_{\frak h}$ is an algebraic boundary condition. 

Return now to the linearized KW operator, and assume that some $j_\sigma = 0$, so that we can augment
the operator $\LKW$ with the boundary condition $\mathcal B_{\frak h}$. Fix $\lambda_0 \in (-1/2, 0)$ and consider 
\begin{equation}
\mathcal H_{\lambda_0} = \{\Psi \in y^{\lambda_0 + 1/2} H^1_0( \d\vec x \d y): {\LKW} \Psi \in y^{\lambda_0 + 1/2} L^2\}.
\end{equation} 
As before, the expected behavior for $\Psi \in y^{\lambda_0 + 1/2}L^2$ is that ${\LKW} \Psi \in y^{\lambda_0 - 1/2}L^2$.
This means that fields in $\mathcal H_{\lambda_0}$ must possess some special properties to ensure that ${\LKW}\Psi$ 
is one order less singular than $y^{\lambda_0 - 1/2}L^2$.  Although we no longer have ODE arguments to fall
back upon, it is still possible to show that any $\Psi \in \mathcal H_{\lambda_0}$ has a {\it weak} partial expansion
\begin{equation*}
\Psi \underset{w}{\sim} \Psi_0 \, y^0 + \tilde{\Psi},
\end{equation*}
where $\Psi_0$ is a distribution of negative order which lies in the Sobolev space $H^{\lambda_0}(\partial M)$. 
The remainder term $\tilde{\Psi}$ vanishes like $y^{\lambda_0 + 1}$ in a similar distributional sense. We do not 
pause to make this more precise (see section 7 of \cite{M-edge}), but note only that the actual meaning of the 
weak expansion above is that if we `test' $\Psi$ against some $\chi \in \mathcal C^\infty(\del M)$, then 
\begin{equation*}
\int \Psi(y,\vec x) \chi(\vec x) \, \d \vec x = \langle \Psi_0, \chi \rangle y^0 + \langle \tilde{\Psi}, \chi \rangle,
\end{equation*}
where the second term on the right vanishes like $y^{\lambda_0 + 3/2}$.  

The point of belaboring all of this is that it is possible to make sense of the leading coefficient $\Psi_0$ of 
a general element $\Psi \in \mathcal H$ as a $\frak c^{\oplus 8}$-valued distribution of negative
order on $\del M$.  Because the boundary condition is algebraic, we can then make sense of the projection 
$\mathcal B_{\frak h}(\Psi_0)$, again as a distribution. In particular, when $\Psi \in \mathcal H$, 
there is now a good meaning of the condition $\mathcal B_{\frak h}(\mathcal C(\Psi)) = 0$.

We can now state analogs of all the main results.
\begin{prop}\label{fred2}
Let ${\LKW}$ be the linearized KW operator on a compact manifold $M$ with boundary and suppose that $\rho$ is
not quasiregular, so that some $j_\sigma = 0$. Using the elliptic boundary condition given by the bundle map $\mathcal B_{\frak h}$,
then for $-1/2 < \lambda_0 < 0$, the mapping 
\begin{equation}
{\LKW}: \{\Psi \in \mathcal H_{\lambda_0}: \mathcal B_{\frak h}(\Psi_0) = 0\} \longrightarrow s^{\lambda_0 + 1/2} L^2(M)
\label{map12}
\end{equation}
is Fredholm. 
\end{prop}
\begin{prop}\label{regularity2}
With all notation as above, suppose that ${\LKW}\Psi = f$ where $f$ is smooth in a neighborhood of some point $q \in \del M$, 
or slightly more generally, where $f$ has an asymptotic expansion as $y \to 0$, and in addition $\mathcal B_{\frak h}(\Psi_0) = 0$
near $q$. Then in that neighborhood, $\Psi$ has a polyhomogeneous expansion 
\begin{equation}
\Psi \sim \sum y^{\gamma_j} \Psi_j,
\end{equation}
where the exponents $\gamma_j$ are of the form $\lambda + \ell$, $\ell \in \mathbb N$, where either $\lambda$ 
is an indicial root of ${\LKW}$ or else is an exponent occurring in the expansion of $f$. As before, this expansion may include log terms. 
\end{prop}
\begin{prop}
Let ${\LKW}$ be the linearized KW operator. If some $j_\sigma = 0$, then for $-1/2 < \lambda_0 < 0$ and relative 
to any choice of subalgebra $\frak h$, the mapping 
\begin{equation}
{\LKW}: y^{\lambda_0} \mathcal C^{k,\ell, \gamma}_0 \longrightarrow y^{\lambda_0 - 1} \mathcal C^{k-1, \ell, \gamma}_0
\label{Fredholder2}
\end{equation}
is Fredholm when $0 \leq \ell \leq k-1$ and $k \geq 1$. 
\label{maphold2}
\end{prop}

These results are proved, as before, using parametrix methods. Unlike the earlier case, this is a slightly
more involved process which requires the introduction of generalized Poisson and boundary trace operators;
see \cite{MVer}.   We find along the way that the analog of the refined mapping property of the generalized
inverse \eqref{refreghold2} still holds. 

\subsection{Regularity of Solutions of the KW Equations}\label{nonlinreg}

We have now described enough of the linear theory that we can formulate and prove the main result
needed earlier in the paper, that solutions of the full nonlinear gauge-fixed KW equations ${\KW}(A,\phi) = 0$
are also polyhomogeneous at $\partial M$. The implication of this regularity is that all the calculations
in section 2 which led to the uniqueness theorem when $M = \R^4_+$ are fully justified. (Of course,
this polyhomogeneity is far more than is really needed to carry out those calculations, but it is very
useful to have this sharp regularity for other purposes too.) 

\begin{prop} \label{propnlreg}
Let ${\KW}(A,\phi) = 0$, and suppose that near $y = 0$, $A = A_{(0)} + a$, $\phi = \phi_{(0)} + \varphi$,
where $(a,\varphi)$ satisfy the Nahm pole boundary conditions. Then $(a,\varphi)$ is polyhomogeneous.
\end{prop}

The proof begins by using \eqref{texpkw} to write ${\KW}(A, \phi) = 0$ as ${\LKW}(a,\varphi) = -Q(a,\varphi)$. 
Observe that since the terms in $\KW$ are at most quadratic, $Q$ is a bilinear form in $(a,\varphi)$.  We 
suppose from the beginning that $a$ and $\varphi$ lie in $y^{\lambda_0} \mathcal C^{1,\gamma}_0$, where the rate
of blowup (or decay) $\lambda_0$ is as dictated by the Nahm pole boundary condition.   The details 
of the proof are essentially the same in the simpler quasiregular case and in the more general 
case where some $j_\sigma = 0$. In  the former, $\overline{\lambda} =1$, $\underline \lambda=-1$, and we take
any elliptic weight $\lambda_0 \in 
(\underline{\lambda}, \overline{\lambda})$, while in the latter, generically we take $\lambda_0 \in (-1/2, 0)$ (if
$j_\sigma=1/2$ does not occur in the decomposition of $\frak g_\C$,  we can take $\lambda_0\in (-1,0)$),
and the generalized inverse is constructed using the more elaborate considerations of section \ref{genebc}.
The key facts that we use below, however, are the existence of a generalized inverse satisfying \eqref{geninv},
in particular so that the remainder term $R_1$ maps into a finite dimensional space of polyhomogeneous
fields and the fact that $G$ satisfies \eqref{refreghold2}. Although the construction of $G$ is more
complicated in the second case, we still end up with the result that both these properties hold then too.

There are two main steps. The first is to prove that $(a,\varphi)$ is conormal of order $\overline{\lambda}$, i.e. 
$(a,\varphi) \in \mathcal A^{\overline{\lambda}}$, which we recall means that 
\begin{equation}
(y\partial_y)^j \partial_{\vec x}^\alpha (a,\varphi) \in \bigcap_k y^{\overline{\lambda_0}} \mathcal C^{k,\gamma}_0
\label{conormal}
\end{equation}
for all $j$ and all multi-indices $\alpha$. In the second step we improve this to the existence of a polyhomogeneous expansion. 

Since $\lambda_0$ is an elliptic weight, there exists a generalized inverse $G$ for ${\LKW}$ which provides an 
inverse to \eqref{Fredholder} up to finite rank errors. Applying $G$ to ${\KW}(A,\phi) = 0$ gives
\begin{equation}
(a,\varphi) = - G Q(a, \varphi) + R_1 (a,\varphi).
\label{inteqn}
\end{equation}
The finite rank operator $R_1$ has range in the space of polyhomogeneous functions (with leading
term $y^{\overline{\lambda}}$), 
maps into a finite dimensional space of polyhomogeneous functions, so the second term on the right is 
polyhomogeneous and hence negligible. We are thus free to concentrate on proving the regularity of the first term
on the right in \eqref{inteqn}.

We first assert that $(a,\varphi) \in y^{\overline{\lambda}} \mathcal C^{k,\gamma}_0$ for all $k \geq 0$. (Note that
this is not conormality since we are not yet applying the tangential vector fields $\del_{x^a}$ without the
extra factor of $y$.)   Since $Q$ is bilinear, we see first that $Q(a,\varphi) \in y^{2\lambda_0} \mathcal C^{1,\gamma}_0$.
so that from \eqref{refreghold2} with $k = 1$, and since $2\lambda_0 > \lambda_0 - 1$, we obtain
$(a,\varphi) \in y^{2\lambda_0 + 1} \mathcal C^{2,\gamma}_0 + y^{\overline{\lambda}}\mathcal C^{m,\gamma}_0$ for all $m$. 
After a finite number of iterations, the right hand side is contained in $y^{\overline{\lambda}}\mathcal C^{m,\gamma}_0$
for some $m$, and bootstrapping further shows that it lies in this space for all $m$. 

Revisiting this iteration, we can improve the regularity with respect to the vector fields $\partial_{x^a}$ as well. 
This relies on a structural fact about $0$-pseudodifferential operators already quoted at the end of section \ref{strgeninv},
nmely that 
\begin{equation*}
[ \partial_{x^a} , G] \in \Psi^{-1, 1, \overline{\lambda}, b}_0,
\end{equation*}
and hence this commutator satisfies the same mapping properties as $G$ itself. We apply this as follows. Write
\begin{equation*}
\partial_{x^a} (a,\varphi) = - G( \partial_{x^a} Q(a,\varphi) ) - [ \partial_{x^\alpha} , G] Q(a,\varphi).
\end{equation*}
We are discarding the term $R_1(a,\varphi)$ since it is already fully regular. It is convenient now to
regard $(a,\varphi)$ as lying in $y^{\lambda_0}\mathcal C^{m,\gamma}_0$ for $\lambda_0 = \overline{\lambda} - \epsilon$,
since we want to use the mapping properties of $G$ at a nonindicial weight. By the mapping properties 
of the commutator, the second term on the right lies in $\cap_m y^{\lambda_0} \mathcal C^{m,\gamma}_0$.
On the other hand, we write $\partial_{x^a} Q(a,\varphi) = y\partial_{x^a} (y^{-1} Q(a,\varphi))$. This lies
in $\cap_m y^{2\lambda_0 - 1} \mathcal C^{m,\gamma}_0$, since $y^{-1}Q(a,\varphi) \in \cap_m y^{2\lambda_0 - 1} \mathcal C^{m,\gamma}_0$
and $y\del_{x^a}$ preserves this property.  Since $G$ acts on this space, the entire first term lies in $\cap_m y^{\lambda_0} 
\mathcal C^{m,\gamma}_0$. This proves that $(a,\varphi) \in \cap_m \mathcal C^{m,1, \gamma}_0$. The same argument improves the tangential 
regularity incrementally, so $(a,\varphi) \in y^{\lambda_0} \mathcal C^{k,\ell,\gamma}_0$ for all $0 \leq \ell \leq k < \infty$.
Recalling \eqref{refreghold2} again, we can now replace $\lambda_0$ by $\overline{\lambda}$. This proves that
$(a,\varphi) \in \mathcal A^{\overline{\lambda}}$. 

The second main step of the proof is easier. We wish to prove that $(a,\varphi)$ has an expansion. This relies on
the observation that we can treat the nonlinear equation ${\KW}(A_{(0)}+a, \phi_{(0)}+\varphi) = 0$ as a nonlinear 
ODE in $y$, regarding the dependence on $\vec x$ as parametric. This is reasonable since $(a,\varphi)$ 
is completely smooth in this tangential variable. Thus we can decompose the linear term in \eqref{texpkw} 
further using the indicial operator $I({\LKW})$ (introduced in eqn.\ (\ref{indicialop})) at any given boundary point to get
\begin{equation}
I({\LKW}) (a, \varphi) = (I({\LKW})-\LKW) (a,\varphi) - Q(a,\varphi).
\label{iter}
\end{equation}
The two terms on the right lie in $y^{\lambda_0 } \mathcal A$ and $y^{2\lambda_0} \mathcal A$, respectively. 
(We recall that $I({\LKW}) - {\LKW}$ has no $\partial_y$ or $1/y$ terms; the coefficients of this
difference are smooth to $y=0$.)  Integrating this ODE shows that $(a,\varphi)$ is a finite sum of terms 
$(a_j, \varphi_j) y^{\lambda_j}$, where the $\lambda_j$ are the indicial roots of ${\LKW}$ which lie between $\lambda_0$ 
and $\mu = \min\{2\lambda_0+1, \lambda_0 + 1\}$, and an error term vanishing at this faster rate $y^\mu$. 
At the next step, inserting this new information into \eqref{iter} shows that this is now an ODE where the right 
side has a partial expansion up to order $\min\{\mu, 2\mu\}$ plus an error term vanishing at that rate,
and so $(a,\varphi)$ has a partial expansion up to order $\mu_2 = \min \{\mu + 1, 2\mu + 1\}$. 
This completes the proof of the existence of the expansion of $(a,\varphi)$ in the case
where $\lambda_0 \in (0, \overline{\lambda})$ is an elliptic weight. 

\appendix
\section{Some Group Theory}\label{groups}
\def\frak{\mathfrak}

The purpose of this appendix is to describe some basic facts and examples in group theory as background to the paper.  

First of all, up to isomorphism, the group $SU(2)$ or equivalently the Lie algebra $\frak{su}(2)$ has precisely
one irreducible representation of dimension $n$, for each positive integer $n$.  
It is convenient to write $n=2j+1$ where $j$ is a non-negative integer or half-integer called the spin.
If $v_j$ denotes an irreducible $\frak{su}(2)$ representation of spin $j$, then for $j\geq j'$, we have
\begin{equation}\label{decomp}v_j\otimes v_{j'}\cong \oplus_{j''=j-j'}^{j+j'}v_{j''}. \end{equation}

Now, for $N\geq 2$,  we will describe group 
homomorphisms $\varrho:SU(2)\to G=SU(N)$, or equivalently Lie algebra homomorphisms $\varrho:\frak{su}(2)\to \frak \frak{su}(N)$.
To describe such a homomorphism amounts to describing how the fundamental $N$-dimensional representation of $SU(N)$,
which we denote $V$, transforms under $\varrho(SU(2))$.  As an $SU(2)$-module, $V$ 
will have to be the direct sum of a number of irreducible $SU(2)$ modules $v_{j_i}$ of dimension $n_i=2j_i+1$, for some $j_i$.
The possibilities simply correspond to partitions of $N$, that is to ways of writing
$N$ as an (unordered) sum of positive integers,
\begin{equation}\label{ecomp}N=n_1+n_2+\dots+n_s. \end{equation}
For example, the trivial homomorphism $\varrho:\frak{su}(2)\to\frak{su}(N)$, which maps $\frak{su}(2)$ to 0,
 corresponds to the partition $N=1+1+\dots+1$ with $N$ terms. 
 At the other extreme, a principal embedding of $\frak{su}(2) $ in $\frak{su}(N)$
(which is the most important example for the present paper) corresponds to the partition with only one term, the integer $N$.

In general, we define the commutant $C$ of $\varrho$ as the subgroup of $SU(N)$ that commutes with $\varrho(SU(2))$;
its Lie algebra $\frak c$ is the subalgebra of $\frak{su}(n)$ that commutes with $\varrho(\frak{su}(2))$.  
If $\varrho$ corresponds as in eqn.\ (\ref{ecomp}) to a partition with $s$ terms, then $C$ is a Lie group of rank $s-1$.
(It is abelian if and only if the $n_i$ are all distinct.)  In particular, for $G=SU(N)$, the only case that $C$ is a finite group (or equivalently
a group of rank 0) is
that $s=1$, meaning that $\varrho$ is a principal embedding.  In this case, $C$ is simply the center of $G$.  Whenever $s>1$,
$C$ has a non-trivial Lie algebra, and this means, in the language of section \ref{second}, that $j_\sigma=0$ occurs
in the decomposition of $\frak g_\C$ under $\frak{su}(2)$.  Thus, for $G=SU(N)$, the only case that $j_\sigma=0$ does
not occur in this decomposition is the  case that $\varrho$ is a principal embedding. 

To explicitly decompose $\frak{su}(N)$ under $\varrho(\frak{su}(2))$, we use the fact that $\frak{su}(N)$ is the traceless
part of $\mathrm{Hom}(V,V)$; equivalently it can be obtained from $V\otimes V^\vee$ by omitting a 1-dimensional
trivial representation.  (Here $V^\vee$ is the dual of $V$.)  Any $\frak{su}(2)$-module is isomorphic to its own dual,
so as a $\frak{su}(2)$ module, $\frak{su}(N)$ is $V\otimes V$ with a copy of the trivial module $v_0$ removed.  
For example, if $\varrho$ is a principal embedding, so that $V$ is an irreducible $\varrho(\frak{su}(2))$ 
module $v_j$ with $N=2j+1$, then we use (\ref{decomp}) to learn that $\frak{su}(N)$ is the direct sum of
$\frak{su}(2)$-modules of spins $j_\sigma=1,2,\dots, N-1$.  

As a corollary, we note that if $j_\sigma=0$ does not occur in the decomposition of $\frak{su}(n)$ (which happens only
if $\varrho$ is principal), then the $j_\sigma$'s are integers and in particular $j_\sigma=1/2$ does not occur in the decomposition
of $\frak{su}(n)$.  As explained below, this statement has an analog for any simple Lie group $G$.
A few additional facts that follow from the above discussion of $SU(N)$ hold for all $G$.
The number of summands in the decomposition of $\frak g$ under a principal $\frak{su}(2)$ subalgebra 
 is always the rank of $G$ (this rank is $N-1$ for $G=SU(N)$).  Also, the minimum
value of $j_\sigma$ for a principal embedding is always
 $j_\sigma=1$, and this value occurs with multiplicity 1,  corresponding to the $\frak{su}(2)$ submodule
$\varrho(\frak{su}(2))\subset \frak g$.

With similar elementary methods, we can analyze the other classical groups $SO(N)$ and $Sp(2k)$. 
Here the following is useful.  An $SU(2)$ module $v$ is said to be real, or to admit a real structure, if there
is a symmetric, non-degenerate, and $SU(2)$-invariant map $v\otimes v\to \C$; it is said to be pseudoreal, or to admit
a pseudoreal structure, if there is an antisymmetric, non-degenerate, and $SU(2)$-invariant map $v\otimes v\to \C$.
The representation $v_j$ is real (but not pseudoreal) if $j$ is an integer, or equivalently the dimension
 $n=2j+1$ is odd, and pseudoreal  (but not real) if $j$ is a half-integer, or equivalently the dimension $n=2j+1$ is even. 
 If $w$ is a 2-dimensional complex vector space (with trivial $SU(2)$ action), then $w$ admits both a symmetric nondegenerate
map $w\otimes w\to \C$ and an antisymmetric one.  So if $v$ is either real or pseudoreal, then $v\oplus v\cong v\otimes w$
admits both a real structure and a pseudreal one.
Suppose that
\begin{equation}\label{morex}v=\oplus_{j\geq 0}a_j v_j,~~~   a_j\in \Bbb Z\end{equation}
is an $SU(2)$ module that is the direct sum of $a_j$ copies of $v_j$ (with almost all $a_j$ vanishing).  
The criterion
for $v$ to be real or pseudoreal reduces to separate conditions on each $a_j$:
$v$ is real precisely if $a_j$ is even for half-integer $j$ (with no restriction for integer $j$), and $v$ is pseudoreal precisely if
$a_j$ is even for integer $j$ (with no restriction for half-integer $j$).  

Now let us consider homomorphisms $\varrho:SU(2)\to G$ for $G=SO(N)$.
Such a homomorphism can be described by giving the decomposition of the fundamental
$N$-dimensional representation $V$ of $SO(N)$ as an $SU(2)$-module.  Thus, such a homomorphism
again determines a partition of the integer $N$, as in (\ref{ecomp}).  However, now we must impose the condition that
the representation $V$ of $SO(N)$ is real.  In view of the statements in the last paragraph,
the condition that this imposes on the partition is simply that even integers $n_i$ in (\ref{ecomp}) must occur with even
multiplicity.

The condition that the commutant $C$ of $SU(2)$ -- or more precisely of $\varrho(SU(2))$ --
is a finite group, and hence that $j_\sigma=0$ does not
occur in the decomposition of the Lie algebra $\frak{so}(N)$, is\footnote{If an  integer $n_i$ appears with multiplicity $d_i$ in the partition of $N$,
then the commutant of $\frak{su}(2)$ in $SO(N)$ contains a factor of $SO(d_i)$ if $n_i$ is odd or $Sp(d_i)$ if $n_i$ is even.  (The last statement
makes sense because $d_i$ is always even when $n_i$ is even.)  
Hence the group $C$ is finite if and only if the $n_i$ are all distinct, so that the $d_i$ are 0 or 1.
The last statement is true for $G=Sp(2k)$ for similar reasons: if an integer $n_i$ appears with multiplicity $d_i$ in the partition of $2k$, then the commutant contains a factor of
$SO(d_i)$ if $n_i$ is even and of $Sp(d_i)$ if $n_i$ is odd.  (For $G=Sp(2k)$,  $d_i$ is even when $n_i$ is odd.)}  that the integers $n_i$ in the partition must be all distinct.
But since even integers must occur with even multiplicity, this   implies that the $n_i$ must be odd.  For example, if $N$ is odd,
a principal embedding of $\frak{su}(2)$ in $\frak{so}(N)$ corresponds to a partition with only one term, the integer $N$.
But if $N$ is even, a principal embedding corresponds to the two-term partition $N=1+(N-1)$.  For $SO(N)$, in contrast to $SU(N)$,
 an embedding that is not principal can still have a trivial commutant.  For example, for $N=9$,
the partition $9=1+3+5$ represents 9 as the sum of distinct odd integers; this embedding is not principal, but its commutant
 is a finite group.  When the commutant is not a finite group,  it has a Lie algebra of positive dimension and hence $j_\sigma=0$ occurs in the
decomposition of $\frak{so}(n)$ under $\frak{su}(2)$.

For $G=SO(N)$, rather as we found for $SU(N)$,
 there is also a useful elementary criterion that ensures that $j_\sigma=1/2$ does not appear in the decomposition
of $\frak{so}(n)$. In fact, there is a useful criterion that ensures that no half-integer value of $j_\sigma$ occurs
in this decomposition.  For this, we first recall that $\frak{so}(N)\cong \wedge^2 V\subset V\otimes V$.  
For a given $\frak{su}(2)$ embedding, the decomposition of $V\otimes V$ under $\frak{su}(2)$ can be worked out using (\ref{decomp}).
One finds that half-integer values of $j$ occur in $V\otimes V$ (and also in $\wedge^2V$) if and only if both
odd and even integers $n_i$ occur in the chosen partition  of $N$.  But we have already observed that if even integers
appear in this partition, then $j_\sigma=0$ occurs in the decomposition of $\frak{so}(n)$ under $\frak{su}(2)$.
So if $j_\sigma=0$ does not occur in the decomposition of $\frak{so}(n)$, then $j_\sigma=1/2$ also does not occur.

The case that $G=Sp(2k)$ for some $k$ can be analyzed similarly, with the words ``even'' and ``odd'' exchanged in some statements.
A homomorphism from $SU(2)$ to $Sp(2k)$ can be described by giving the decomposition of the $2k$-dimensional representation
$V$ of $Sp(2k)$ as a direct sum of $SU(2)$ modules.  Thus, such a homomorphism determines a partition $2k=n_1+n_2+\dots+n_s$.
Now the fact the representation $V$ of $Sp(2k)$ is pseudoreal implies that odd integers occur in this partition with even multiplicity.
The condition that the commutant $C$ is a finite group, so that $j_\sigma=0$ does not occur in the decomposition
of $\frak{sp}(2k)$ under $\frak{su}(2)$,  is again that the integers appearing in the partition should be distinct.
But now this implies that these integers are all even.  A principal embedding is the case that the partition consists of only of a single
integer $2k$.   Just as for $SO(N)$, there are non-principal embeddings with the property that $j_\sigma=0$ does not occur
in the decomposition of $\frak{sp}(2n)$; these correspond to partitions of $2k$ as the sum of distinct even integers, for example
$6=2+4$.

By the same argument as for $SO(N)$, one can show that if $j_\sigma=0$ does not occur in the decomposition of $\frak{sp}(2k)$,
then half-integer values of $j_\sigma$ do not occur and in particular $j_\sigma=1/2$ does not occur.  For this, one uses
the fact that $\frak{sp}(2k)\cong {\mathrm{Sym}}^2V\subset V\otimes V$, along with the rule (\ref{decomp}) for decomposition of tensor
products.

To understand homomorphisms from $SU(2)$ to an exceptional Lie group $G$, it is probably best to use less elementary
methods, and we will not explore this here.  We remark, however,  that the following feature of the above examples
is actually true for any simple Lie group $G$:  if $j_\sigma=0$ does not occur in the decomposition of $\frak g_\C$ under $\frak{su}(2)$,
then only integer values of $j_\sigma$ occur in this decomposition and in particular $j_\sigma=1/2$ does not occur.  (For a proof,
see the next paragraph.)
Given this, it follows from the formulas of section \ref{throots} that if $j_\sigma=0$ does not occur in the decomposition
of  $\frak g_\C$, then there are no indicial roots in the gap between $\underline\lambda=-1$ and $\overline\lambda=1$.
Both $-1$ and $1$ always are indicial roots in this situation, 
since $j_\sigma=1$ always occurs in the decomposition of $\frak g_\C$, the corresponding subspace of $\frak g_\C$ being
$\varrho(\frak{su}(2))$.

A proof of a claim in the last paragraph was sketched for us by B. Kostant.  The complexification of  $\varrho$ is a homomorphism
of complex Lie algebras $\varrho:\frak{sl}_2(\C)\to \frak g_\C$.  We take a standard basis $(h,e,f)$ of $\frak{sl}_2(\C)$ and
write simply $(h,e,f)$ for their images in $\frak g_\C$.    The hypothesis that half-integer values of $j_\sigma$ occur in the decomposition
of $\frak g_\C$ means, in the terminology of \cite{Carter}, p. 165, that $e$ is not even.
 In this case, according to Proposition 5.7.6 of that
reference, $e$ is not distinguished, and therefore, according to Proposition 5.7.4 of the same reference, the homomorphism $\mathrm{ad}(e):
\frak g(0)\to \frak g(2)$ has a non-trivial kernel.  This kernel is the commutant $\frak c$ of $\varrho(\frak{sl}_2(\C))$.
\vskip1cm
 \noindent {\it {Acknowledgements}}  Research of RM supported in part by NSF Grant DMS-1105050.
 Research of EW supported in part by NSF Grant PHY-1314311.
 
  \bibliographystyle{unsrt}

\end{document}